%% file: DRW_130121_Marc210215_depotjournal.tex
\documentclass[a4paper,11pt,reqno]{amsart}
 
\usepackage{amsmath}
\usepackage{amssymb}
\usepackage{xspace,epsfig}
\usepackage{mathptmx}
 
\usepackage{graphicx}
\usepackage{a4wide}
\usepackage{eucal}
\usepackage[pdftex,colorlinks]{hyperref}

\makeatletter
\let\NAT@parse\undefined
\makeatother

\usepackage{pgfplots}
\pgfplotsset{compat=newest}
\usetikzlibrary{plotmarks}
\usetikzlibrary{arrows.meta}
\usepgfplotslibrary{patchplots}
\usepackage{grffile}

\usepackage{graphicx}
\usepackage{tikz}
\usepackage{bm}
\usepackage{cuted} 
\usepackage{datetime}
\usepackage{hyperref}

\usepackage{amsthm}

\numberwithin{equation}{section}

\newtheorem{theorem}{Theorem}[section]
\newtheorem{proposition}[theorem]{Proposition}
\newtheorem{lemma}[theorem]{Lemma}

\newtheorem{remark}{Remark}
\theoremstyle{remark}

\newtheorem{definition}{Definition}
 
\usepackage{enumitem}

 \newcommand{\myeq}[2]{\begin{equation}#1\label{eq:#2}\end{equation}}
\newcommand{\mymatrix}[1]{\begin{bmatrix}#1\end{bmatrix}}
\newcommand{\lemmaspace}{\vspace{.5cm}}
\newcommand{\lemmaspacetwo}{\vspace{.01cm}}

\newcommand{\plotkernelwidth}{3.5in}
\newcommand{\plotkernelheight}{3in}
\newcommand{\plottimewidth}{4in}
\newcommand{\plottimeheight}{3.5in}

\newcommand{\halpha}{{\nu_2}}
\newcommand{\hbeta}{{\nu_1}}
\renewcommand{\psi}{\halpha}
\renewcommand{\zeta}{\hbeta}
\newcommand{\mathrmd}{\mathrm{d}} 
 
\newcommand{\R}{\mathbb R}

\newcommand{\be}{\begin{equation}}
\newcommand{\ee}{\end{equation}}
\newcommand{\ba}{\begin{eqnarray}}
\newcommand{\ea}{\end{eqnarray}}

\usepackage{hyperref}

\numberwithin{equation}{section}

\begin{document}
\title[Finite-time stabilization of an overhead crane with a flexible cable]{Finite-time stabilization of an overhead crane with a \\ flexible cable submitted to an affine tension}
 
\author[Wijnand]{Marc Wijnand}
\address{Sorbonne Universit\'e, UMR 9912 STMS (Ircam -- CNRS -- Sorbonne Universit\'e), 
1 place Igor  Stravinsky, 
75004 Paris, France}
\email{marc.wijnand@ircam.fr}

\author[d'Andr\'ea-Novel]{Brigitte d'Andr\'ea-Novel}
\address{Sorbonne Universit\'e, UMR 9912 STMS (Ircam -- CNRS -- Sorbonne Universit\'e), 
1 place Igor  Stravinsky, 
75004 Paris, France}
\email{brigitte.dandreanovel@ircam.fr}

\author[Rosier]{Lionel Rosier}
\address{Universit\'e du Littoral C\^ote d'Opale, Laboratoire de Math\'ematiques Pures et Appliqu\'ees J. Liouville, 50 Rue Ferdinand Buisson, 62100 Calais, France}
\email{lionel.rosier@univ-littoral.fr}

\date{}
 
\maketitle
 
\begin{abstract}
The paper is concerned with the finite-time stabilization of a hybrid PDE-ODE system describing the motion of an overhead crane with a flexible cable. 
The dynamics of the flexible cable is described by the wave equation with a variable coefficient which is 
an affine function of  the curvilinear abscissa along the cable. Using several changes of variables, 
a backstepping transformation, and 
a finite-time stable second-order ODE for the dynamics of a  conveniently chosen variable,  
we prove that a global finite-time stabilization occurs for the full system constituted of the platform and the cable. 
The kernel equations and the finite-time stable ODE are numerically solved in order to compute 
the nonlinear feedback law, and
numerical simulations validating our finite-time stabilization approach are presented.
\end{abstract}
 
\vspace{0.3cm}
 
\textbf{2010 Mathematics Subject Classification:} 93C20, 93D15.
 
\vspace{0.5cm}
 
\textbf{Keywords:} Finite-time stability; PDE-ODE system; Volterra integral transformation; cascade systems;
nonlinear feedback law; transparent boundary conditions; wave equation.
 

 \section{INTRODUCTION}

\subsection{Stabilization of hybrid PDE-ODE systems}

The stabilization of hybrid PDE-ODE systems has attracted the attention of the control community since several decades.
We can mention applications such as the control of a rotating body-beam without natural damping \cite{BBeam} as well as the case of a rotating and translating body-beam \cite{mattioni2020stabilisation}, a slide-flute \cite{dandrea2010acoustic}, a switched power converter with a transmission line \cite{daafouz2014nonlinear}, turbulent fluid motion and traffic flow \cite{trinh2017design}, and the overhead crane taking into account the flexibility of the cable, as discussed below.

\subsection{Crane model}

In \cite{ABCR}, the authors derived and investigated a model for the dynamics of a motorized platform  of mass $M$ moving along a horizontal bench. A flexible (and nonstretching) cable of length $\bm{1}$ was attached to the platform and was holding a load mass $m$. Assuming that the transversal and angular displacements were small and that the acceleration of the load mass could be neglected, they obtained the following system:
\begin{align}
y_{tt}- (d(s)y_{s})_s&=0, \label{AA1} \\
y_s(0, t)&=0, \label{AA2} \\
y(1, t) &= X_\mathrm{p}(t),\label{AA3}  \\ 
\ddot X_\mathrm{p}(t) &= \lambda (dy_s)(1, t) + \frac{V(t)}{M}, \label{AA4}
\end{align}
where
\begin{align}
d(s)        &:=  g s+ \frac{gm} { \rho},  \label{AAAA1}\\
\lambda &:= \frac{ (m+\rho)g } {Md(1) } \cdot  \label{AAAA2}
\end{align}
The initial conditions are $y(s,0) = y^0(s)$, $y_t(s,0) = y^1(s), X_\mathrm{p}(0)=X_\mathrm{p}^0, \dot{X}_\mathrm{p}(0) = X_\mathrm{p}^1$.
In the above system, $s$ denotes the curvilinear abscissa (i.e. the arclength) along the cable, 
$y=y(s,t)$ is the horizontal displacement at time $t$ of the point on the cable of 
curvilinear abscissa $s$, $X_\mathrm{p}$ is the abscissa of the platform, $\rho$ the mass per unit length
of the cable, and $V$ the force applied to the platform.  As usual, 
$y_{tt}=\partial ^2 y /\partial t^2$, $y_{ss}=\partial ^2 y /\partial s^2$, etc., 
and $\ddot X_\mathrm{p}=\mathrm{d}^2X_\mathrm{p}/\mathrm{d}t^2$.

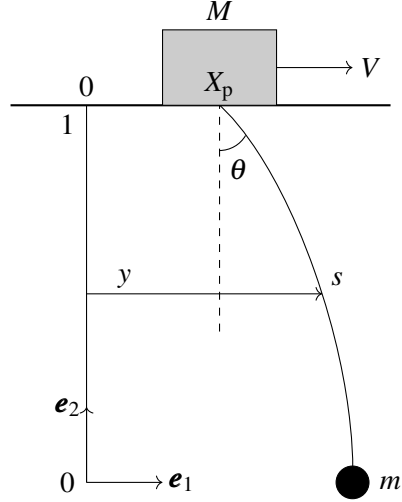
\begin{figure}[htpb]
   \begin{center}
   \begin{tikzpicture}
       \node[left of=, xshift=.75cm] at (0,0) {$0$};
       \node[right of=, xshift=-.75cm] at (1,0) {$\bm{e}_1$};
       \node[left of=, xshift=.75cm] at (0,1) {$\bm{e}_2$};
       \node[left of=, xshift=.75cm, yshift=-.25cm] at (0,5) {$1$};
       \node[above of=, yshift=-.75cm] at (0,5) {$0$};
	   \draw[->] (0,0) -- (1, 0);
	   \draw[->] (0, 0) -- (0, 1);
	   \draw (0,0) -- (0,5);
	   \draw[thick=3] (-1,5) -- (4,5);
	   \draw[fill=black!20!white] (1,5) -- (2.5,5) -- (2.5,6) -- (1,6)-- (1,5);
	   \node[above of=] at (1.75,5.25) {$M$};
	   \node[above of=] at (1.75, 4.25) {$X_\mathrm{p}$};
	   \draw[->] (2.5,5.5) -- (3.5,5.5);
	   \node[right of=, xshift=-.75cm] at (3.5,5.5) {$V$};
	   \draw [black, fill=black] (3.5,0) circle (1.25ex);
	   \node[right of=] at (3, 0) {$m$};
	   \draw[dashed] (1.75,2) -- (1.75,5);
	   \draw plot [smooth, tension=1.2] coordinates {(1.75,5) (3, 2.8) (3.5, 0)};
	   \draw[->] (0,2.5) -- (3.1,2.5);
	   \node[] at (.5,2.7) {$y$};
	   \node at (3.3,2.7) {$s$};
	   \draw (1.4, 4.6) ++(-30:.4) arc (-90:-30:.4);
	   \node at (2, 4.15) {$\theta$};
   \end{tikzpicture}
   \end{center}
   \caption{The overhead crane with flexible cable}
   \label{fig:grue}
\end{figure}

In \cite{DMR}, the authors supposed that $m\gg \rho$, so that $ gm / \rho \gg gs$ for $s\in (0, 1)$ and  it could be assumed that the function $d=d(s)$ was constant. In the present paper, we go back to the original problem without this assumption, so that the tension $d(s)$ is given by its affine expression \eqref{AAAA1}. After the following intermediate feedback law:
\be
\label{eqfeedbackuv}
V (t) = M U(t) - (m+\rho) g \theta(t) , \; \mbox{with} \; \theta(t) := y_s (1, t),
\ee
where the angular deviation $\theta(t)$ of the cable with respect to the vertical axis, at the curvilinear abscissa $s=1$ (i.e. at the connection point to the platform), is supposed to be measured (see Fig. \ref{fig:grue}), we obtain the following system:
\ba
y_{tt}-(d(s)y_{s})_s =0,&&  (s,t)\in (0, 1) \times (0, +\infty), \label{A1}\\
y_s(0, t) =0,&& t\in (0, +\infty), \label{A2}\\
y(1, t) = X_\mathrm{p}(t),&& t\in (0, +\infty), \label{A3}\\ 
\ddot X_\mathrm{p}(t) = U(t),&& t\in (0, +\infty). \label{A4}
\ea
with initial conditions $y(s,0) = y^0(s)$, $y_t(s,0) = y^1(s), X_\mathrm{p}(0)=X_\mathrm{p}^0, \dot{X}_\mathrm{p}(0) = X_\mathrm{p}^1$.

\subsection{Previously obtained control results for the crane}

An asymptotic (but not exponential) stabilization of \eqref{AA1}-\eqref{AA4} was established in \cite{ABCR}, while
an exponential stabilization was subsequently derived in \cite{AC} by using the cascaded structure of the system and a backstepping approach.
A similar result was obtained for system \eqref{A1}-\eqref{A4} with a constant tension (but with Dirichlet boundary conditions) in \cite{MRC}.
The dynamics of the load mass was taken into account in \cite {mifdal}.

The backstepping approach is a powerful tool for the design of stabilizing controllers in the context of finite dimensional systems 
(see for example \cite{Sepulchre}), but the cascaded structure of flexible mechanical systems coupling ODE and PDE is also a useful property in regard to stabilization, as for the overhead crane with flexible cable. One can also refer to \cite{BBeam}, where the authors proposed a class of nonlinear asymptotically stabilizing boundary feedback laws for a rotating body-beam without natural damping.
Let us also mention that in \cite{VLFC} the authors considered the case of a variable length flexible cable.

In \cite{AC}, if we restrict ourselves to system \eqref{A1}-\eqref{A4} with a constant tension (rescaled to be $1$), the authors considered the linear feedback law 
\[
\label{A5}
U(t) = -K^{-1} \left( y_{st} (1, t) +k y_t(1, t) \right) \\ -\mu \left( \dot X_\mathrm{p}(t) +K^{-1} \left( y_s(1, t) + ky(1, t)\right) \right),
\]
where $k,K>0$ and $\mu > K/2$ were some constants, and proved that system \eqref{A1}-\eqref{A4} 
was exponentially stable.

In \cite{DMR}, the same system with constant tension was stabilized in finite time by the nonlinear feedback law
\be
U(t) = -y_{st}(1, t) - \left\lfloor y_t(1, t) + y_s(1, t) \right\rceil^\halpha \\ - \left\lfloor y(1, t) + \int_0^1 y_t(\xi,t)\, \mathrm{d} \xi \right\rceil^\hbeta,
\ee
where $0<\halpha<1$, $\hbeta>\frac{\halpha}{2-\halpha}$ 
and $\lfloor x\rceil^\hbeta := \mathrm{sign}(x)|x|^\hbeta$. (Actually, the finite-time stability was fully justified in \cite{DMR} when $0<\halpha<1$ and $\hbeta=1$, although numerical simulations suggested that it was valid for $0<\halpha<1$ and
$\hbeta>\frac{\halpha}{2-\halpha} \cdot$)

\subsection{Finite-time stabilization of $2\times 2$ hyperbolic systems}
Solutions of certain hyperbolic systems with {\em transparent  boundary conditions} can reach the  equilibrium state in {\em finite time}. Such a property, 
called {\em finite-time stability} in \cite{APR,PR} or {\em super-stability} in \cite{SLX}, was first noticed in \cite{MP} for the usual wave equation. The extension of such a 
property to the wave equation on a network was investigated in \cite{APR, SLX}. 
It turns out that the finite-time stability also holds for one-dimensional first order quasilinear hyperbolic systems of diagonal form without source terms
\cite{LS,PR}. If source terms are incorporated in such systems, the finite-time stability is in general lost when using transparent boundary conditions, 
but a {\em rapid} stabilization still occurs \cite{GPR}. In the linear case, however, the use of a boundary feedback law based on a backstepping transformation  allows to recover the finite-time stability for hyperbolic systems with source terms \cite{CVKB}.   

\subsection{Finite-stability of an abstract evolution system}
Let $\Phi =\Phi (x, t) $, $x\in H$, $t\in \R _+$, be the flow associated with an evolution system in a Hilbert space $H$. We assume 
 that $\Phi$ satisfies the 
semigroup property: $\Phi (\Phi(x,s),t)=\Phi(x,t+s)$ for all $x\in H$ and all $t, s\in \R_+$, and that $x=0$ is an equilibrium point; that is, $\Phi (0,t)=0$ for all $t\in \R _+$. 

We say that the flow $\big( \Phi(x, t)\big) _{(x, t)\in H \times \R_+}$ is {\em globally finite-time stable} if there exists a nondecreasing function $T  : (0, +\infty) \to  (0, + \infty)$
such that $\Phi (x, t)=0$ for all $x\in H \setminus \{ 0 \}$ and all $t\ge T(\Vert x\Vert _H)$,  and if the equilibrium point $x=0$ is  Lyapunov stable; that is, for each $\varepsilon >0$, there  exists some $\delta >0$ such that 
\[
\Vert x\Vert_H < \delta \Rightarrow \Vert \Phi(x, t)\Vert_H <\epsilon \quad \forall t\ge 0.  
\] 
The finite-time stability of the second-order ODE $\ddot x= -\lfloor \dot x\rceil^\psi  -\lfloor x\rceil^\zeta$, 
or equivalently of the first-order system 
\ba
\dot x_1&=& x_2, \\
\dot x_2&=&  -\lfloor x_2\rceil^\psi  -\lfloor x_1\rceil^\zeta ,
\ea 
was established for $0<\psi <1$ and $\zeta > \frac{\psi}{2-\psi}$  
by Haimo (see \cite{haimo}), and for  $0<\psi < 1$ and  $\zeta =\frac{\psi}{2-\psi}$ by 
Bhat and Bernstein (see \cite{BR,BB}).
      
\subsection{Aim and structure of the present paper}

The aim of the paper is to design a boundary feedback law $U(t)$ leading to the finite-time stability of the system with affine tension.
For that purpose, we first use several changes of variables to transform the original system 
\eqref{A1}-\eqref{A4} into a 2$\times$2 hyperbolic system with coupling terms.
Next, following  \cite{CVKB}, we define a target system for which the application of transparent 
boundary conditions gives a finite-time stability, and we define a backstepping transformation leading to this target system.  
Finally, we show that the finite-time stabilization of a certain quantity $\phi$ defined in terms of  $X_\mathrm{p}(t)$ and $z=z(x, t)$
(see below \eqref{eq:defphi})
yields the finite-time stabilization of both the platform and the cable.

The paper is scheduled as follows. 
In Section \ref{sec:transforms}, the original system is transformed into a 2$\times$2 hyperbolic system with coupling terms. In Section \ref{sec:back}, the 2$\times$2 hyperbolic system is transformed into the target system by using  some Volterra integral transformation borrowed from \cite{CVKB}.  Methods for numerically solving the kernel equations are described. Section \ref{sec:res} is dedicated to the proof of our main result concerning the finite-time stabilization of the hybrid PDE-ODE system. Our theoretical result is illustrated by a simulation of the finite-time stabilization of the overhead crane in Section  \ref{sec:simu}.
Finally, we give some words of conclusion in Section \ref{sec:confut}.


\section{Derivation of the $2\times 2$  hyperbolic system} \label{sec:transforms}
The second order hyperbolic equation \eqref{A1} with boundary conditions \eqref{A2}-\eqref{A3} and initial conditions is rewritten as a $2\times2$ quasilinear hyperbolic system by performing two transformations.
\subsection{First transformation}
Following \cite{MRR}, we set  $z(x, t) := y(s, t)$, where
\myeq{ x(s): = \frac{1}{J} \int_0^s d^{-1}(\sigma) \, \mathrm{d} \sigma = \frac{1}{gJ} \ln \left( 1 + \frac{\rho}{m} s \right) \ge 0 }{tfyz}
and
\[ J := \int_0^1 d^{-1}(\sigma) \, \mathrm{d} \sigma = \frac{1}{g} \ln \left( 1 + \frac{\rho}{m} \right) > 0, \]
so that $x(0)=0$ and $x(1)=1$.
With this transformation, the PDE \eqref{A1} can be rewritten as
\be 
z_{tt} (x, t) - \lambda^2(x) z_{xx} (x, t) = 0, \label{eq:zlambda}
\ee
where
\[ \begin{aligned} \lambda (x) &:= \frac{1}{J \sqrt{\tilde{d}(x)}}, & \tilde{d}(x) &:= d(s(x)) = \frac{gm}{\rho} e^{gJx}.\end{aligned} \]
The PDE \eqref{eq:zlambda} has to be supplemented with  the boundary conditions 
\be
 \left\{ \begin{aligned}
z_x (0, t) &= 0, \\
z (1, t) &= X_\mathrm{p}(t), \\
\ddot{X}_\mathrm{p} (t) &= U(t),
\end{aligned} \right. 
\label{eq:zlambda2}
\ee
and the  initial conditions
\be
\left\{ \begin{array}{r@{\ }c@{\ }l}
z (x, 0) &= z^0 (x) &= y^0 (s), \\
z_t (x, 0) &= z^1 (x) &= y^1(s).
\end{array} \right. 
\label{eq:zlambda3}
\ee

\subsection{Second transformation}
Equation \eqref{eq:zlambda} is subsequently rewritten as a system of two first-order PDEs. Noticing that $(\partial_t - \lambda \partial_x)(\partial_t + \lambda \partial_x) z = -\lambda \lambda' z_x = (\partial_t + \lambda \partial_x)(\partial_t - \lambda \partial_x) z$, we infer that the Riemann invariants  $S := z_t + \lambda z_x$ and $D := z_t - \lambda z_x$ satisfy the system
\ba
S_t(x, t) - \lambda(x) S_x(x, t) &=& -\frac{\lambda'(x)}{2} (S(x, t)-D(x, t)),  \label{AB1}\\
D_t(x, t) + \lambda(x) D_x(x, t) &=& -\frac{\lambda'(x)}{2} (S(x, t)-D(x, t)). \label{AB2}
\ea
Setting
\be \bm{w}(x, t) = \mymatrix{u  (x,t)\\ v (x,t) }  := \mymatrix{e^{-\frac{1}{2}\ln(\lambda(x))} D(x, t) \\ e^{-\frac{1}{2}\ln(\lambda(x))} S(x, t)} \\ = \frac{1}{\sqrt{\lambda(x)}} \mymatrix{D (x, t) \\ S(x, t)}, \label{eq:wtf} \ee
we obtain the system
\be
 \bm{w}_t = \mymatrix{-\lambda(x) & 0 \\ 0 & \lambda(x)} \bm{w}_x + \mymatrix{0 & -\frac{\lambda'(x)}{2} \\ \frac{\lambda'(x)}{2} & 0} \bm{w} 
 \label{AB3}
 \ee
with boundary conditions 
\ba
u(0,t)&=&v(0,t), \label{AB3bis}\\
v(1,t)&=&W(t), \label{AB3ter} 
\ea
$W(t)$ denoting a feedback law that can be expressed in terms of the original control $U(t)=\ddot X_\mathrm{p} (t)$ and the angle $\theta (t) =y_s(1,t)$ as  
\[ W(t) = \sqrt{J} \sqrt[4]{d(1)} \left( \dot X_\mathrm{p}(t) + \sqrt{d(1)} \theta (t) \right). \]
System \eqref{AB3}-\eqref{AB3ter} is in the form required to apply the backstepping transformation from 
\cite{CVKB}. 

\section{Backstepping transformation and kernel calculation}\label{sec:back}

\subsection{Finite-time stabilization of the PDE}

In \cite{CVKB}, the following $2 \times 2$ linear hyperbolic PDE system is considered:
\ba
\bm{w}_t(x, t) 
 &=& \underbrace{\mymatrix{-\epsilon_1(x) & 0 \\ 0 & \epsilon_2(x)}}_{\bm{\Sigma}(x)} \bm{w}_x (x,t) + \underbrace{\mymatrix{0 & c_1(x) \\ c_2(x) & 0}}_{\bm{C} (x) } \bm{w} (x,t) \label{systw1}\\
u(0, t)&=& qv(0, t) \hspace{.5cm} v(1, t) = W(t), \label{systw2}
\ea
where $\bm{w} (x, t) = \mymatrix{u & v}^\intercal (x, t)$, $x \in [0, 1]$, $t \ge 0$, $c_1(x)$ and $c_2(x)$ are functions in $C^0 ([0, 1])$, $\epsilon_1(x)$ and $\epsilon_2(x)$ are strictly positive functions in 
$C^1([0, 1])$, and $q \in \R^*$. (For a discussion of the case $q=0$, see \cite[\textsection 3.5]{CVKB}.)

Clearly, system \eqref{AB3}-\eqref{AB3ter} is of the form \eqref{systw1}-\eqref{systw2} if we pick
\begin{eqnarray*}
&&\epsilon_1(x)=\epsilon_2(x)=\lambda(x)\quad \forall x\in [0,1],\\
&&c_2(x)=\frac{\lambda '(x)}{2} = -c_1(x)\quad \forall x\in [0,1], \\
&&q=1.
\end{eqnarray*}

A backstepping transform is defined between the original state variables $\bm{w} (x, t) = \mymatrix{u & v}^\intercal (x, t)$ and the target state variables $\bm{\gamma} (x, t) = \mymatrix{\alpha & \beta}^\intercal (x, t)$ that satisfy the target system 
\ba
\bm{\gamma}_t (x, t) &=& \bm{\Sigma}(x) \bm{\gamma}_x (x, t), \label{pdegamma1} \\
\alpha (0, t) &=& q\beta(0, t), \hspace{.6cm}  \beta (1, t) = 0. \label{pdegamma2}
\ea
Note that the target system \eqref{pdegamma1}-\eqref{pdegamma2}
is finite-time stable, since after some time we have transparent boundary conditions for both $\alpha$ and $\beta$.

\begin{definition}[Backstepping transform \cite{CVKB}]
The backstepping transform between the original state variables $\bm{w} (x, t)$ and the target state variables $\bm{\gamma} (x, t)$ is defined through a Volterra integral equation
\myeq{
\bm{\gamma}(x, t) = \bm{w}(x, t) - \int_0^x \bm{K}(x,\xi) \bm{w}(\xi,t)\;\mathrm{d}\xi \\
}{backsteppingtf}
with the direct kernels $\bm{K}(x,\xi )$ decomposed as
\[ \bm{K}(x,\xi) = \mymatrix{K^{uu} & K^{uv} \\ K^{vu} & K^{vv}} (x,\xi) .\]
The inverse transformation is given by
\myeq{
\bm{w}(x, t) = \bm{\gamma}(x, t) + \int_0^x \bm{L}(x,\xi)\bm{\gamma}(\xi,t)\;\mathrm{d}\xi
}{backsteppingtfinv}
with the inverse kernels
\[ \bm{L}(x,\xi) = \mymatrix{L^{\alpha\alpha} & L^{\alpha\beta} \\ L^{\beta\alpha} & L^{\beta\beta}} (x,\xi) . \]
\end{definition}

\begin{figure}[hbtp]
	\begin{center}
	\begin{tikzpicture}[scale=3]
		\draw[->] (-0.2, 0) -- (1.2, 0);
		\draw[->] (0, -0.2) -- (0, 1.2);
		\draw[fill=gray!50] (0, 0) -- (1, 1) -- (1,0) -- (0, 0);
		\node[right] at (1.2, 0) {$x$};
		\node[above] at (0, 1.2) {$\xi$};
		\draw (-0.05, 1) -- (0.05, 1);
		\node[below] at (-0.18, 0) {$0$};
		\node[left] at (0, 1) {$1$};
		\node[below] at (1, 0) {$1$};
		\node at (.833, .166) {$\mathcal{T}$};
		\draw [black, fill=black] (0,0) circle (0.15ex);
		\draw [black, fill=black] (0.333,0) circle (0.15ex);
		\draw [black, fill=black] (0.666,0) circle (0.15ex);
		\draw [black, fill=black] (1,0) circle (0.15ex);
		\draw [black, fill=black] (0.333,0.333) circle (0.15ex);
		\draw [black, fill=black] (0.666,0.333) circle (0.15ex);
		\draw [black, fill=black] (1,0.333) circle (0.15ex);
		\draw [black, fill=black] (0.666,0.666) circle (0.15ex);
		\draw [black, fill=black] (1,0.666) circle (0.15ex);
		\draw [black, fill=black] (1,1) circle (0.15ex);
		\draw (0.333,-.2) -- (0.333,0);
		\draw (0.666,-.2) -- (0.666,0);
		\draw (1,0.333) -- (1.2,0.333);
		\draw (1,0.666) -- (1.2,0.666);
		\draw[<->] (0.333,-.2) -- (0.666,-.2);
		\draw[<->] (1.2,0.333) -- (1.2,0.666);
		\node[above] at (.5,-.23) {$\Delta x$};
		\node[left] at (1.25,.5) {\rotatebox{90}{$\Delta \xi$}};
	\end{tikzpicture}
	\end{center}
	\caption{Domain $\mathcal{T}$ with uniform discretization grid}
	\label{fig:tdom} 
\end{figure}
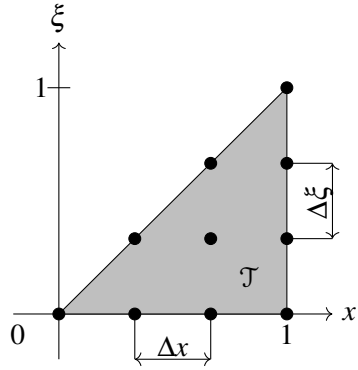

By a direct calculation, it is proven in \cite{CVKB} that the direct kernels $\bm{K}$ are the solution of the  following Goursat system of two $2\times2$ first-order hyperbolic PDEs
\myeq{ \left\{ \begin{aligned}
\epsilon_1(x)K^{uu}_x + \epsilon_1(\xi)K^{uu}_\xi &= -\epsilon_1'(\xi)K^{uu} - c_2(\xi)K^{uv}, \\
\epsilon_1(x)K^{uv}_x - \epsilon_2(\xi)K^{uv}_\xi &= \epsilon_2'(\xi)K^{uv}-c_1(\xi)K^{uu}, \\
\epsilon_2(x) K^{vu}_x - \epsilon_1(\xi)K^{vu}_\xi &= \epsilon_1'(\xi)K^{vu}+c_2(\xi)K^{vv}, \\
\epsilon_2(x) K^{vv}_x + \epsilon_2(\xi)K^{vv}_\xi &= -\epsilon_2'(\xi)K^{vv} + c_1(\xi) K^{vu}
\end{aligned} \right.
}{goursatone}
on the triangular domain $\mathcal{T} = \{ (x,\xi) \; | \; 0\le\xi\le x\le1 \}$ (see Fig. \ref{fig:tdom})
with the boundary conditions
\[ \left\{ \begin{aligned}
K^{uu}(x, 0) &= \frac{\epsilon_2(0)}{q\epsilon_1(0)}K^{uv}(x, 0) ,\\
K^{uv}(x, x) &= \frac{c_1(x)}{\epsilon_1(x)+\epsilon_2(x)} ,\\
K^{vu}(x, x) &= -\frac{c_2(x)}{\epsilon_1(x)+\epsilon_2(x)} ,\\
K^{vv}(x, 0) &= \frac{q\epsilon_1(0)}{\epsilon_2(0)} K^{vu}(x, 0)
\end{aligned} \right. \]
for $x\in [0,1]$. 

It is also shown in \cite{CVKB} that the indirect kernels $\bm{L}$ satisfy a similar Goursat system of $4\times 4$ first-order hyperbolic PDEs defined 
on $\mathcal{T}$. 

Note that it is also possible to compute the inverse kernels $\bm{L}$ from the direct kernels $\bm{K}$ \cite{AA}:
\myeq{\bm{L}(x, \xi) = \bm{K}(x, \xi) + \int_\xi^x \bm{K} (x, \sigma) \bm{L}(\sigma, \xi) \, \mathrm{d} \sigma.}{relKL}

If the coefficients $\epsilon_i(x)$ and $c_i(x)$ in the Goursat systems are constant (independent of $x$ or $\xi$), 
explicit solutions for the kernels have been obtained  in \cite{VK}. In the general case, numerical methods have been designed  in \cite{AA}. These methods approximate the left-hand side of the equations of the Goursat system (for instance \eqref{eq:goursatone}) by a directional derivative and interpolate between values at points of a triangular half of an equidistant square grid as shown in Figure \ref{fig:tdom}.

\lemmaspace
\begin{lemma}[Finite-time feedback law for system \eqref{systw1}-\eqref{systw2}
 \cite{CVKB}] System \eqref{systw1}-\eqref{systw2} is mapped on \eqref{pdegamma1}-\eqref{pdegamma2}
  with finite-time dynamics by the backstepping transform \eqref{eq:backsteppingtf} with kernels satisfying \eqref{eq:goursatone} and by using the feedback law
\myeq{ W(t) = \int_0^1 K^{vu}(1, \xi) u(\xi,t) \, \mathrm{d}\xi + \int_0^1 K^{vv}(1,\xi) v(\xi,t) \, \mathrm{d} \xi. }{dummycoronlaw}
\end{lemma}

\subsection{Finite-time stabilization of the hybrid PDE-ODE system}

In the case of the hybrid (PDE-ODE) crane model \eqref{A1}-\eqref{A4}, one does not control
the dynamics of the system cable--platform through the
boundary condition $v(1,t)=W(t)$, but rather through the acceleration $\ddot{X}_\mathrm{p}(t) =U(t)$ of the platform.

In that case, the transparent boundary condition $\beta (1,t) = 0$ in \eqref{pdegamma2} is not 
necessarily satisfied. 

Our goal is to design a feedback law  $U(t)$ in such a way that both $\beta(1, t)$ and a new  variable $\phi (t)$
vanish after a finite time depending on the initial data. Next, both the target variable $\bm{\gamma}(x, t)$ 
and $X_\mathrm{p}(t)$ will be stabilized in a finite time, after which the entire system will be at the equilibrium position. This is the subject of the next section.

\section{Finite-time stability of the complete system}\label{sec:res}

A finite-time stabilizer $U(t)$ for the crane and the cable with affine tension is constructed as follows. 
Our goal is  to have the condition $\beta (1, t)=0$ satisfied for $t$ large enough.
 We write 
 \ba
 2\frac{\dot X_\mathrm{p}}{\sqrt{\lambda (1) }}
 &=& 2\frac{z_t (1,t)}{\sqrt{\lambda (1)}} \nonumber \\
 &=& \mymatrix{1 & 1} \bm{w}(1, t) \nonumber \\
 &=& \mymatrix{1 & 1}  \bm{\gamma} (1, t) + \int_0^1 \mymatrix{1 & 1} \bm{L}(1, \xi ) \bm{\gamma} (\xi , t) \, \mathrm{d}\xi \nonumber \\
 &=& \alpha (1, t) + \beta (1, t) \nonumber \\
&&\  + \int_0^1 [ \big( L^{\alpha \alpha} (1, \xi) + L^{\beta \alpha }(1, \xi) \big) 
 \alpha (\xi , t) +  \big( L^{\alpha \beta} (1, \xi) + L^{\beta \beta }(1, \xi)\big) \beta (\xi , t) ]\, \mathrm{d}\xi 
 \label{KK1}
 \ea  
Let
\begin{equation} 
\phi (t) := \frac{2}{\sqrt{\lambda(1)}} X_\mathrm{p}(t) + \int_0^1 \left[ a(x) \alpha(x, t) + b(x) \beta(x, t) \right] \, \mathrm{d} x \label{eq:defphi} 
\end{equation}
where $a(x)$ and $b(x)$ are two functions to be defined. Then we have
\ba
\dot \phi (t) &=&\frac{2}{\sqrt{\lambda (1)}} \dot X_\mathrm{p}(t)  + \int_0^1 [a\alpha _t + b\beta _t] \, \mathrm{d}x \nonumber\\
&=&\frac{2}{\sqrt{\lambda (1)}} \dot X_\mathrm{p}(t)  + 
 \int_0^1 [ -a\lambda \alpha _x + b\lambda \beta _x]\, \mathrm{d}x \nonumber\\
&=&  \frac{2}{\sqrt{\lambda (1)}} \dot X_\mathrm{p}(t)  
+\left[ 
-a(x) \lambda (x) \alpha (x, t) +b(x) \lambda (x) \beta (x, t) 
\right]_{x=0}^1 
 +\int_0^1 [ (a\lambda)_x \alpha  - (b\lambda )_x \beta  ] \, \mathrm{d}x\quad \label{KK2}
 \ea
where we used an integration by parts in the last line.
Replacing in \eqref{KK2}  $\dot X_\mathrm{p}(t)$ by its expression in \eqref{KK1} results in
 \ba
 \dot \phi (t)
 &=& \alpha (1, t) + \beta (1, t) + \left[ -a\lambda \alpha + b\lambda \beta \right] _{x=0}^1  \nonumber\\
 &&+\int_0^1 \left[ (L^{\alpha \alpha} (1,x) +L^{\beta \alpha} (1,x) + (a\lambda )_x )\alpha 
 +(L^{\alpha \beta} (1,x)+L^{\beta\beta} (1,x) -(b\lambda )_x)\beta \right] \, \mathrm{d}x.  \label{KK3}
 \ea
 We define the functions $a$ and $b$ as 
\ba
a(x) &:=& \frac{1}{\lambda(x)} \left(a_0 - \int_0^x \left( L^{\alpha\alpha} (1,\xi) + L^{\beta\alpha} (1,\xi) \right) \, \mathrm{d} \xi \right),  \label{PPP1}\\
b(x) &:=&  \frac{1}{\lambda(x)} \left(b_0 + \int_0^x \left( L^{\alpha\beta} (1,\xi) + L^{\beta\beta} (1,\xi) \right) \, \mathrm{d} \xi \right)     \label{PPP2}
\ea
with $a_0,b_0$ some constants still to choose.
Then  
\[
L^{\alpha \alpha} (1,x) +L^{\beta \alpha} (1,x) + (a\lambda )_x = 0 =
 L^{\alpha \beta} (1,x)+L^{\beta\beta} (1,x) -(b\lambda )_x, \]
 and hence \eqref{KK3} becomes
 \be
 \dot \phi = [ 1-(a\lambda)(1)] \alpha (1, t) + [1+(b\lambda) (1)] \beta (1, t)
 +(a_0-b_0)\alpha (0, t) 
 \ee
 for $\alpha (0, t)=\beta (0, t)$, ($q=1$). Now we let
 \[ a_0 = b_0 := 1 + \int_0^1 \left( L^{\alpha\alpha} (1, \xi) + L^{\beta \alpha} (1, \xi) \right) \, \mathrm{d} \xi, \]
so that 
\[
1-(a\lambda )(1)=0=a_0-b_0. 
\]
We obtain
\ba
\dot{\phi}(t) &=&[1+ ( b \lambda )(1)]\beta (1, t) \nonumber\\
&=& \left[ 1+b_0+\int_0^1[ L^{\alpha \beta} (1,x)+L^{\beta\beta}(1,x) ] \, \mathrm{d}x \right] \beta (1, t) \nonumber\\
&=& \mu \beta (1, t) \label{eq:phidotmubeta}
\ea
with 
\be 
\mu  := 2 + \int_0^1 \left( L^{\alpha\alpha}(1,x) + L^{\beta\alpha}(1,x) \right. \\ \left. + L^{\alpha\beta}(1,x) + L^{\beta\beta}(1,x) \right) \, \mathrm{d} x \label{eq:defmu}.
\ee

Later on, we shall define $\beta (1,t)$ from $\dot \phi (t)$ using \eqref{eq:phidotmubeta}. It is thus important
to prove that $\mu\ne 0$. This is done in the following proposition. 

\begin{proposition}
\label{propo1}
The kernels $L$ satisfy
\begin{eqnarray*}
&&L^{\alpha \alpha } (x, \xi ) = L^{\beta \beta } (x, \xi )  \ge 0\quad \forall (x, \xi ) \in {\mathcal T},\\
&&L^{\alpha \beta } (x, \xi )   = L^{\beta \alpha } (x, \xi)  \ge 0 \quad \forall (x, \xi ) \in {\mathcal T},
\end{eqnarray*}
and hence $\mu \ge 2$. 
\end{proposition}
\begin{proof}
Recall that $\epsilon _1(x)=\epsilon _2 (x) =\lambda (x)$ and that $c_1(x) = - c_2(x) = -\frac{\lambda ' (x) }{2}$. 
Let $C_1>0$ and $C_2>0$ be defined by 
\[
\lambda(x)=\frac{1}{ J \sqrt{\frac{gm}{\rho}e^{gJx}}} =:C_1e^{-C_2x}.
\]
Then the Goursat system satisfied
by the kernels $L$ reads (see \cite{CVKB})
\ba
\lambda (x) L^{\alpha\alpha}_x + \lambda (\xi ) L^{\alpha\alpha} _\xi &=& -\lambda ' ( \xi ) L^{\alpha \alpha}  - \frac{\lambda ' (x) }{2} L^{\beta \alpha}, 
\quad (x,\xi ) \in {\mathcal T}, \label{TT1}\\
\lambda (x) L^{\alpha\beta}_x - \lambda (\xi ) L^{\alpha\beta} _\xi &=& \lambda ' ( \xi ) L^{\alpha \beta}  - \frac{\lambda ' (x) }{2} L^{\beta \beta},
\quad (x,\xi ) \in {\mathcal T}, \label{TT2} \\
\lambda (x) L^{\beta\alpha}_x - \lambda (\xi ) L^{\beta\alpha} _\xi &=& \lambda ' ( \xi ) L^{\beta \alpha}  -  \frac{\lambda ' (x) }{2} L^{\alpha \alpha},
\quad (x,\xi ) \in {\mathcal T}, \label{TT3} \\
\lambda (x) L^{\beta\beta}_x + \lambda (\xi ) L^{\beta\beta} _\xi &=& - \lambda ' ( \xi ) L^{\beta \beta}  -  \frac{\lambda ' (x) }{2} L^{\alpha \beta}, 
\quad (x,\xi ) \in {\mathcal T}, \label{TT4}
\ea
with the boundary conditions 
\ba
L^{\alpha\alpha} (x,0)&=& L^{ \alpha \beta } (x,0), \quad 0\le x\le 1, \label{TT5}\\
L^{\alpha\beta} (x,x)&=& -\frac{\lambda '(x)}{4\lambda (x) } = \frac{C_2}{4}, \quad 0\le x\le 1, \label{TT6}\\
L^{\beta\alpha} (x,x)&=& \frac{C_2}{4}, \quad 0\le x\le 1, \label{TT7}\\
L^{\beta\beta} (x,0)&=& L^{\beta \alpha } (x,0), \quad 0\le x\le 1. \label{TT8}
\ea
The existence and uniqueness of the solution $(L^{\alpha \alpha}, L^{\alpha \beta}, L^{\beta \alpha}, L^{\beta\beta})\in [C^0( {\mathcal T} )]^4$
was already established in \cite[Theorem A.1]{CVKB}. However, the claims in Proposition \ref{propo1} are still to be justified. 

We introduce the following $2\times 2$ system:
\ba
\lambda (x) f_x + \lambda (\xi ) f_\xi &=& -\lambda ' (\xi ) f + \mu (x) g, \quad (x,\xi) \in {\mathcal T}, \label{TT11}\\
\lambda (x) g_x - \lambda (\xi ) g_\xi &=& \lambda ' (\xi ) f + \mu (x) f, \quad (x,\xi) \in {\mathcal T}, \label{TT12}
\ea
with the boundary conditions
\ba
f(x,0)&=&g(x,0), \quad 0\le x\le 1, \label{TT13}\\
g(x,x)&=& C, \quad 0\le x\le 1, \label{TT14}
\ea
where $C:=C_2/4>0$ and $\mu (x)=-\lambda '(x)/2>0$ for $x\in [0,1]$. If we prove the existence of a solution $(f,g)$ of \eqref{TT11}-\eqref{TT14}
in $[C^0( {\mathcal T} )]^2$, then
setting $L^{\alpha\alpha}=L^{\beta\beta}:=f$, $L^{\alpha\beta}=L^{\beta\alpha}:=g$, we obtain a solution of \eqref{TT1}-\eqref{TT8} in 
$[C^0( {\mathcal T} )]^4$. The uniqueness of such a solution reduces the study of system \eqref{TT1}-\eqref{TT8} to those of system \eqref{TT11}-\eqref{TT14}. Therefore, Proposition \ref{propo1} is a direct consequence of the following 
lemma.
\begin{lemma}
\label{lemm1}
Let $C>0$, $\lambda \in C^1([0,1])$ and $\mu \in C^0([0,1])$ with $\lambda (x) >0$ and $\mu (x)>0$ for all $x\in [0,1]$. Then there exists a unique solution $(f,g)\in [C^0({\mathcal T}, \R_+ )]^2$ of \eqref{TT11}-\eqref{TT14}, where $\R _+=[0,+\infty )$. 
\end{lemma}
\noindent
{\em Proof of Lemma \ref{lemm1}.} We first express $(f,g)$ as a fixed-point of a map from $[C^0({\mathcal T}, \R_+)]^2$ into itself. Next, 
reducing the domain for $(x,\xi)$ to a ``trapezoid'', we prove that the above map is a contraction. The full domain $\mathcal T$ is covered by
iterating the above construction. 

The characteristic system for \eqref{TT11} reads 
\begin{eqnarray*}
\frac{\mathrmd x}{\mathrmd s} &=& \lambda (x), \\
\frac{\mathrmd\xi}{\mathrmd s} &=& \lambda (\xi ),\\
\frac{\mathrmd z}{\mathrmd s} &=& -\lambda ' (\xi ) z + \mu (x) g  
\end{eqnarray*}
with the boundary conditions $x(0)=x_0$, $\xi (0)=0$, and $z(0) = g(x_0,0)$.

Let $\Lambda (x) := \int_0^x \frac{\mathrmd s}{\lambda (s)}$. Then one readily obtains 
\[
x(s)=\Lambda ^{-1} (\Lambda (x_0)+s), \ \ \xi (s) =\Lambda ^{-1} (s)
\]
and 
\[
f(x(s),\xi (s))= z(s) = e^{-\int_0^s \lambda ' (\xi (\sigma )) \,\mathrmd\sigma } 
\left( 
g(x_0,0)+ 
\int_0^s e^{\int_0^\sigma \lambda ' (\xi (\tau )) \,\mathrmd\tau }  
 \mu (x(\sigma ))\, g(x(\sigma ) , \xi (\sigma )) \,\mathrmd\sigma
 \right) .
\] 
Since $s=\Lambda (\xi )$, $x_0=\Lambda ^{-1} (\Lambda (x)-\Lambda (\xi ))$, we have that 
$x(\sigma ) =\Lambda ^{-1} (\Lambda (x) -\Lambda (\xi ) +\sigma )$, $\xi (\sigma ) =\Lambda ^{-1} (\sigma )$ and 
\begin{multline}
f(x,\xi ) = e^{-\int_0^{\Lambda (\xi)} \lambda '(\Lambda ^{-1} (\sigma )) \,\mathrmd\sigma  }
\left[ g(\Lambda ^{-1} (\Lambda (x) -\Lambda (\xi)), 0) \right. \\
+\left. \int_0^{\Lambda (\xi )} e^{\int_0^\sigma \lambda ' (\Lambda ^{-1} (\tau ) ) \,\mathrmd\tau } 
\mu (\Lambda ^{-1} (\Lambda (x) -\Lambda (\xi ) +\sigma )) 
g(\Lambda ^{-1} (\Lambda (x) -\Lambda (\xi ) +\sigma ), \Lambda ^{-1} (\sigma )) \,\mathrmd\sigma 
\right]
\label{W1} 
\end{multline}
Similarly, the characteristic system for \eqref{TT12} reads
\begin{eqnarray*}
\frac{\mathrmd x}{\mathrmd s} &=& \lambda (x), \\
\frac{\mathrmd\xi}{\mathrmd s} &=& - \lambda (\xi ),\\
\frac{\mathrmd z}{\mathrmd s} &=& \lambda ' (\xi ) z + \mu (x) f  
\end{eqnarray*}
with the boundary conditions $x(0)=x_0$, $\xi (0)=x_0$, and $z(0) = C$.

We readily obtain
\[
x(s)=\Lambda ^{-1} (\Lambda (x_0)+s), \ \ \xi (s) =\Lambda ^{-1} (\Lambda (x_0) -s)
\]
so that 
\[
s=\frac{\Lambda (x) -\Lambda (\xi ) }{2}, \ \  x_0=\Lambda ^{-1} \left(\frac{\Lambda (x)  +\Lambda (\xi )}{2}\right), 
\]
and 
\[
g(x(s),\xi (s))= z(s) = e^{\int_0^s \lambda ' (\xi (\sigma )) \,\mathrmd\sigma } 
\left( 
C+ 
\int_0^s e^{-\int_0^\sigma \lambda ' (\xi (\tau )) \,\mathrmd\tau }  
 \mu (x(\sigma ))\, f(x(\sigma ) , \xi (\sigma )) \,\mathrmd\sigma
 \right) .
\] 
It follows that 
\begin{multline}
g(x,\xi )=  e^{\int_0^{\frac{\Lambda (x) -\Lambda (\xi)}{2}}  \lambda ' \left( \Lambda ^{-1}\left( \frac{\Lambda (x ) + \Lambda (\xi )}{2} - \sigma \right)\right) \,\mathrmd\sigma } 
\left(  C  + 
\int_0^{\frac{\Lambda (x) -\Lambda (\xi ) }{2}} e^{-\int_0^\sigma \lambda ' \left( \Lambda ^{-1} \left(\frac{\Lambda (x ) + \Lambda (\xi )}{2}  - \tau \right)\right) \,\mathrmd\tau }  
 \mu \left(\Lambda ^{-1} \left(\frac{\Lambda (x) + \Lambda (\xi) }{2}   + \sigma \right)\right)     \right. \\
  \left.  \times f\left( \Lambda ^{-1} \left(\frac{\Lambda (x ) + \Lambda (\xi ) }{2} + \sigma \right) , 
 \Lambda ^{-1} \left(\frac{\Lambda (x ) + \Lambda (\xi ) }{2} - \sigma \right) 
 \right)\, \,\mathrmd\sigma
 \right) .
 \label{W1bis}
\end{multline}
Replacing in \eqref{W1bis} $f$ by its expression in \eqref{W1}, we arrive to $g=\Gamma (g)$, where 
\begin{multline}
\Gamma (g)(x,\xi ):=  \\
e^{\int_0^{\frac{\Lambda (x) -\Lambda (\xi)}{2}}  \lambda ' \left( \Lambda ^{-1} \left( \frac{\Lambda (x ) + \Lambda (\xi )}{2} 
-\sigma \right)\right) \,\mathrmd\sigma } 
\bigg(  C  + 
\int_0^{\frac{\Lambda (x) -\Lambda (\xi ) }{2}} e^{-\int_0^\sigma \lambda ' \left( \Lambda ^{-1} \left(\frac{\Lambda (x ) + \Lambda (\xi )}{2}  -\tau \right)\right) \,\mathrmd\tau }  
 \mu \left(\Lambda ^{-1} \left(\frac{\Lambda (x) + \Lambda (\xi) }{2}   + \sigma \right)\right)      \\
 \left. 
\times e^{-\int_0^{\frac{\Lambda (x) + \Lambda (\xi ) }{2} -\sigma }   \lambda '(\Lambda ^{-1} (\sigma )) \,\mathrmd\sigma  }
\bigg\{ g(\Lambda ^{-1} (2\sigma ), 0) \right. +
\int_0^{\frac{\Lambda (x) + \Lambda (\xi )}{2} -\sigma  } e^{\int_0^\kappa \lambda ' (\Lambda ^{-1} (\tau ) ) \,\mathrmd\tau } 
\mu (\Lambda ^{-1} (2\sigma + \kappa ))  \\
 \times g(\Lambda ^{-1} (2\sigma + \kappa ), \Lambda ^{-1} (\kappa ))\,\mathrmd\kappa \bigg\} \,\mathrmd\sigma \bigg).    
 \label{W3}
\end{multline}
For $0\le \epsilon_1\le \epsilon _2\le 1$, let 
\begin{eqnarray*}
{\mathcal T} _{\epsilon _1, \epsilon _2} &:=& \{ (x, \xi ) \in {\mathcal T}; \  \ \Lambda (\epsilon _1) \le \Lambda (x)-\Lambda (\xi ) \le 
\Lambda (\epsilon _2 )\}, \\
E_{\epsilon _1, \epsilon _2} &:=&   C^0( {\mathcal T} _{\epsilon _1, \epsilon _2} , \R _+ ).
\end{eqnarray*}
Let $\epsilon \in (0,1)$. Noticing that the restriction of $\Gamma (g)$ to ${\mathcal T} _{0, \epsilon }$ depends only on the restriction of $g$ to 
${\mathcal T} _{0, \epsilon }$, we can define a map $g\in E_{0, \epsilon }\to \Gamma (g) \in E_{0, \epsilon }$
and find a constant $K= K\big(\Vert \mu\Vert _{L^\infty (0,1) }, \Vert \lambda ' \Vert _{L^\infty (0,1)}, 
\Vert \Lambda \Vert  _{L^\infty (0,1)}\big)>0$ such that 
\[
\Vert \Gamma (g_1)-\Gamma (g_2) \Vert _{L^\infty({\mathcal T} _{0, \epsilon })}
\le K \sup _{(x,\xi) \in {\mathcal T} _{0, \epsilon }} \frac{\Lambda (x) -\Lambda (\xi )}{2} \Vert g_1-g_2\Vert_{L^\infty ({\mathcal T} _{0, \epsilon })}
\le  \frac{K \Lambda (\epsilon )}{2} \Vert g_1-g_2\Vert_{L^\infty( {\mathcal T} _{0, \epsilon })}.
\]
For $\epsilon >0$ small enough, we have that $K\Lambda (\epsilon)/2<1$, so that the map $g\in E_{0, \epsilon }\to \Gamma (g) \in E_{0, \epsilon }$
is a contraction, and hence it has a unique fixed point by the contraction mapping theorem. Using \eqref{W1}, this yields a unique solution $(f,g)\in E_{0,\epsilon}^2$. 
Proceeding in the same way (using the values of $g$  on the characteristic curve $\Lambda (x)-\Lambda (\xi ) = \Lambda (\epsilon)$ 
computed in the previous step), we can extend $g$ and $f$ on the trapezoids ${\mathcal T}_{\epsilon, 2\epsilon}$,
${\mathcal T}_{2\epsilon, 3\epsilon}$, etc. Finally, the functions $g$ and $f$ are defined on the whole domain ${\mathcal T}$ and they take nonnegative values. 
The proofs of Lemma \ref{lemm1} and of Proposition \ref{propo1} are complete.  
\end{proof} 
\begin{remark}
Actually, it follows from \eqref{W1} and \eqref{W1bis} that $f$ and $g$ take strictly positive values on $\mathcal T$. The same is true for
$L^{\alpha\alpha} = L^{\beta\beta}$ and $L^{\alpha \beta}=L^{\beta \alpha}$.
\end{remark}

The control input $ \ddot{X}_\mathrm{p}(t) =U(t)$ is designed in such a way that 
$\phi(t)$ obeys the following dynamics with finite-time stability 
(see \cite{BB,haimo})
\myeq{ \ddot{\phi} (t) + \lfloor \dot{\phi} (t) \rceil^\halpha + \lfloor \phi(t) \rceil^\hbeta = 0, }{haimo}
where $0<\halpha<1$, $\hbeta \ge\frac{\halpha}{2-\halpha}$ and $\lfloor x\rceil^\hbeta := \mathrm{sign}(x)|x|^\hbeta$. 
After substitution of the expression \eqref{eq:defphi} in \eqref{eq:haimo}, one obtains the feedback law
\begin{multline}
U(t):= -\frac{\sqrt{\lambda(1)}}{2} \left( \int_0^1 \left( a(x) \alpha_{tt}(x, t) + b(x) \beta_{tt}(x, t) \right) \, \mathrm{d} x
\right. \\
 + \left\lfloor \frac{2}{\sqrt{\lambda(1)}} \dot{X}_\mathrm{p}(t) + \int_0^1 \left( a(x) \alpha_t(x, t) + b(x) \beta_t(x, t)\right) \, \mathrm{d} x \right\rceil^\halpha
\\ \left. + \left\lfloor \frac{2}{\sqrt{\lambda(1)}} X_\mathrm{p}(t) + \int_0^1 \left( a(x) \alpha (x, t) + b(x) \beta (x, t) \right) \, \mathrm{d} x \right\rceil^\hbeta \right) .
\label{eq:laloi}\end{multline}
\normalsize

We claim that with this feedback law, the crane and the cable with affine tension are stabilized in finite time. 
We first establish the finite-stability of the system in the new variables $ ( \bm{\gamma}, \phi)$, and next we go back to the original variables $(z, X_\mathrm{p})$.  
Let $H := [L^2(0, 1)]^2\times \R^2$ be endowed with the norm 
\[
\Vert (\alpha , \beta , \phi ^0, \phi ^1)\Vert _H^2 :=
\int_0^1[ \alpha (x)^2+\beta (x)^2 ]\, \mathrm{d}x + |\phi ^0|^2 + |\phi ^1|^2. 
\]
\lemmaspacetwo
\begin{theorem}
\label{thm1}
Let $\psi, \zeta\in \R$ with $0<\psi <1$ and $\zeta \ge \frac{\psi}{2-\psi}$. Then 
the system \eqref{pdegamma1}, \eqref{eq:haimo}, with the boundary conditions
\be
\alpha (0, t)= \beta (0, t), \quad \beta (1, t) =\mu  ^{-1} \dot \phi (t),
\label{U1}
\ee
 is well-posed and  globally finite-time stable in the Hilbert space 
 $H$. More precisely, there exists a nondecreasing 
 function $T_1:(0, +\infty)\to (0, +\infty)$ such that for all $(\alpha ^0,\beta ^0, \phi ^0, \phi ^1)
 \in H$, the solution of \eqref{pdegamma1}, \eqref{eq:haimo}, \eqref{U1}
 and $(\alpha (.,0), \beta(.,0), \phi(0), \dot \phi (0))=(\alpha ^0, \beta ^0, \phi ^0, \phi ^1)$
satisfies 
\ba
\alpha (x, t)=\beta (x, t)=0,&&  \forall x\in [0, 1],\ 
\forall t\ge T_1(\Vert (\alpha ^0, \beta ^0, \phi ^0, \phi ^1)\Vert _H), \label{U2}\\
\phi (t)=0,&& \forall t\ge T_1(\Vert (\alpha ^0, \beta ^0, \phi ^0, \phi ^1)\Vert _H).
\label{U3}
\ea 
\end{theorem}
\begin{proof}
We use the cascaded structure of the PDE-ODE system. 
The Cauchy problem  
\ba
\ddot{\phi} (t) + \lfloor \dot{\phi} (t) \rceil^\halpha + \lfloor \phi(t) \rceil^\hbeta &=& 0, \quad t\in [0, T], \label{PP1}\\
 (\phi (0), \dot \phi (0)) &=& (\phi ^0, \phi ^1) \label{PP2}
\ea
 admits solutions $\phi \in C^2([0, T])$ for some (possibly small) $T>0$ by Peano's theorem. 
The solution is actually unique forward, defined on $[0, T]$ for all $T>0$, and a global finite-time stability occurs for system \eqref{eq:haimo} (see \cite{BR,BB,haimo}).
Let $T_0:(0, +\infty)\to (0, +\infty)$ be a nondecreasing function such that
\[ \phi (t) =0\  \textrm{ if } \  t\ge T_0( ||(\phi ^0, \phi ^1) || ). \] 

where $||(\phi ^0, \phi ^1) ||  =(|\phi ^0|^2 + |\phi ^1|^2) ^\frac{1}{2}$.  Next we consider the system 
\ba
\beta_t=\lambda (x)\beta_x,&& x\in (0, 1), \ t>0,\label{U11}\\
\beta (1, t)=\mu ^{-1} \dot \phi (t),&& t>0, \label{U12}\\
\beta (x, 0)=\beta ^0(x),&& x\in (0, 1).\label{U13}    
\ea
Following \cite[Definition 2.1 p. 25]{coron}, for given $T>0$, $\beta ^0\in L^2(0, 1)$ and $\phi \in H^1(0, T)$, we say that a solution of the Cauchy problem \eqref{U11}-\eqref{U13} is a function
$\beta\in C^0([0, T],L^2(0, 1))$ such that, for every $\tau\in [0, T]$ and every 
$\eta \in C^1([0, 1]\times [0,\tau ])$ such that 
$\eta (0, t)=0$ for all $t\in [0, T]$, we have 
\[
-\int_0^\tau\int_0^1 (\eta _t-(\lambda \eta)_x)\beta \, \mathrm{d}x\mathrm{d}t 
- \frac{\lambda (1)}{\mu} \int_0^\tau \eta (1, t) \dot \phi (t) \, \mathrm{d}t
+\int_0^1[\eta (x,\tau)\beta (x,\tau)-\eta (x, 0)\beta ^0(x)]\, \mathrm{d}x=0. 
\]  
Then, following closely \cite[Theorem 2.4 p. 27]{coron}, it can be proved that
there exists a unique solution $\beta\in C^0([0, T],L^2(0, 1))$ of the Cauchy problem \eqref{U11}-\eqref{U13}. 
Furthermore, using the writing $\lambda(x)=C_1e^{-C_2x}$,
we readily obtain with the method of characteristics that $\beta$ is given by 
\be
\beta(x, t) = 
\left\{ 
\begin{array}{ll} 
\beta ^0 \left( \frac{1}{C_2}\ln (e^{C_2x} +C_1C_2t )  \right) 
&\textrm{if } t < \frac{e^{C_2} -e^{C_2x }}{ C_1C_2} , \\
 \mu ^{-1} \dot\phi \left( t+ \frac{e^{C_2x} - e^{C_2} }{C_1C_2} \right) 
&\textrm{if } t\ge \frac{e^{C_2} -e^{C_2x }}{ C_1C_2} \cdot 
\end{array} 
\right. 
\label{VV1}
\ee 
 In particular, since $\phi (t)=0$ for $t\ge T_0 ( \Vert (\phi ^0, \phi ^1) \Vert )$, we infer that 
 $\beta (x, t)=0$ for $x\in [0, 1]$ and $t\ge T_0 ( \Vert (\phi ^0, \phi ^1) \Vert )
+\frac{e^{C_2}-e^{C_2x}}{C_1C_2}$. 

Note also that $\beta$ has a trace $\beta (0,\cdot )\in L^2(0,T)$ for all $T>0$, which is given by 
\[
\beta(0, t) = 
\left\{ 
\begin{array}{ll} 
\beta ^0 \left( \frac{1}{C_2}\ln (1 +C_1C_2t )  \right) 
&\textrm{if} \ t < \frac{e^{C_2} - 1 }{ C_1C_2} , \\
 \mu ^{-1} \dot\phi \left( t+ \frac{1 - e^{C_2} }{C_1C_2} \right) 
&\textrm{if} \ t\ge \frac{e^{C_2} - 1 }{ C_1C_2} \cdot 
\end{array} 
\right. 
\]

Similarly, it can be proved that for any $\alpha ^0\in L^2(0, 1)$, there exists for all $T>0$ a unique solution $\alpha
\in C^0 ( [0, T],L^2(0, 1))$ of the Cauchy problem
\ba
\alpha _t= - \lambda (x) \alpha _x,&& x\in (0, 1), \ t>0,\label{U21}\\
\alpha (0, t)=\beta  (0, t),&& t>0,                                  \label{U22}\\
\alpha  (x, 0)=\alpha  ^0(x),&& x\in (0, 1).                     \label{U23}    
\ea
Furthermore, it is given by
\be
\alpha (x, t) = 
\left\{ 
\begin{array}{ll} 
\alpha ^0 \left( \frac{1}{C_2}\ln (e^{C_2x} -C_1C_2t )  \right) 
&\textrm{if } t < \frac{e^{C_2x} - 1 }{ C_1C_2} ,   \\
 \beta  \left( 0, t - \frac{ e^{C_2x} -1 }{C_1C_2} \right) 
&\textrm{if } t\ge \frac{e^{C_2x} -1 }{ C_1C_2} \cdot 
\end{array} 
\right. 
\label{VV2}
\ee
 Since $\alpha (0, t) = \beta (0, t)=0$ for $t\ge T_0 ( \Vert (\phi ^0, \phi ^1) \Vert ) +\frac{e^{C_2}  -1 }{ C_1C_2 }$, we infer that 
 $\alpha  (x, t)=0$ for $x\in [0, 1]$ and $t\ge T_0 ( \Vert (\phi ^0, \phi ^1) \Vert )
 + \frac{ e^{C_2} + e^{C_2x} -2 }{ C_1C_2}  $. 
Thus 
\begin{eqnarray*}
 \phi(t)=0, && \textrm {for } \ 
t \ge T_1(\Vert (\alpha ^0, \beta ^0, \phi ^0, \phi ^1)\Vert _H) :=
 T_0 ( \Vert (\phi ^0, \phi ^1) \Vert )
+ \frac{ 2 e^{C_2}  -2 }{ C_1C_2} ,\\
\alpha(x, t)=\beta (x, t) = 0, && \textrm{for }\    
t \ge T_1(\Vert (\alpha ^0, \beta ^0, \phi ^0, \phi ^1)\Vert _H)  \textrm{ and }   x\in [0, 1].
\end{eqnarray*}

Using the stability of the origin in $\R^2$ for \eqref{eq:haimo} and the formulas \eqref{VV1}, \eqref{VV2}, we infer the stability of  the origin in $H$ for system \eqref{pdegamma1}, \eqref{eq:haimo} and
\eqref{U1}, which is thus finite-time stable. 
\end{proof}

Introduce the space 
\begin{eqnarray*}
  {\mathcal H}  &:=& 
   \{ (z^0,z^1,\Omega ^0,\Omega ^1) \in H^2(0, 1) \times H^1(0, 1)\times \R \times \R;   \\
&&\  z^0_x (0)=0, \ z^0(1)=\Omega ^0,\\
&&\  \ \mu \beta ^0(1)= \frac{2}{\sqrt{\lambda (1)}} \Omega ^1 -a(1)\lambda (1)\alpha ^0(1) + b(1) \lambda (1)\beta ^0(1) 
+\int_0^1 [(a\lambda )_x \alpha ^0 -(b\lambda )_x \beta ^0] \, \mathrm{d}x, \\
&&\  \ \textrm{ with }\  \mymatrix{\alpha^0 \\ \beta^0} (x)
:= \mymatrix{u^0 \\ v^0} (x)
-\int_0^x {\bf{K}} (x,\xi )     \mymatrix{u^0 \\ v^0} (\xi ) \, \mathrm{d}\xi \quad \forall x\in [0, 1]\\
&&\  \ \textrm{ and } \  
\mymatrix{u^0 \\ v^0} (x) := \frac{1}{\sqrt{\lambda (x)}}  \mymatrix{
z^1(x) -\lambda (x) z^0_x(x) \\[2mm] 
 z^1(x) + \lambda (x) z^0_x (x) 
}
\quad \forall x\in [0, 1]
  \} 
  \end{eqnarray*}
endowed with the norm $\Vert (z^0,z^1,\Omega ^0, \Omega ^1) \Vert _{\mathcal H} ^2 :=\Vert z^0\Vert _{H^2(0, 1)} ^2 
+ \Vert z^1\Vert _{H^1(0, 1)} ^2 + |\Omega ^0|^2+|\Omega ^1|^2.$

We are in a position to state the main result in this paper. 
\begin{theorem}
\label{thm2}
Let $\psi, \zeta\in \R$ with $0<\psi <1$ and $\zeta \ge \frac{\psi}{2-\psi}$. Then 
the system \eqref{eq:zlambda}-\eqref{eq:zlambda2} with the feedback law \eqref{eq:laloi} 
is well-posed and  globally finite-time stable in the Hilbert space $\mathcal H$. 
 More precisely, there exists a nondecreasing 
 function $T:(0, +\infty)\to (0, +\infty)$ such that for all $(z^0, z^1, \Omega ^0, \Omega ^1) \in {\mathcal H}$, 
 the solution of  \eqref{eq:zlambda}-\eqref{eq:zlambda2}, \eqref{eq:laloi} 
and $(z(.,0),z_t(.,0), X_\mathrm{p}(0), \dot X_\mathrm{p} (0) )=(z^0,z^1,\Omega ^0, \Omega ^1)$
satisfies 
\ba
z(x, t)=0,&&  \forall x\in [0, 1],\ \forall t\ge T(\Vert (z^0, z^1,\Omega ^0, \Omega ^1)\Vert _{\mathcal H}), \label{W2} \\
X_\mathrm{p}(t)=0,&& \forall t\ge T(\Vert (z^0,z^1,\Omega ^0, \Omega ^1)\Vert _{\mathcal H}) .
\label{W3}
\ea 
\end{theorem}

\begin{proof}
Let $u^0,v^0,\alpha ^0, \beta ^0\in H^1(0, 1)$ be as in the definition of the space $\mathcal H$. 
Note that $\alpha^0 (0)=\beta ^0 (0)$, for $z^0_x(0)=0$. 
Let 
\begin{eqnarray}
\phi^0 &:=&  \frac{2}{\sqrt{\lambda(1)}} \Omega ^0 + \int_0^1 \left[ a(x) \alpha ^0(x) + b(x) \beta ^0(x) \right] \, \mathrm{d} x,  \label{W4}\\
\phi ^1 &:=& \mu \beta ^0(1). \label{W5}
\end{eqnarray}
We infer from Theorem \ref{thm1}  the existence and uniqueness of a solution $(\alpha, \beta, \phi )$ on $\R_+$
of \eqref{pdegamma1}, \eqref{eq:haimo},  \eqref{U1}, and  $(\alpha (.,0),\beta (.,0), \phi, \dot \phi) =(\alpha ^0, \beta ^0, \phi ^0, \phi ^1)$, 
and this solution satisfies \eqref{U2}-\eqref{U3}. We know also from Theorem \ref{thm1} that $\alpha, \beta \in C^0( [0, T], L^2(0, 1))$ and that  $\phi \in C^2( [0, T] )$ for all $T>0$. \\

Here the initial data are more regular, and therefore  the trajectories are expected to be  more regular. \\

\noindent
{\sc Claim 1.} $\alpha, \beta\in C^0([0, T], H^1(0, 1)) \cap C^0([0, 1], H^1(0, T))$ for all $T>0$.\\[3mm]
Indeed, using the properties $\beta ^0\in H^1(0, 1)$, $\dot \phi \in C^1([0, T])\subset H^1(0, T)$, the compatibility condition  \eqref{W5}
and the formula \eqref{VV1}, we infer that  $\beta \in C^0([0, T], H^1(0, 1)) \cap C^0([0, 1],H^1(0, T))$. In particular,  $\beta ( 0 , . ) \in H^1(0, T)$.
In the same way, using the properties $\alpha ^0\in H^1(0, 1)$, $\beta (0,.)\in H^1(0, T)$, the compatibility condition $\alpha ^0 (0)=\beta ^0(0)$ and 
the formula  \eqref{VV2}, we conclude that $\alpha \in C^0([0, T], H^1(0, 1)) \cap C^0([0, 1],H^1(0, T))$.

Next, we define respectively $X_\mathrm{p}(t)$ as 
\be
X_\mathrm{p}(t):= \frac{\sqrt{\lambda (1)}}{2} \left( \phi (t) -\int_0^1 [a(x) \alpha (x, t) + b(x) \beta (x, t) ]\, \mathrm{d}x \right) ,
\label{W6}
\ee
$\bm{w}(x, t)$ as 
\be
\label{W6bis}
\bm{w}(x, t) = 
\mymatrix{u(x, t) \\ v(x, t)}
:=\bm{\gamma} (x, t) +\int_0^x \bm{L}(x, \xi ) \bm{\gamma} ( \xi , t ) \, \mathrm{d} \xi 
\ee
(with 
$\bm{\gamma} (x, t):= 
\mymatrix{ \alpha (x, t) \ \beta (x, t) }^\intercal $),  and $[ D(x, t)\   S(x, t)] ^T$ as 
\be
\label{W6ter}
\mymatrix{D(x, t) \\ S(x, t)}
:=\sqrt{ \lambda (x) } \, w(x, t).   
\ee
Clearly, we also have that
\be
u, v, D, S \in C^0([0, T], H^1(0, 1)) \cap C^0([0, 1], H^1(0, T))\quad \forall T>0.
\label{W10}
\ee
Furthermore, since $\bm{\gamma}$ satisfies \eqref{pdegamma1}, we obtain that $\bm{w}$ satisfies \eqref{AB3} and that 
the functions $D, S$ satisfy system \eqref{AB1}-\eqref{AB2}. 
We are in a position to define the function $z(x, t)$. \\

\noindent
{\sc Claim 2. }  For every $T>0$, there exists a unique function $z\in C^0([0, T], H^2(0, 1))\cap C^0([0, 1], H^2(0, T))$ 
of the system 
\begin{eqnarray}
z_t &=& \frac{S+D}{2},\quad (x, t)\in (0, 1)\times (0, T),\label{W11}\\
z_x &=& \frac{S-D}{2\lambda },\quad (x, t)\in (0, 1)\times (0, T),\label{W12}\\
z(1,0) &=& z^0(1)=X_\mathrm{p}(0). \label{W13}
\end{eqnarray}
Indeed, setting $f:=(S+D)/2$ and $g=(S-D)/(2\lambda)$, we notice that Schwarz' condition $f_x=g_t$ is satisfied, since by
\eqref{AB1}-\eqref{AB2}
\[
g_t= \frac{S_t-D_t}{2\lambda } = \lambda \frac{S_x+D_x}{2\lambda} + 0=f_x.
\]
On the other hand, the compatibility condition  $z^0(1)=X_\mathrm{p}(0)$ is fulfilled, for 
\[
X_\mathrm{p}(0)=\frac{\sqrt{\lambda (1)}}{2} \left( \phi ^0 -\int_0^1 [a(x) \alpha ^0(x) + b(x) \beta ^0(x)] \,\mathrmd x \right) 
=\Omega ^0 =z^0(1), 
\]
 where we used  \eqref{W6},  \eqref{W4} and the 
property $z^0(1)=\Omega ^0$ from the definition of $\mathcal H$. It follows that there exists a unique solution $z=z(x, t)$ of system \eqref{W11}-\eqref{W13}
which is given explicitly by 
\be
\label{W15}
z(x, t)= X_\mathrm{p} (0) + \int_0^t f(1,s) \, \mathrm{d}s+\int_1^x g(s,t)\, \mathrm{d}s
= z^0(1) + \int_1^x g(s,0)\, \mathrm{d}s + \int_0^t f(x,s)\, \mathrm{d}s. 
\ee
Combined with \eqref{W10}, this yields
\be
\label{W16}
z\in C^0([0, T], H^2(0, 1)) \cap C^1([0, T], H^1(0, 1) )   \cap C^0([0, 1],H^2(0, T)). 
\ee

We claim that $z(x,0)=z^0(x)$ for all $x\in [0,1]$. Indeed, we have by \eqref{W15}, \eqref{W6ter},  and some condition in the definition of $\mathcal H$ that for all $x\in [0,1]$
\begin{eqnarray*}
z(x,0)&=&z^0(1)+\int_1^x g(s,0)\, \mathrmd s\\
&=& z^0(1)+\int_1^x \frac{S(s,0)-D(s,0)}{2\lambda (s)}\, \mathrmd s \\
&=& z^0(1) + \int_1^x \frac{1}{2\sqrt{\lambda (s)}} (v^0(s)-u^0(s))\, \mathrmd s \\
&=& z^0(1)+\int_1^x z^0_x(s)\, \mathrmd s\\
&=& z^0(x)
\end{eqnarray*}

We also claim that 
\be
\label{ABC1} 
z(1,t)=X_\mathrm{p} (t) \qquad \forall t\ge 0. 
\ee
Indeed, we infer from \eqref{W15} and \eqref{W6ter} that 
\[
z(1,t)=\int_0^t f(1,s) \mathrmd s + X_\mathrm{p}(0) = \int_0^t \frac{S(1,s)+D(1,s)}{2}\, \mathrmd s + X_\mathrm{p}(0)=\frac{\sqrt{\lambda (1)}}{2}  \int_0^t (u(1,s)+v(1,s))\, \mathrmd s + X_\mathrm{p}(0)
\]
so that with \eqref{W6bis} 
\[
z_t(1,t) = \frac{\sqrt{\lambda (1)}}{2}   \big( u(1,t)+v(1,t)\big)  =   \frac{\sqrt{\lambda (1)}}{2} \left( \alpha (1,t) + \beta (1,t) + \int_0^1 \mymatrix{1 & 1} \bm{L}(1,\xi ) 
\mymatrix{ \alpha (\xi , t) \\  \beta (\xi , t)} \, \mathrmd\xi 
\right) .
\]
On the other hand, \eqref{W6} gives \eqref{eq:defphi} and hence \eqref{KK2}. Using \eqref{PPP1}-\eqref{PPP2}, we obtain
\begin{multline*}
\dot \phi (t) = \frac{2}{\sqrt{\lambda (1)}} \dot X_\mathrm{p} (t)  \\
+\underbrace{\left( -a_0 + \int_0^1 \big( L^{\alpha\alpha}(1, \xi) + L^{\beta \alpha} (1,\xi )  \big) \mathrmd \xi  \right) }_{=-1} \alpha (1,t) 
+  \underbrace{\left( b_0 + \int_0^1 \big( L^{\alpha\beta}(1, \xi) + L^{\beta \beta} (1,\xi )  \big) \mathrmd\xi  \right) }_{= \mu-1} \beta (1,t) \\ 
 - \int_0^1 [(L^{\alpha\alpha}(1,x) + L^{\beta \alpha} (1,x))  \alpha (x,t) + (L^{\alpha \beta }(1,x)+ L^{\beta\beta} (1,x))\beta (x,t)] \mathrmd x. 
\end{multline*}
Combined with \eqref{U12}, this yields
\begin{eqnarray*}
\dot X_\mathrm{p} (t) 
&=& \frac{\sqrt{\lambda (1)} }{2} 
\left(    \alpha (1,t) + \beta (1,t) + 
\int_0^1 [(L^{\alpha\alpha}(1,x) + L^{\beta \alpha} (1,x))  \alpha (x,t) + (L^{\alpha \beta }(1,x)+ L^{\beta\beta} (1,x))\beta (x,t)] \, \mathrmd x \right) \\
&=& z_t(1,t).
\end{eqnarray*}
Integrating w.r.t $t$ and using $X_\mathrm{p} (0)=z^0(1)=z(1,0)$, we  obtain \eqref{ABC1}.

Let us check that $z$ solves system 
\eqref{eq:zlambda}-\eqref{eq:zlambda3}. 

Replacing $z_t$ and $z_x$ by their expressions in terms of $S,D$ and using \eqref{AB1}-\eqref{AB2}, we obtain
\[
z_{tt} - \lambda ^2 z_{xx}=\frac{S_t+D_t}{2}-\lambda ^2 \left( \frac{S_x-D_x}{2\lambda } -\lambda ' \frac{S-D}{2\lambda ^2}  \right) =0.
\] 

For the boundary conditions \eqref{eq:zlambda2}, we have that $z(1, t)=X_\mathrm{p}(t)$ (by construction of $z$), and that 
\[
z_x(0, t)=\frac{ S(0, t) - D(0, t) }{ 2 \lambda (0) } = \frac{ \sqrt{\lambda (0)} }{2 \lambda (0) }\left( v(0, t) -u(0, t)  \right) =
\frac{1}{2\sqrt{\lambda (0) }} (\beta (0, t)-\alpha (0, t))=   0.
\]
For the initial conditions \eqref{eq:zlambda3}, we have that $z(\cdot ,0)=z^0$ (by construction of $z$) and that 
\[
z_t(x, 0) = \frac{ S(x, 0)+D(x, 0) }{2} =\frac{\sqrt{\lambda (x)} }{2} (u(x, 0)+v(x, 0)) 
=\frac{\sqrt{\lambda (x)} }{2} (u^0(x)+v^0(x))=z^1(x)\quad \forall x\in [0, 1].  
\]
Let us investigate the dynamics of $X_\mathrm{p}$. Set $U(t):=\ddot X_\mathrm{p}(t)$ for $t\in \R_+$.\\

\noindent
{\sc Claim 3.} $U\in L^2(0, T)$ for all $T>0$. \\[3mm]
We infer from  \eqref{eq:haimo} and \eqref{W6} that 
\begin{eqnarray*}
U(t)&=& \ddot X_\mathrm{p}(t)\\
&=& \frac{\sqrt{\lambda (1)}}{2} 
\left( 
\ddot \phi (t) -\int_0^1 [a\alpha _{tt} +b\beta _{tt}]\, \mathrm{d}x
\right) \\
&=&  - \frac{\sqrt{\lambda (1)}}{2}  \left( 
 \lfloor \dot{\phi} (t) \rceil^\psi + \lfloor \phi(t) \rceil^\zeta - \int_0^1 [a\alpha _{tt} +b\beta _{tt}]\, \mathrm{d}x
\right) .  
\end{eqnarray*}
As $\phi \in C^2([0, T])$, the two first terms in the last equation are in $C^0([0, T])$, and hence in $L^2(0, T)$. 
It remains to show that the map $t\to \int_0^1 [a\alpha _{tt} +b\beta _{tt}]\, \mathrm{d}x$ is in $L^2(0, T)$
for all $T>0$. Using \eqref{U11} and \eqref{U21} and next an integration by parts, we obtain
\begin{eqnarray*}
\int_0^1[a\alpha_{tt} + b\beta _{tt}]\, \mathrm{d}x 
&=& \int_0^1 [a(-\lambda \alpha _{xt} )+b\lambda \beta _{xt} ] \, \mathrm{d}x \\
&=& \int_0^1 [(a\lambda )_x \alpha _t - (b\lambda )_x \beta _t]\, \mathrm{d}x  + 
[-a\lambda \alpha _t +b\lambda \beta _t  ]_0^1.
\end{eqnarray*}
Replacing  in the last integral term $\alpha _t $ and $\beta_t$ by $-\lambda \alpha _x$ and $\lambda \beta _x$, respectively, and integrating by parts again, we obtain
\begin{eqnarray}
\int_0^1[a\alpha_{tt} + b\beta _{tt}]\, \mathrm{d}x 
&=& \int_0^1 [ ((a\lambda )_x\lambda )_x \alpha  + ((b\lambda )_x\lambda)_x \beta ]\, \mathrm{d}x   -[(a\lambda )_x\lambda \alpha        +(b\lambda )_x\lambda \beta ]_0^1 
\label{UUU} \\
&&\qquad + [ - a\lambda \alpha _t +b\lambda \beta _t  ]_0^1.\nonumber
\end{eqnarray}
(Note that  $(a\lambda )_x = -L^{\alpha \alpha }(1,x)-L^{\beta\alpha}(1,x)$ and  $(b\lambda)_x=L^{\alpha\beta}(1,x)+L^{\beta\beta}(1,x)$.) 
It is clear that the two first terms in the right hand side of \eqref{UUU} are in $C^0([0,T])$ (and thus in $L^2(0,T)$), for $\alpha , \beta  \in C^0( [0, T], H^1(0, 1))$.
On the other hand, using  \eqref{eq:haimo}, \eqref{U12} and \eqref{U22}, we obtain that 
\begin{eqnarray*}
[-a\lambda \alpha _t +b\lambda \beta _t  ]_0^1 &=& -a(1) \lambda (1) \alpha _t (1, t) + b(1) \lambda (1)\mu ^{-1}  \ddot \phi (t)\\
&=& -a(1) \lambda (1) {\color{black} \alpha _t (1, t)}  +  b(1) \lambda (1)\mu ^{-1} 
(\lfloor \dot{\phi} (t) \rceil^\halpha + \lfloor \phi(t) \rceil^\hbeta). 
\end{eqnarray*}
As $\phi \in C^2([0, T])$ and $\alpha \in C^0([0, 1],H^1(0, T))$ (and hence,  $\alpha _t(1,.)\in L^2(0, T)$), we infer that the map $t\to  [-a\lambda \alpha _t +b\lambda \beta _t  ]_0^1$ is also in $L^2(0, T)$ for all $T>0$. Claim 3 is proved. 
 
Let us have a look at the initial conditions for $X_\mathrm{p}$. By \eqref{PP2}, \eqref{W4} and \eqref{W6}, we have that 
\[ X_\mathrm{p}(0)= 
\frac{\sqrt{\lambda (1)}}{2} \left( \phi ^0 -\int_0^1 [a(x) \alpha ^0(x) + b(x) \beta ^0(x) ]\, \mathrm{d}x \right)  =\Omega ^0.
\] 
On the other hand, \eqref{W6} gives \eqref{eq:defphi} and \eqref{KK2}. Picking $t=0$ and using  \eqref{U12}, we arrive at
\[
\mu \beta (1,0)=\dot \phi (0)  =  \frac{2}{\sqrt{\lambda (1)}} \dot X_\mathrm{p}(0) -a(1)\lambda (1)\alpha ^0(1) + b(1) \lambda (1)\beta ^0(1) 
+\int_0^1 [(a\lambda )_x \alpha ^0 -(b\lambda )_x \beta ^0]\, \mathrm{d}x. 
\]
Comparing with the condition 
\[
\mu \beta ^0(1)= \frac{2}{\sqrt{\lambda (1)}} \Omega ^1 -a(1)\lambda (1)\alpha ^0(1) + b(1) \lambda (1)\beta ^0(1) 
+\int_0^1 [(a\lambda )_x \alpha ^0 -(b\lambda )_x \beta ^0]\, \mathrm{d}x
\]
present in the definition of $\mathcal H$, we infer that $\dot X_\mathrm{p}(0)=\Omega ^1$. 

Using Claim 2, we see  that the uniqueness of $(z, X_\mathrm{p})$ follows from those of 
$(\alpha, \beta, \phi )$.  

We know that $\alpha(x, t)=\beta(x, t)=\phi (t)=0$ for $x\in [0, 1]$ and 
$t\ge T_1(\Vert (\alpha ^0, \beta ^0, \phi ^0, \phi ^1)\Vert _H)$. We infer 
 from \eqref{W6}  that $X_\mathrm{p}(t) = 0$ for $t\ge T_1(\Vert (\alpha ^0, \beta ^0, \phi ^0, \phi ^1)\Vert _H)$.

Finally, from  \eqref{W6bis} and \eqref{W6ter}, we infer that 
\[
u(x, t)=v(x, t)=D(x, t)=S(x, t)=0, \quad \forall x\in [0, 1],\ \forall t\ge T_1(\Vert (\alpha ^0, \beta ^0, \phi ^0, \phi ^1)\Vert _H).
\]
It follows that 
\[
z_t(x, t)=z_x(x, t)=0, \quad \forall x\in [0, 1], \ \forall t\ge T_1(\Vert (\alpha ^0, \beta ^0, \phi ^0, \phi ^1)\Vert _H).
\]
Since $z(1, t)=X_\mathrm{p}(t)=0$ for $t\ge T_1(\Vert (\alpha ^0, \beta ^0, \phi ^0, \phi ^1)\Vert _H)$, 
we arrive at the conclusion that 
\[
z(x, t)=0, \quad \forall x\in [0, 1], \ \forall t\ge T_1(\Vert (\alpha ^0, \beta ^0, \phi ^0, \phi ^1)\Vert _H).
\]
One can pick 
\[
T( R) := \sup \{ T_1(\Vert (\alpha ^0, \beta ^0, \phi ^0, \phi ^1)\Vert _H ); \ 
\Vert (z^0,z^1,\Omega ^0, \Omega ^1)\Vert_{\mathcal H} \le R \}. 
\]
It remains to prove the stability of the origin in $\mathcal H$ for system \eqref{eq:zlambda}-\eqref{eq:zlambda2} and \eqref{eq:laloi}. 
Assume that $\Vert (z^0,z^1,\Omega ^0, \Omega ^1)\Vert_{\mathcal H}$ is small, and (at least)  less than 1. Then $\Vert (\alpha ^0, \beta ^0, \phi ^0, \phi ^1)\Vert _H$ is small, and by 
Theorem \ref{thm1}, $\phi (t)$ and $\dot \phi (t)$ remain small. We know that $z(., t), X_\mathrm{p} (t),\alpha (.,t), \beta (., t)$, and $\phi (t)$ vanish for $t\ge T(1)$. 
Using  \eqref{VV1} and \eqref{VV2}, we see that $\Vert \beta (.,t)\Vert_{H^1(0,1)}$ and $\Vert \alpha (.,t)\Vert _{H^1(0,1)}$ remain small. 
Using \eqref{W6} and \eqref{KK2}, we conclude that $X_\mathrm{p} (t) $ and $\dot X_\mathrm{p} (t)$ also remain small. 
Clearly, $\Vert D(.,t)\Vert_{H^1(0,1)}$ and $\Vert S(.,t)\Vert _{H^1(0,1)}$ also remain small by \eqref{W6bis}-\eqref{W6ter}. 
This yields that $\Vert z_t (.,t)\Vert _{H^1(0,1)}$ remains small, by \eqref{W11}.  Using \eqref{W12} and \eqref{ABC1}, we infer that
$\Vert z(.,t)\Vert _{H^2(0,1)}$ remains small.  
\end{proof}


\section{Simulation}\label{sec:simu}
As numerical illustration, system \eqref{A1}-\eqref{A4} controlled by the feedback law \eqref{eq:laloi} is simulated. The following system parameters are used: 
$m=2\, \mathrm{kg},\  \rho = 2\, \mathrm{kg}/\mathrm{m}, \ g=9.81\, \mathrm{m}/\mathrm{s}^2$. 
\subsection{Numerical calculation of the kernels}

\begin{figure}[thpb]
	\centering
	\input{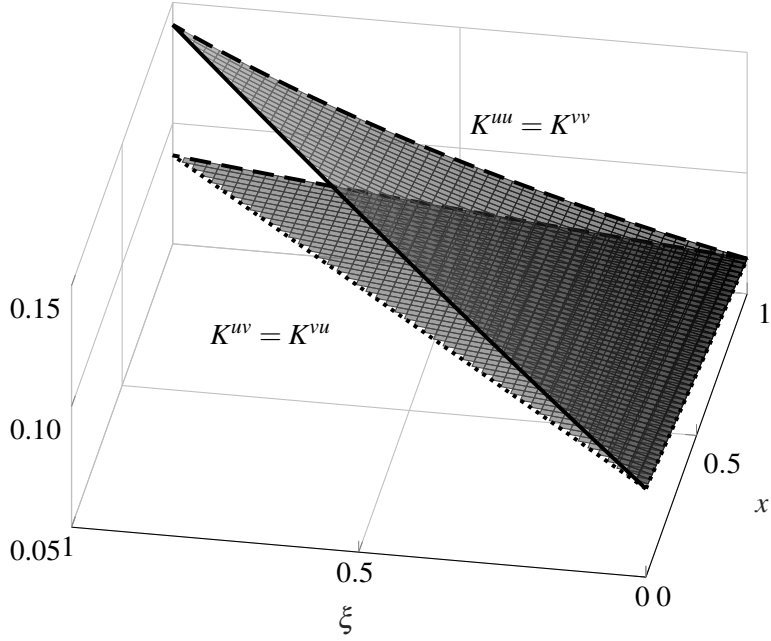}
	\caption{Direct kernels $\bm{K}$. Dotted boundaries correspond to boundary conditions at $(x, 0)$ and $(x, x)$ of the Goursat system \eqref{eq:goursatone}. Dashed boundaries correspond to controller gains at $(1, \xi)$ used in the control law \eqref{eq:dummycoronlaw}.}
	\label{fig:kernelsK}
\end{figure}

\begin{figure}[thpb]
	\centering
	\input{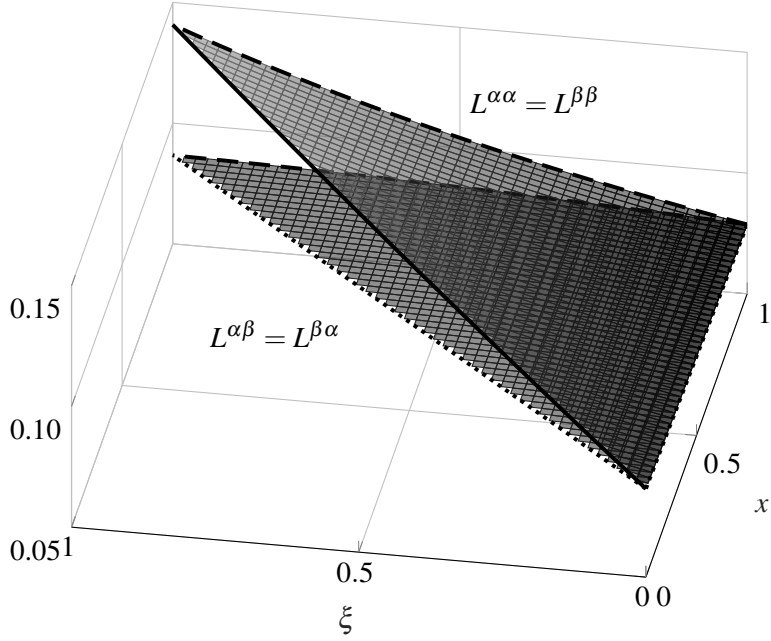}
	\caption{Inverse kernels $\bm{L}$. Dotted boundaries correspond to boundary conditions at $(x, 0)$ and $(x, x)$ of the Goursat system. Dashed boundaries correspond to controller gains at $(1, \xi)$ used in the equations \eqref{PPP1}, \eqref{PPP2}, \eqref{eq:defmu}.}
	\label{fig:kernelsL}
\end{figure}

\begin{figure}[thpb]
	\centering
	\input{pics/kernelL1article.tex}
	\caption{Controller gains $\bm{L}(1, \xi)$.}
	\label{fig:kernelL1}
\end{figure}
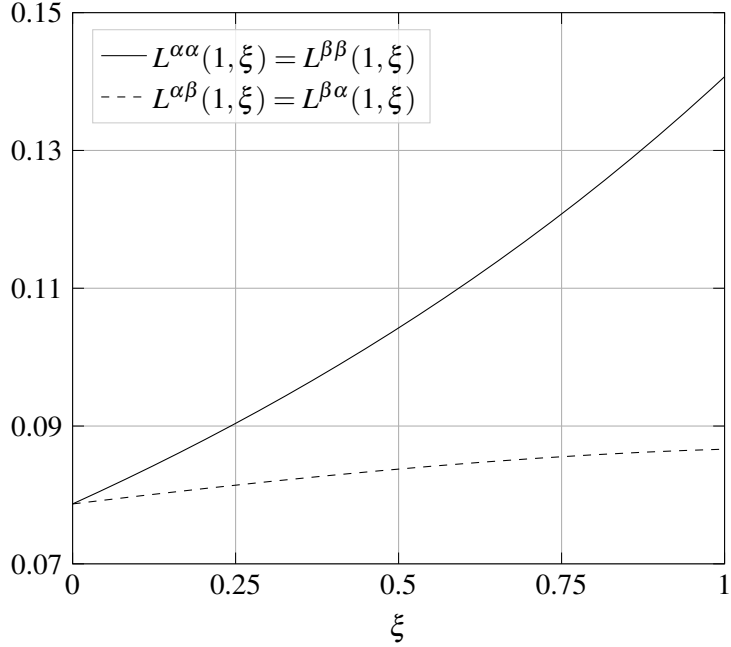

The direct kernels $\bm{K}$ are calculated by numerically solving system \eqref{eq:goursatone} \cite{AA}, where the domain $\mathcal{T}$ is discretized using a uniform grid with $\Delta x = \Delta \xi = 0.005$ (see Fig. \ref{fig:tdom}). Knowing the direct kernels $\bm{K}$, the inverse kernels $\bm{L}$ are calculated numerically using \eqref{eq:relKL}. The obtained kernels are shown in Figures \ref{fig:kernelsK}-\ref{fig:kernelsL}. The used controller gains correspond to the kernels $\bm{L}$ evaluated at the boundary $(x=1, \xi)$, that are shown in detail in Fig \ref{fig:kernelL1}. 

\subsection{Time evolution of the controlled system}

Once the direct and inverse kernels are known, the system \eqref{A1}-\eqref{A4} controlled by the law \eqref{eq:laloi} can be simulated in the coordinates $(X_\mathrm{p}(t), \alpha(x, t), \beta(x, t))$. Finally, the inverse state transformation will be applied in order to obtain the simulated time evolution in the original coordinates $(X_\mathrm{p}(t), y(s))$.

We consider zero initial conditions, except for $X_\mathrm{p}(0) = 0.5$, i.e., the system is initially at rest with the platform at a nonzero position.

\subsubsection{Time evolution of $\phi(t)$}

The finite-time dynamics of the variable $\phi(t)$ \eqref{eq:defphi} is simulated using the numerical method proposed in \cite{PEB}. This method supposes that the system's vector field $\bm{F}$ is $\bm{d}$-homogeneous \cite[Def. 3.6]{PEB}, which corresponds to imposing that $\hbeta=\frac{\halpha}{2-\halpha}$ in \eqref{eq:haimo}. A nonlinear state transformation is provided with which an alternative continuous time system representation $\dot{\bm{z}}=\tilde{\bm{F}}(\bm{z})$ is obtained, that admits an implicit discretization scheme preserving the finite-time stability property in discrete time. Performing the inverse state transformation then yields a discretization scheme 
preserving the finite-time stability in the original coordinates.

Rewrite the finite-time ODE dynamics \eqref{eq:haimo} as
\myeq{ \frac{\mathrm{d}}{\mathrm{d}t} \mymatrix{\phi \\ \dot{\phi}} = \mymatrix{\dot{\phi} \\ -\lfloor \phi \rceil^{\frac{\halpha}{2-\halpha}} -  \lfloor \dot{\phi} \rceil^\halpha} =:               \bm{F}\left(\mymatrix{\phi \\ \dot \phi} \right). }{defFfield}

Let us check that the assumptions in \cite[Theorem 4.1]{PEB} are fulfilled. 
We have that the vector field $\bm{F}(\bm{x})$, with $\bm{x}:=\mymatrix{x_1 & x_2}^\intercal$,
is uniformly continuous on the unit sphere, $\bm{d}$-homogeneous with homogeneity degree $\nu_{\bm{d}}=-1$ for the weighted dilation $\bm{d}(s)=\mymatrix{e^{r_1s} & 0 \\ 0 & e^{r_2s}}$ where $\mymatrix{r_1 & r_2} = \mymatrix{\frac{2-\halpha}{1-\halpha} & \frac{1}{1-\halpha}}$, and $\bm{F}(-\bm{x}) = -\bm{F}(\bm{x})$.
Let $\bm{G_d}=\text{diag}(r_1, r_2)$ be the generator of the dilatation $\bm{d}(s)$ and take the symmetric matrix $\bm{P} := \bm{I}_2$ that satisfies $\bm{PG_d}+\bm{G_d}^\intercal\bm{P}\succ 0$.

Lastly, it is required that the condition
\myeq{  \bm{z}^\intercal \bm{\Xi}^\intercal(\bm{z}) \bm{P} \bm{\Xi} (\bm{z}) \left[\frac{(\bm{I}_2-\bm{G_d}) \bm{z} \bm{z}^\intercal \bm{P}}{\bm{z}^\intercal\bm{P}\bm{G_d}\bm{z}} + \bm{I}_2 \right] \bm{F} \left( \frac{\bm{z}}{ \lVert \bm{z} \rVert} \right) < 0 }{condPEB}
be satisfied for all $\bm{z}\in\mathbb{R}^2\backslash\{\bm{0}\}$, where $\bm{\Xi} = \bm{I}_2$ and $\lVert \bm{z} \rVert = \sqrt{\bm{z}^\intercal \bm{P} \bm{z}}$, which according to \cite[Thm. 3.8]{PEB}-\cite[Thm. 4]{polyakov2019sliding} guarantees that the quadratic form $V(\bm{z}) = \bm{z}^\intercal \bm{P} \bm{z} = \Vert \bm{z}\Vert ^2 $ is a Lyapunov function for the system
\myeq{ \dot{\bm{z}} = \left( \frac{\left(\bm{I}_2 - \bm{G_d}\right) \bm{z}\bm{z}^\intercal\bm{P}}{\bm{z}^\intercal\bm{PG_d z}} + \bm{I}_2 \right) \bm{F} \left( \frac{\bm{z}}{\lVert \bm{z} \rVert} \right) := \tilde{\bm{F}}(\bm{z}), }{thmpolyakov2019sliding}
which is the transformed version of the original system $\dot{\bm{x}}=\bm{F}(\bm{x})$ using the nonlinear state transformation \eqref{eq:thm41statetf1}-\eqref{eq:thm41statetf2} (see below).

Condition \eqref{eq:condPEB} simplifies to
\myeq{ \frac{\Vert \bm{z} \Vert ^3}{r_1z_1^2 + r_2z_2^2} 
\left[ \frac{z_1z_2}{\Vert \bm{z}\Vert ^2}  - \frac{z_2}{\Vert \bm{z}\Vert } \left\lfloor \frac{z_1}{\Vert \bm{z}\Vert } 
\right\rceil^\zeta -  \frac{z_2}{\Vert \bm{z}\Vert} \left\lfloor \frac{z_2}{\Vert \bm{z}\Vert } \right\rceil^\psi\right] <0. }{strictcond2}
We note that the strict inequality \eqref{eq:strictcond2} is not satisfied for $z_2=0$, since the left-hand expression equals $0$.
One notes that for $z_2=0=\dot{z}_2$, \eqref{eq:thmpolyakov2019sliding} reduces to
\[ \mymatrix{\dot z_1 \\ 0} = \mymatrix{0 \\ -\left \lfloor \text{sign}(z_1) \right\rceil^{\frac{\psi}{2-\psi}}}, \]
which implies that $z_1=\dot{z}_1=0$. Then, the asymptotic stability of the origin for \eqref{eq:thmpolyakov2019sliding}  is obtained by LaSalle's invariance theorem if we prove the strict inequality \eqref{eq:strictcond2} for the case $z_2\neq 0$. It suffices to consider the strict negativity of the expression between brackets.

\begin{lemma}
Let $\psi\in (0,1)$ and $\zeta= \psi/(2-\psi)\in (0,1)$. Then for all $\bm{z}=(z_1, z_2)\in \R \times \R ^*$, it holds
\[
\frac{z_2}{\Vert \bm{z}\Vert}    \left[  \frac{z_1}{\Vert \bm{z}\Vert } -  \left\lfloor \frac{z_1}{\Vert \bm{z}\Vert }  \right\rceil ^\zeta  -    \left\lfloor \frac{z_2}{\Vert \bm{z}\Vert}  \right\rceil ^\psi      \right]  <0.
\]
\end{lemma} 
\begin{proof}
Let
\[
f(\bm{z}) = \frac{z_2}{\Vert \bm{z}\Vert}    \left[  \frac{z_1}{\Vert \bm{z}\Vert } -  \left\lfloor \frac{z_1}{\Vert \bm{z}\Vert }  \right\rceil ^\zeta  -    \left\lfloor \frac{z_2}{\Vert \bm{z}\Vert}  \right\rceil ^\psi      \right], 
\quad \bm{z}\in \R \times \R ^*.  
\]
Note first that $f(0,z_2)=-\left\vert \frac{z_2}{\Vert \bm{z}\Vert } \right\vert ^{1+\psi} <0$ for $z_2\ne 0$. Thus we can assume that $z_1\ne 0$, as well. Let $\bm{y}:=(y_1,y_2)=
\bm{z}/\Vert \bm{z}\Vert$. Then $y_1,y_2\in (-1,0)\cup (0,1)$, $\Vert \bm{y}\Vert =1$ and 
$f(\bm{z})=g(\bm{y})=  y_2\left( y_1 -  \left\lfloor y_1  \right\rceil ^\zeta -\left\lfloor y_2 \right\rceil ^\psi \right)$. 
From $\Vert \bm{y}\Vert =1$, we infer that $y_2=\pm \sqrt{1-y_1^2}$. We have to consider four cases:
\begin{enumerate}[label=(\roman*)]
\item $y_1\in (0,1)$ and $y_2=\sqrt{1-y_1^2}$;
\item $y_1\in (0,1)$ and $y_2=-\sqrt{1-y_1^2}$;
\item $y_1\in (-1,0)$ and $y_2=\sqrt{1-y_1^2}$;
\item $y_1\in (-1,0)$ and $y_2=-\sqrt{1-y_1^2}$.
\end{enumerate}

In case (i), we have $g(\bm{y})<0$, for $ y_1 -  \left\lfloor y_1  \right\rceil ^\zeta <0$. In case (iv), we also have that $g(\bm{y})<0$, for 
\[
y_1 -  \left\lfloor y_1  \right\rceil ^\zeta -\left\lfloor y_2 \right\rceil ^\psi  = - ( |  y_1 |  -  | y_1 |  ^\zeta ) + |  y_2 |  ^\psi  >0. 
\]
Let $h(s)=-s +s^\zeta -(1-s^2)^\frac{\psi}{2}$ for $s\in (0,1)$. Then in case (ii), we have that $g(\bm{y})=\sqrt{1-y_1^2} h(y_1)$, while in case (iii) we have that 
\[
g(\bm{y})=\sqrt{1-y_1^2} \left(-|y_1|+|y_1|^\zeta - \left(1-y_1^2\right)^\frac{\psi}{2}\right) 
= \sqrt{1-y_1^2}\, h(|y_1 | ).  
\]
Thus, it remains to show that $h(s)<0$ for $s\in (0,1)$. Denoting $\sigma:=s^\zeta \in (0,1)$, we note that 
\begin{eqnarray}
h(s)<0 &\iff& -\sigma ^\frac{1}{\zeta} + \sigma < \left(1-\sigma ^\frac{2}{\zeta} \right) ^\frac{\psi }{2} \nonumber \\
&\iff& \sigma\left( 1- \sigma ^{ \frac{1}{\zeta} -1} \right) < \left( 1+ \sigma ^\frac{1}{\zeta}\right) ^\frac{\psi}{2} \left(1-\sigma ^\frac{1}{\zeta}\right) ^\frac{\psi}{2} \nonumber \\
&\iff& \sigma  ^\frac{2}{\psi}  \left( 1- \sigma ^{ \frac{1}{\zeta} -1} \right) ^ \frac{2}{\psi} < \left( 1+ \sigma ^\frac{1}{\zeta}\right) \left(1-\sigma ^\frac{1}{\zeta}\right). \label{toto}
\end{eqnarray} 
But for $\sigma\in (0,1)$, we have that $ \left( 1- \sigma ^{ \frac{1}{\zeta} -1} \right) ^ \frac{2}{\psi}  <   1- \sigma ^{ \frac{1}{\zeta} -1} <   1- \sigma ^ \frac{1}{\zeta} $ (since $2 / \psi >2$ and 
$1/ \zeta >1$)
and 
$\sigma ^\frac{2}{\psi} < 1< 1+ \sigma ^\frac{1}{\zeta}$, so that \eqref{toto} holds true. 
\end{proof}

Applying \cite[Theorem 4.1]{PEB} to the two-dimensional vector field \eqref{eq:defFfield}, we define the nonlinear state transformation between the 
original state $\bm{x}=\mymatrix{x_1 & x_2}^\intercal:=\mymatrix{\phi & \dot{\phi}}^\intercal$ and a new state $\bm{z}=\mymatrix{z_1 & z_2}^\intercal$:
\myeq{ \bm{z} = \bm{\Phi}(\bm{x}) = \lVert \bm{x} \rVert_{\bm{d}} \bm{d}\left(-\ln \lVert \bm{x} \rVert_{\bm{d}} \right) \bm{x}, }{thm41statetf1}
with inverse transformation
\myeq{ \bm{x} = \bm{\Phi}^{-1}(\bm{z})=\bm{d}\left(\ln\lVert \bm{z} \rVert\right)\frac{\bm{z}}{\lVert\bm{z}\rVert}. }{thm41statetf2}
Here, the canonical homogeneous "norm" \cite[Eq. 3.4]{PEB} is defined as
\[ \lVert \bm{x} \rVert_{\bm{d}} = e^{s_x} \]
where $s_x\in\mathbb{R}$ satisfies
\myeq{\lVert \bm{d}(-s_x)\bm{x} \rVert = 1.}{cond_sx}
More explicitly, \eqref{eq:thm41statetf1} corresponds to
\[ \mymatrix{z_1 \\ z_2} = \begin{cases} \mymatrix{e^{\frac{-s_x}{1-\nu_2}}x_1 \\ e^{\frac{-\nu_2 s_x}{1-\nu_2}}x_2} &\text{if\ } (x_1, x_2) \neq (0, 0) \\ \mymatrix{0 \\ 0} &\text{else},\end{cases} \]
where $s_x$ is the solution of the implicit condition \eqref{eq:cond_sx}
\[ \sqrt{e^{-\frac{2(2-\nu_2)s_x}{1-\nu_2}}x_1^2 + e^{-\frac{2s_x}{1-\nu_2}}x_2^2} = 1 \]
that is solved numerically \cite{more}.

Finally, the obtained implicit numerical scheme in the transformed coordinates $\bm{z}$ is
\[ \frac{\bm{z}(t+\Delta t) - \bm{z}(t)}{\Delta t} = \bm{\tilde{F}}(\bm{z}(t+\Delta t)), \]
where
\[ \bm{\tilde{F}}(\bm{z}) = \left( \frac{\left(\bm{I}_2 - \bm{G_d}\right)\bm{zz}^\intercal}{\bm{z}^\intercal \bm{G_d z}} + \bm{I}_2 \right) \bm{F} \left( \frac{\bm{z}}{\sqrt{\bm{z}^\intercal \bm{z}}} \right). \]

\begin{figure}[thpb]
	\centering
	\input{pics/simphiarticle.tex}\\
	\hspace{.192cm}\input{pics/simphiarticle_log.tex}
	\caption{Simulated time evolution of $\phi(t)$ and $\dot{\phi}(t)$.}
	\label{fig:simphi}
\end{figure}
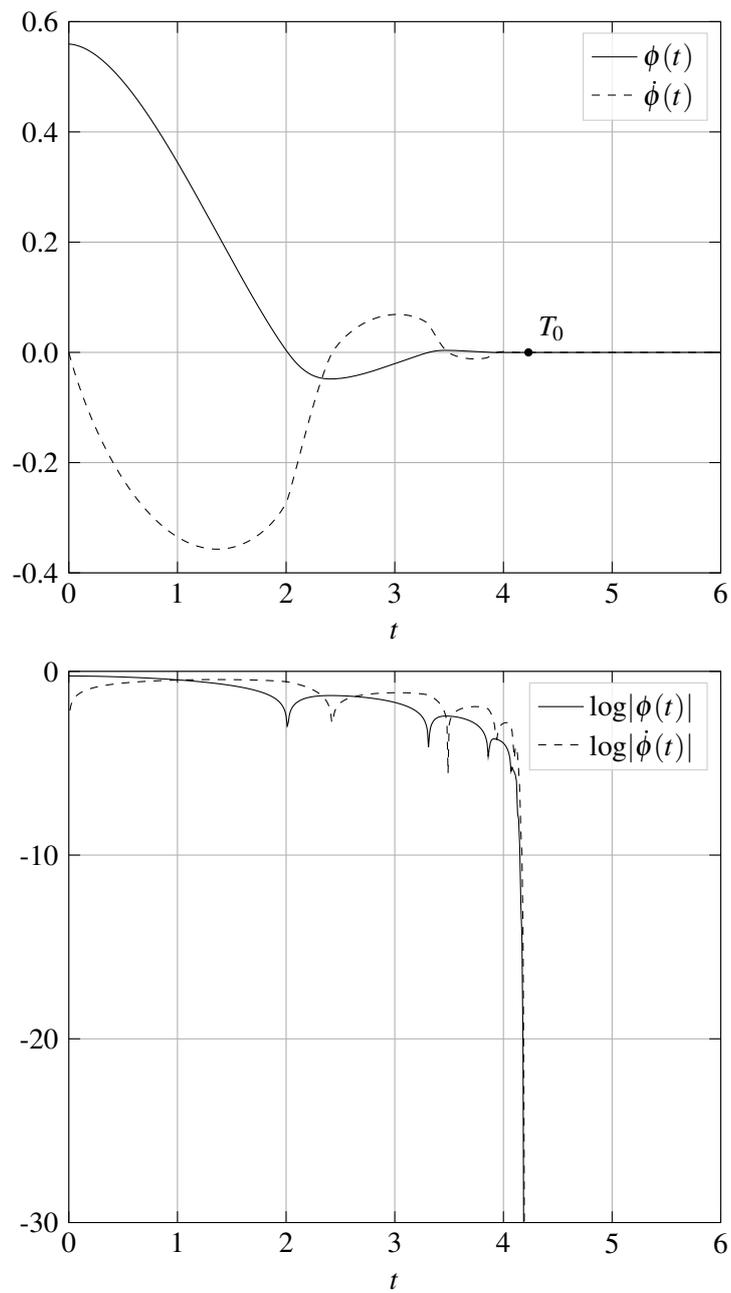

Using this method, the time evolution of $\phi(t)$ with control parameters $\psi=\frac{1}{2}$, $\zeta=\frac{1}{3}$ is simulated for a time step $\Delta t = 0.01$ (Fig. \ref{fig:simphi}). We observe that it is stabilized in a finite time $T_0\approx 4.23$.

\subsubsection{Time evolution of $\alpha(x, t)$ and $\beta(x, t)$}

The two first order hyperbolic PDEs for $\beta(x, t)$, \eqref{U11}, and $\alpha(x, t)$, \eqref{U21}, are solved numerically using a first order downwind-upwind scheme \cite{Sainsaulieu}. Uniform spatial $(\Delta x=0.05)$ and temporal $(\Delta t = 0.01)$ discretization steps are chosen, respecting the CFL stability condition $\text{max}(\lambda(x))\frac{\Delta t}{\Delta x} \leq 1$.

At each time $t$, the following steps are executed. First the boundary condition $\beta(1, t)$ is evaluated from  \eqref{eq:phidotmubeta} knowing the value of $\dot \phi(t)$, where $\mu=2.379$ is calculated from \eqref{eq:defmu} by the trapezoidal rule on the spatial grid. Then, the first order downwind scheme
\[ \beta(x, t) = \beta(x, t-\Delta t) + \lambda(x) \frac{\Delta t}{\Delta x} \left( \beta(x+\Delta x, t-\Delta t) - \beta(x, t-\Delta t) \right) \]
translates the information for $x$ from $1$ to $0$ on the spatial grid. Next, the boundary condition $\alpha(0, t) = \beta(0, t)$ is evaluated. Finally, the first order upwind scheme
\[ \alpha(x, t) = \alpha(x, t-\Delta t) - \lambda(x) \frac{\Delta t}{\Delta x} \left(\alpha(x, t-\Delta t) - \alpha(x-\Delta x, t-\Delta t) \right) \]
translates the information for $x$ from $0$ to $1$ on the spatial grid.

\begin{figure}[thpb]
	\centering
	\input{pics/simalphabetaarticle.tex}\\
	\hspace{.38cm}\input{pics/simalphabetaarticle_log.tex}
	\caption{Simulated time evolution of $\beta(x, t)$ and $\alpha(x, t)$, plotted for $x=0$ and $x=1$.}
	\label{fig:simalphabeta}
\end{figure}
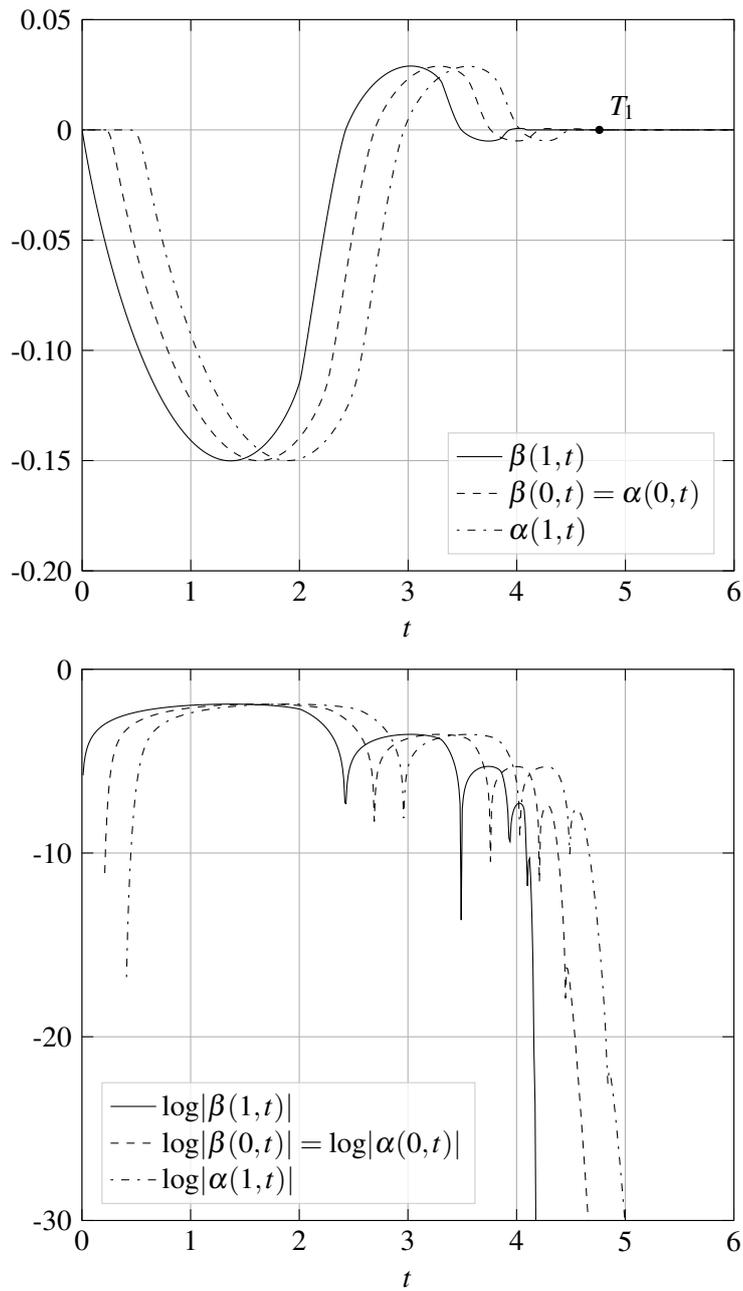

The simulated time evolution of $\beta(x, t)$ and $\alpha(x, t)$ is shown in Figure \ref{fig:simalphabeta}. We observe that they are stabilized in a finite time whose value can be verified using the expression \[T_1=T_0+2\int_0^1 \frac{1}{\lambda(x)}\, \mathrm{d} x = T_0 + 4\sqrt{\frac{m}{\rho g}} \left( e^{\frac{gJ}{2}} - 1\right) \approx 4.76. \]

\subsubsection{Time evolution of $X_\text{p}(t)$}
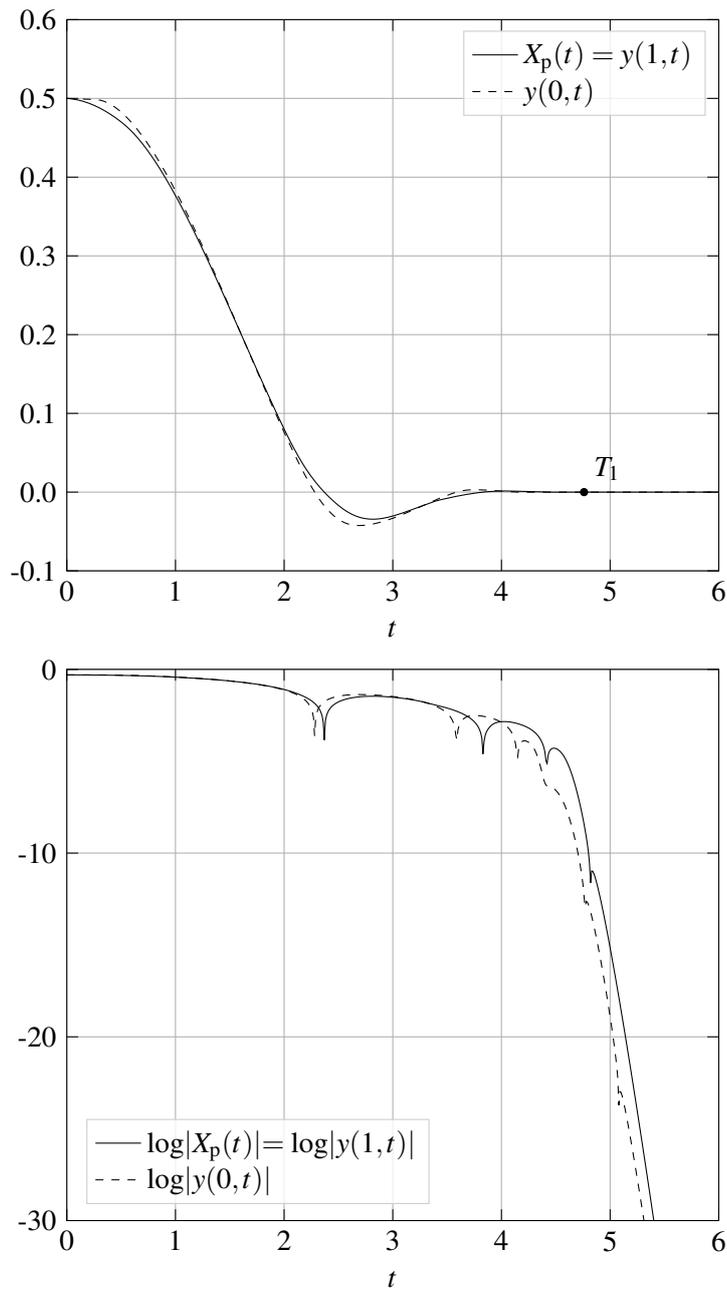
\begin{figure}[thpb]
	\centering
	\input{pics/simXparticle.tex}\\
	\hspace{.192cm}\input{pics/simXparticle_log.tex}
	\caption{Simulated time evolution of the platform position $X_\text{p}(t)$ and cable position $y(s, t)$ at $s=0$ and $s=1$.}
	\label{fig:simXp}
\end{figure}

At each discrete time $t$, $\phi(t)$, $\alpha(x, t)$ and $\beta(x, t)$ are known on a grid for $x$. The corresponding numerical value for $X_\text{p}(t)$ is computed from the definition for $\phi(t)$ (see 
   \eqref{eq:defphi}), numerically evaluating the integral by the trapezoidal rule on the spatial grid. The simulated time evolution of $X_\text{p}(t)$ is shown in Figure \ref{fig:simXp}.

\subsubsection{Time evolution of $y(s, t)$}

Given the simulated time evolutions of $(\alpha(x, t), \beta(x, t))$, one can perform the inverse transformations in order to express the movement of the cable in the original coordinates $y(s, t)$:
\[ \mymatrix{\alpha \\ \beta}(x, t) \xrightarrow{\eqref{eq:backsteppingtfinv}} \mymatrix{u \\ v}(x, t) \xrightarrow{\eqref{eq:wtf}} \frac{1}{\sqrt{\lambda(x)}} \mymatrix{D(x, t) \\ S(x, t)} \\ \xrightarrow[\text{combination}]{\text{linear}} z_x(x, t) \xrightarrow{\int \cdot \, \mathrm{d} x} z(x, t) \xrightarrow{\eqref{eq:tfyz}} y(s, t). \]

For the numerical integration of $z_x(x, t)$ at a given time $t$, $z(1, t)$ is set to $X_\text{p} (t)$ in order to satisfy the compatibility condition corresponding to \eqref{A3}. The spatial integration is then evaluated numerically using the trapezoidal rule for $x$ from $1$ to $0$ on the spatial grid.

The obtained time evolution of $y(s, t)$ is shown in Figure \ref{fig:simXp}.
In conclusion, both the platform $X_\mathrm{p}(t)$ and  the cable $y(s,t)$ have been stabilized at the origin in a finite time, namely $T_1\approx 4.76$.

\section{Conclusions and future works} \label{sec:confut}


A finite-time controller for the motion of an overhead crane described by a hybrid PDE-ODE system with varying tension along the cable is derived. Our feedback law incorporates existing results for the finite-time stabilization of $2 \times 2$ quasilinear hyperbolic systems, that involve kernels that are defined by a Goursat-type system of PDEs. Recent numerical methods are used to calculate the kernels and a simulation of the finite-time stabilization of the overhead crane is shown.


The computation of the proposed control law $u(t)$ \eqref{eq:laloi} needs the knowledge of the position $y(s, t)$ of the cable on its entire length. An observer for this position would be needed if one wants to obtain a more realistic implementation using only the measurement of the position $X_\mathrm{p}(t)$ of the platform and the angle $\theta(t)$ at the top of the cable.

Furthermore, the proposed crane model could be extended to the case with a variable length cable, and with an additional degree of freedom for the bench on which the platform moves (a perpendicular translation of the bench in the case of an overhead crane, or a rotation of the bench in the case of a tower crane).

\section*{Acknowledgments}


The authors thank Henrik Anfinsen for useful discussions related to the numerical calculation of the kernels as presented in \cite{AA}, and for providing code for a numerical implementation of Eq. \eqref{eq:relKL}.

The authors were supported by ANR project Finite4SoS (ANR 15 CE230007). 





\end{document}

%% file: pics/kernelL1article.tex
\begin{tikzpicture}

\begin{axis}[
width=\plottimewidth,
height=\plottimeheight,
legend cell align={left},
legend style={fill opacity=0.8, draw opacity=1, text opacity=1, at={(0.03,0.97)}, anchor=north west, draw=white!80!black},
tick align=inside,
tick pos=both,
x grid style={white!69.0196078431373!black},
xlabel={$\xi$},
xmajorgrids,
xmin=0, xmax=1,
xtick style={color=black},
xtick={0,0.25,0.5, 0.75, 1},
xticklabels={0,0.25,0.5, 0.75, 1},
y grid style={white!69.0196078431373!black},
ymajorgrids,
ymin=0.07, ymax=0.15,
ytick style={color=black},
ytick={0.07,0.09,0.11,0.13,0.15},
yticklabels={0.07,0.09,0.11,0.13,0.15}
]
\addplot [black]
table {%
0 0.0786961317062378
0.0408163070678711 0.0804749727249146
0.0816326141357422 0.0823019742965698
0.122448921203613 0.0841789245605469
0.163265347480774 0.0861074924468994
0.204081654548645 0.0880893468856812
0.244897961616516 0.0901265144348145
0.285714268684387 0.0922207832336426
0.326530575752258 0.0943740606307983
0.367347002029419 0.0965886116027832
0.40816330909729 0.0988665819168091
0.448979616165161 0.101210236549377
0.489795923233032 0.103621959686279
0.530612230300903 0.106104373931885
0.571428537368774 0.108659982681274
0.612244844436646 0.111291646957397
0.653061151504517 0.114002108573914
0.693877577781677 0.116794586181641
0.734693884849548 0.119672179222107
0.775510191917419 0.12263822555542
0.81632661819458 0.125696420669556
0.857142925262451 0.128850340843201
0.897959232330322 0.1321040391922
0.938775539398193 0.135461568832397
0.979591846466064 0.138927102088928
1 0.140703201293945
};
\addlegendentry{$L^{\alpha\alpha}(1, \xi)=L^{\beta\beta}(1, \xi)$}
\addplot [black, dashed]
table {%
0 0.0786961317062378
0.102040767669678 0.0798498392105103
0.204081654548645 0.0809489488601685
0.306122422218323 0.0819859504699707
0.387755155563354 0.0827653408050537
0.469387769699097 0.083495020866394
0.551020383834839 0.0841699838638306
0.632652997970581 0.084784984588623
0.714285731315613 0.0853341817855835
0.795918345451355 0.0858112573623657
0.877551078796387 0.0862096548080444
0.938775539398193 0.0864510536193848
1 0.0866434574127197
};
\addlegendentry{$L^{\alpha\beta}(1, \xi)=L^{\beta\alpha}(1, \xi)$}
\end{axis}

\end{tikzpicture}

%% file: pics/simphiarticle.tex
\begin{tikzpicture}

\begin{axis}[
width=\plottimewidth,
height=\plottimeheight,
legend cell align={left},
legend style={fill opacity=0.8, draw opacity=1, text opacity=1, draw=white!80!black},
tick align=inside,
tick pos=both,
x grid style={white!69.0196078431373!black},
xlabel={$t$},
xmajorgrids,
xmin=0, xmax=6,
xtick style={color=black},
y grid style={white!69.0196078431373!black},
ymajorgrids,
ymin=-0.4, ymax=0.6,
ytick style={color=black},
ytick={-0.4,-0.2,0,0.2,0.4,0.6},
yticklabels={-0.4,-0.2,0.0,0.2,0.4,0.6}
]
\addplot [black]
table {%
0 0.559439301490784
0.0299999713897705 0.55882716178894
0.059999942779541 0.557654976844788
0.100000023841858 0.555276393890381
0.139999985694885 0.552027225494385
0.180000066757202 0.547964811325073
0.220000028610229 0.543138742446899
0.259999990463257 0.537593603134155
0.299999952316284 0.531369805335999
0.340000033378601 0.524504661560059
0.379999995231628 0.517032980918884
0.430000066757202 0.506889820098877
0.480000019073486 0.495909690856934
0.529999971389771 0.48414933681488
0.579999923706055 0.471662282943726
0.629999995231628 0.458499908447266
0.680000066757202 0.444711685180664
0.740000009536743 0.427406311035156
0.799999952316284 0.409347534179688
0.860000014305115 0.390612363815308
0.930000066757202 0.368000030517578
1 0.344685554504395
1.08000004291534 0.317326784133911
1.16999995708466 0.285838723182678
1.28999996185303 0.243081092834473
1.58000004291534 0.139354944229126
1.66999995708466 0.108033776283264
1.75 0.0809376239776611
1.82000005245209 0.0579944849014282
1.87999999523163 0.0390450954437256
1.92999994754791 0.0238687992095947
1.98000001907349 0.00937008857727051
2.01999998092651 -0.00162076950073242
2.04999995231628 -0.00927770137786865
2.07999992370605 -0.0162771940231323
2.10999989509583 -0.0225789546966553
2.14000010490417 -0.0281517505645752
2.17000007629395 -0.0329747200012207
2.20000004768372 -0.0370466709136963
2.23000001907349 -0.0403925180435181
2.25999999046326 -0.0430598258972168
2.28999996185303 -0.0451062917709351
2.3199999332428 -0.0465902090072632
2.35999989509583 -0.0477883815765381
2.40000009536743 -0.0482139587402344
2.45000004768372 -0.0479804277420044
2.50999999046326 -0.0469305515289307
2.5699999332428 -0.0451607704162598
2.64000010490417 -0.04234778881073
2.72000002861023 -0.0383374691009521
2.80999994277954 -0.0330511331558228
2.92000007629395 -0.0258355140686035
3.25999999046326 -0.00287675857543945
3.32999992370605 0.000912785530090332
3.36999988555908 0.00245356559753418
3.41000008583069 0.00332820415496826
3.46000003814697 0.00368201732635498
3.54999995231628 0.00341260433197021
3.67000007629395 0.00227391719818115
3.91000008583069 -0.000211954116821289
6 0
};
\addlegendentry{$\phi(t)$}
\addplot [black, dashed]
table {%
0 0
0.0199999809265137 -0.0144171714782715
0.0399999618530273 -0.0277656316757202
0.059999942779541 -0.0403707027435303
0.0900000333786011 -0.0581879615783691
0.120000004768372 -0.0749301910400391
0.149999976158142 -0.0907632112503052
0.180000066757202 -0.105798959732056
0.210000038146973 -0.120118975639343
0.240000009536743 -0.133785963058472
0.269999980926514 -0.146849274635315
0.309999942779541 -0.163397192955017
0.350000023841858 -0.179021000862122
0.389999985694885 -0.193783164024353
0.430000066757202 -0.207734942436218
0.470000028610229 -0.220919013023376
0.509999990463257 -0.233371496200562
0.549999952316284 -0.245123147964478
0.589999914169312 -0.256200551986694
0.629999995231628 -0.266626834869385
0.670000076293945 -0.276422023773193
0.710000038146973 -0.285603642463684
0.75 -0.294186949729919
0.789999961853027 -0.302185535430908
0.829999923706055 -0.309611082077026
0.879999995231628 -0.318102836608887
0.930000066757202 -0.325732469558716
0.980000019073486 -0.332513928413391
1.02999997138977 -0.338458657264709
1.08000004291534 -0.343575358390808
1.12999999523163 -0.347870230674744
1.17999994754791 -0.351346969604492
1.23000001907349 -0.35400652885437
1.27999997138977 -0.355847358703613
1.33000004291534 -0.356864809989929
1.37999999523163 -0.357051134109497
1.42999994754791 -0.356394529342651
1.48000001907349 -0.354879021644592
1.51999998092651 -0.353033781051636
1.55999994277954 -0.350610971450806
1.60000002384186 -0.347592949867249
1.63999998569489 -0.343957781791687
1.67999994754791 -0.339677810668945
1.72000002861023 -0.334718346595764
1.75999999046326 -0.329033970832825
1.79999995231628 -0.322564840316772
1.8400000333786 -0.315228700637817
1.87000000476837 -0.309087991714478
1.89999997615814 -0.302320241928101
1.92999994754791 -0.29481565952301
1.96000003814697 -0.286388039588928
1.98000001907349 -0.280080795288086
2 -0.272914409637451
2.00999999046326 -0.26868736743927
2.02999997138977 -0.255527257919312
2.05999994277954 -0.234171867370605
2.10999989509583 -0.197231411933899
2.17000007629395 -0.152849197387695
2.21000003814697 -0.124103784561157
2.24000000953674 -0.103318810462952
2.26999998092651 -0.0833760499954224
2.29999995231628 -0.0644075870513916
2.32999992370605 -0.0465373992919922
2.34999990463257 -0.0353037118911743
2.36999988555908 -0.0246796607971191
2.39000010490417 -0.0147539377212524
2.41000008583069 -0.0056988000869751
2.42000007629395 -0.00165081024169922
2.44000005722046 0.00455033779144287
2.46000003814697 0.00995922088623047
2.49000000953674 0.0171771049499512
2.51999998092651 0.023597240447998
2.54999995231628 0.02937912940979
2.58999991416931 0.0362541675567627
2.63000011444092 0.0423020124435425
2.67000007629395 0.0476155281066895
2.71000003814697 0.0522593259811401
2.75 0.0562798976898193
2.78999996185303 0.0597114562988281
2.83999991416931 0.0632103681564331
2.89000010490417 0.0658619403839111
2.94000005722046 0.0676857233047485
2.99000000953674 0.0686874389648438
3.03999996185303 0.0688583850860596
3.07999992370605 0.0683794021606445
3.11999988555908 0.0673255920410156
3.16000008583069 0.0656536817550659
3.20000004768372 0.0632916688919067
3.24000000953674 0.06010901927948
3.26999998092651 0.0570229291915894
3.28999996185303 0.0544791221618652
3.29999995231628 0.0529693365097046
3.30999994277954 0.0510678291320801
3.32999992370605 0.0448048114776611
3.41000008583069 0.0188150405883789
3.4300000667572 0.0131099224090576
3.45000004768372 0.00796699523925781
3.47000002861023 0.0034714937210083
3.48000001907349 0.0015251636505127
3.5 -0.00119590759277344
3.52999997138977 -0.00409698486328125
3.5699999332428 -0.00700569152832031
3.60999989509583 -0.00914573669433594
3.65000009536743 -0.0106602907180786
3.70000004768372 -0.0117700099945068
3.75 -0.0120428800582886
3.78999996185303 -0.0116238594055176
3.82999992370605 -0.0105041265487671
3.84999990463257 -0.00954389572143555
3.85999989509583 -0.00883173942565918
3.91000008583069 -0.0019831657409668
3.9300000667572 -0.000238656997680664
3.96000003814697 0.000821709632873535
4 0.00149071216583252
4.03999996185303 0.00158464908599854
4.07000017166138 0.00114655494689941
4.09000015258789 9.82284545898438e-05
4.13000011444092 -2.71797180175781e-05
4.28999996185303 0
6 0
};
\addlegendentry{$\dot\phi(t)$}
\addplot [semithick, black, mark=*, mark size=1.25, mark options={solid}, forget plot]
table {%
4.23 0
};
\draw (axis cs:4.23,0.03) node[
  anchor=base west,
  text=black,
  rotate=0.0
]{$T_0$};
\end{axis}

\end{tikzpicture}

%% file: pics/simphiarticle_log.tex
\begin{tikzpicture}

\begin{axis}[
width=\plottimewidth,
height=\plottimeheight,
legend cell align={left},
legend style={fill opacity=0.8, draw opacity=1, text opacity=1, draw=white!80!black},
tick align=inside,
tick pos=both,
x grid style={white!69.0196078431373!black},
xlabel={$t$},
xmajorgrids,
xmin=0, xmax=6,
xtick style={color=black},
y grid style={white!69.0196078431373!black},
ymajorgrids,
ymin=-30, ymax=0,
ytick style={color=black},
ytick={-30,-20,-10,0},
yticklabels={-30,-20,-10,0}
]
\addplot [black]
table {%
0 -0.252247022406423
0.01 -0.25235565558621
0.02 -0.25251433968417
0.03 -0.252722457335048
0.04 -0.252979111724526
0.05 -0.253283392779988
0.06 -0.253634432891413
0.07 -0.254031419532024
0.08 -0.254473595940957
0.09 -0.254960258192488
0.1 -0.255490751306081
0.11 -0.256064465314518
0.12 -0.25668083160402
0.13 -0.257339319615158
0.14 -0.258039433907624
0.15 -0.258780711560712
0.16 -0.259562719872344
0.17 -0.260385054319245
0.18 -0.261247336744101
0.19 -0.262149213739666
0.2 -0.263090355203958
0.21 -0.264070453044473
0.22 -0.265089220012633
0.23 -0.266146388652515
0.24 -0.267241710350241
0.25 -0.268374954472444
0.26 -0.269545907583848
0.27 -0.27075437273543
0.28 -0.272000168815792
0.29 -0.273283129959366
0.3 -0.27460310500591
0.31 -0.275959957006481
0.32 -0.277353562771654
0.33 -0.278783812458285
0.34 -0.280250609191579
0.35 -0.281753868719578
0.36 -0.283293519097541
0.37 -0.284869500399966
0.38 -0.286481764458263
0.39 -0.288130274622297
0.4 -0.289815005544232
0.41 -0.291535942983258
0.42 -0.293293083629952
0.43 -0.295086434949146
0.44 -0.296916015040303
0.45 -0.298781852514508
0.46 -0.300683986387272
0.47 -0.302622465986442
0.48 -0.304597350874588
0.49 -0.306608710785295
0.5 -0.308656625572878
0.51 -0.310741185175082
0.52 -0.312862489588379
0.53 -0.315020648855542
0.54 -0.317215783065211
0.55 -0.31944802236321
0.56 -0.321717506975421
0.57 -0.324024387242065
0.58 -0.32636882366325
0.59 -0.328750986955721
0.6 -0.331171058120748
0.61 -0.333629228523139
0.62 -0.336125699981385
0.63 -0.338660684868985
0.64 -0.341234406227019
0.65 -0.343847097888078
0.66 -0.346499004611674
0.67 -0.349190382231302
0.68 -0.351921497813338
0.69 -0.354692629827995
0.7 -0.357504068332585
0.71 -0.360356115167377
0.72 -0.363249084164349
0.73 -0.366183301369192
0.74 -0.369159105276934
0.75 -0.372176847081595
0.76 -0.375236890940327
0.77 -0.378339614252513
0.78 -0.381485407954351
0.79 -0.384674676829481
0.8 -0.387907839836266
0.81 -0.391185330452356
0.82 -0.394507597037257
0.83 -0.397875103213617
0.84 -0.401288328268046
0.85 -0.404747767572306
0.86 -0.408253933025789
0.87 -0.411807353520238
0.88 -0.415408575427765
0.89 -0.419058163113245
0.9 -0.422756699472297
0.91 -0.426504786496078
0.92 -0.430303045864257
0.93 -0.434152119567595
0.94 -0.438052670561665
0.95 -0.44200538345334
0.96 -0.446010965221816
0.97 -0.450070145976015
0.98 -0.454183679750385
0.99 -0.458352345341222
1 -0.462576947185805
1.01 -0.466858316286792
1.02 -0.471197311184492
1.03 -0.475594818979833
1.04 -0.480051756411023
1.05 -0.484569070987151
1.06 -0.489147742182184
1.07 -0.493788782693091
1.08 -0.498493239766101
1.09 -0.503262196595396
1.1 -0.508096773798877
1.11 -0.512998130976001
1.12 -0.517967468353062
1.13 -0.523006028521734
1.14 -0.528115098277132
1.15 -0.533296010562168
1.16 -0.538550146525527
1.17 -0.543878937701187
1.18 -0.549283868318067
1.19 -0.554766477749112
1.2 -0.560328363109914
1.21 -0.565971182017844
1.22 -0.571696655523631
1.23 -0.577506571228354
1.24 -0.583402786600017
1.25 -0.589387232505126
1.26 -0.595461916972112
1.27 -0.601628929205027
1.28 -0.607890443867661
1.29 -0.614248725660153
1.3 -0.620706134212334
1.31 -0.627265129320393
1.32 -0.633928276556119
1.33 -0.640698253280958
1.34 -0.647577855100371
1.35 -0.654570002797755
1.36 -0.661677749791281
1.37 -0.668904290161689
1.38 -0.676252967304319
1.39 -0.683727283264555
1.4 -0.691330908822524
1.41 -0.699067694400461
1.42 -0.706941681874641
1.43 -0.714957117383538
1.44 -0.723118465234861
1.45 -0.731430423026759
1.46 -0.739897938112793
1.47 -0.748526225556843
1.48 -0.757320787742948
1.49 -0.76628743582691
1.5 -0.77543231324157
1.51 -0.784761921496754
1.52 -0.794283148548593
1.53 -0.804003300052108
1.54 -0.813930133856775
1.55 -0.824071898158266
1.56 -0.834437373782481
1.57 -0.84503592115204
1.58 -0.855877532572942
1.59 -0.866972890582955
1.6 -0.878333433226955
1.61 -0.889971427272301
1.62 -0.901900050554717
1.63 -0.914133484859138
1.64 -0.926687020998982
1.65 -0.939577178072511
1.66 -0.952821839260318
1.67 -0.9664404070017
1.68 -0.980453980973114
1.69 -0.99488556301972
1.7 -1.00976029410133
1.71 -1.02510572946003
1.72 -1.0409521596694
1.73 -1.05733298708004
1.74 -1.07428516956273
1.75 -1.09184974654775
1.76 -1.11007246641222
1.77 -1.12900453962452
1.78 -1.1487035492062
1.79 -1.16923455972288
1.8 -1.19067147919662
1.81 -1.21309874656265
1.82 -1.23661344286341
1.83 -1.26132796078071
1.84 -1.28737341981463
1.85 -1.31490409212448
1.86 -1.34410322093991
1.87 -1.37519079326711
1.88 -1.40843411223776
1.89 -1.44416247458075
1.9 -1.48278802944315
1.91 -1.52483623399159
1.92 -1.57099174896464
1.93 -1.62217024286482
1.94 -1.67963592430727
1.95 -1.74520493527025
1.96 -1.82162299602063
1.97 -1.913334308923
1.98 -2.02825884865089
1.99 -2.18275126445739
2 -2.4206437666007
2.01 -2.97222485787334
2.02 -2.79027746275813
2.03 -2.3723464636444
2.04 -2.16774905782958
2.05 -2.03255797582204
2.06 -1.93232534781961
2.07 -1.85325599346789
2.08 -1.78841962911419
2.09 -1.73383904770491
2.1 -1.68702029267492
2.11 -1.64629681419234
2.12 -1.61050007934449
2.13 -1.578779053717
2.14 -1.55049437204271
2.15 -1.52515288009546
2.16 -1.50236536375827
2.17 -1.48181829608712
2.18 -1.463254442282
2.19 -1.44645927718798
2.2 -1.43125133415691
2.21 -1.41747526986075
2.22 -1.40499682641373
2.23 -1.3936991211732
2.24 -1.38347986075342
2.25 -1.37424919327594
2.26 -1.36592799910211
2.27 -1.35844648405554
2.28 -1.35174298536747
2.29 -1.34576293275029
2.3 -1.34045792813375
2.31 -1.33578492043368
2.32 -1.33170545861845
2.33 -1.32818500909349
2.34 -1.32519232314723
2.35 -1.32269883716922
2.36 -1.32067808161031
2.37 -1.31910506093456
2.38 -1.31795553695149
2.39 -1.31720507431197
2.4 -1.31682748905349
2.41 -1.31679147526328
2.42 -1.31704776569343
2.43 -1.31745120286178
2.44 -1.31808010790175
2.45 -1.31893553291872
2.46 -1.32001254661419
2.47 -1.32130568397628
2.48 -1.32280980312738
2.49 -1.32452027106692
2.5 -1.32643298665605
2.51 -1.32854435758829
2.52 -1.33085126546307
2.53 -1.33335103010746
2.54 -1.33604137682176
2.55 -1.3389204075724
2.56 -1.34198657618481
2.57 -1.34523866724988
2.58 -1.34867577836581
2.59 -1.35229730534139
2.6 -1.35610293002529
2.61 -1.36009261047494
2.62 -1.36426657322606
2.63 -1.36862530746782
2.64 -1.37316956096691
2.65 -1.37790033761753
2.66 -1.38281889652344
2.67 -1.38792675254397
2.68 -1.39322567825891
2.69 -1.39871770732786
2.7 -1.40440513923914
2.71 -1.41029054546142
2.72 -1.41637677702917
2.73 -1.42266697361062
2.74 -1.42916457412513
2.75 -1.43587332899533
2.76 -1.44279731413982
2.77 -1.44994094683327
2.78 -1.45730900358481
2.79 -1.46490664021159
2.8 -1.47273941431402
2.81 -1.48081331039236
2.82 -1.48913476788276
2.83 -1.49771071243407
2.84 -1.50654859079786
2.85 -1.51565640976242
2.86 -1.52504277963036
2.87 -1.53471696281971
2.88 -1.54468892826359
2.89 -1.55496941239536
2.9 -1.56556998763997
2.91 -1.57650313949145
2.92 -1.5877823534475
2.93 -1.59942221330259
2.94 -1.61143851257952
2.95 -1.62384838121856
2.96 -1.63667043005785
2.97 -1.64992491614829
2.98 -1.6636339325762
2.99 -1.67782162725033
3 -1.69251445609
3.01 -1.7077414772857
3.02 -1.72353469486956
3.03 -1.73992946183396
3.04 -1.75696495561386
3.05 -1.77468474209609
3.06 -1.7931374487061
3.07 -1.81237757292948
3.08 -1.83246646038928
3.09 -1.85347349709999
3.1 -1.87547757488929
3.11 -1.89856890890692
3.12 -1.92285131417265
3.13 -1.94844508816486
3.14 -1.97549070465312
3.15 -2.00415361015931
3.16 -2.03463054471572
3.17 -2.06715801013606
3.18 -2.10202382895592
3.19 -2.13958326006325
3.2 -2.18028202025513
3.21 -2.22469011057583
3.22 -2.27355318759722
3.23 -2.32787370674244
3.24 -2.38904533884942
3.25 -2.4590891485976
3.26 -2.54110091717685
3.27 -2.64018683583262
3.28 -2.7657121042545
3.29 -2.93797517269934
3.3 -3.21711582555425
3.31 -4.11820092055982
3.32 -3.36224953405297
3.33 -3.03965499747381
3.34 -2.86760219763614
3.35 -2.75380233059358
3.36 -2.6717759021792
3.37 -2.61019322402837
3.38 -2.56309509332884
3.39 -2.52685720690083
3.4 -2.49901979780352
3.41 -2.47779368446806
3.42 -2.46183835912851
3.43 -2.45014473304526
3.44 -2.44195461092236
3.45 -2.43669441012098
3.46 -2.43391677764601
3.47 -2.43324400311103
3.48 -2.4342888201573
3.49 -2.4360763101016
3.5 -2.43860310076825
3.51 -2.44219343368679
3.52 -2.44683895749329
3.53 -2.4525184878104
3.54 -2.45921401237673
3.55 -2.46691380308362
3.56 -2.47561305202884
3.57 -2.48531393490815
3.58 -2.49602559622438
3.59 -2.50776421119168
3.6 -2.52055318188236
3.61 -2.53442349695308
3.62 -2.54941427946419
3.63 -2.5655735521262
3.64 -2.58295925965238
3.65 -2.60164060311491
3.66 -2.62169976231315
3.67 -2.64323411146168
3.68 -2.66635907481084
3.69 -2.69121182818661
3.7 -2.71795613957084
3.71 -2.74678877253726
3.72 -2.77794807702337
3.73 -2.81172570800646
3.74 -2.84848292497856
3.75 -2.88867378281597
3.76 -2.93287901488837
3.77 -2.98185711298798
3.78 -3.03662427100664
3.79 -3.09858533866736
3.8 -3.16976085776574
3.81 -3.25321032898337
3.82 -3.35390112454696
3.83 -3.48074852291
3.84 -3.6524767030044
3.85 -3.92154009434521
3.86 -4.62896332637046
3.87 -4.18979172767506
3.88 -3.87006888263828
3.89 -3.74066296341946
3.9 -3.68792197847542
3.91 -3.67368837631474
3.92 -3.67912179634764
3.93 -3.69356705245019
3.94 -3.70528444814546
3.95 -3.72239369909398
3.96 -3.74563497639777
3.97 -3.7753905136284
3.98 -3.81219569278339
3.99 -3.85688596630447
4 -3.9107449891905
4.01 -3.97574605023438
4.02 -4.05500615035114
4.03 -4.1537286034888
4.04 -4.28141873066069
4.05 -4.45814062251992
4.06 -4.73880738786288
4.07 -5.44334499732269
4.08 -5.22894165510594
4.09 -5.37072101187896
4.1 -5.48405584885042
4.11 -5.64953228335847
4.12 -5.98368600619972
4.13 -7.71764683168058
4.14 -8.06170962152654
4.15 -9.57875077444486
4.16 -12.6180725048019
4.17 -14.4706333755157
4.18 -20.3524682437675
4.19 -38.1916315344905
4.2 -80.1498233021185
4.21 -42.4351645842267
4.22 -111.556256935223
4.23 -inf
4.24 -inf
4.25 -inf
4.26 -inf
4.27 -inf
4.28 -inf
4.29 -inf
};
\addlegendentry{$\log\lvert\phi(t)\rvert$}
\addplot [black, dashed]
table {%
0 -inf
0.01 -2.13193828133015
0.02 -1.84111924607472
0.03 -1.67371532498487
0.04 -1.55649271619389
0.05 -1.46660970695381
0.06 -1.39393381417112
0.07 -1.33307993710654
0.08 -1.28084442759787
0.09 -1.2351672695821
0.1 -1.19464596728656
0.11 -1.15828290072075
0.12 -1.12534336349382
0.13 -1.09527078148364
0.14 -1.06763353707586
0.15 -1.04209024403941
0.16 -1.01836629146314
0.17 -0.996237537372553
0.18 -0.975518688287024
0.19 -0.956054837461562
0.2 -0.937715185214808
0.21 -0.920388299641818
0.22 -0.903978485808105
0.23 -0.888402966442502
0.24 -0.87358966596725
0.25 -0.859475449417002
0.26 -0.846004708711389
0.27 -0.833128217264164
0.28 -0.820802194105163
0.29 -0.808987533196616
0.3 -0.797649164183413
0.31 -0.786755518596143
0.32 -0.776278081322463
0.33 -0.766191011527303
0.34 -0.756470820521376
0.35 -0.747096096624152
0.36 -0.738047269038266
0.37 -0.729306404289647
0.38 -0.720857029995889
0.39 -0.712683981681681
0.4 -0.704773269121975
0.41 -0.697111959304517
0.42 -0.689688073596005
0.43 -0.682490497095783
0.44 -0.675508898486841
0.45 -0.668733658960921
0.46 -0.662155809014501
0.47 -0.655766972094358
0.48 -0.649559314222607
0.49 -0.643525498857276
0.5 -0.637658646350185
0.51 -0.63195229745275
0.52 -0.626400380395392
0.53 -0.620997181129757
0.54 -0.615737316376977
0.55 -0.610615709171231
0.56 -0.605627566627276
0.57 -0.600768359694381
0.58 -0.596033804688182
0.59 -0.591419846417035
0.6 -0.586922642741152
0.61 -0.582538550421602
0.62 -0.578264112132616
0.63 -0.574096044524889
0.64 -0.57003122724003
0.65 -0.566066692787181
0.66 -0.562199617202422
0.67 -0.558427311419967
0.68 -0.554747213291559
0.69 -0.551156880197023
0.7 -0.547653982194711
0.71 -0.544236295665693
0.72 -0.540901697410105
0.73 -0.537648159158101
0.74 -0.534473742461453
0.75 -0.531376593935088
0.76 -0.528354940820668
0.77 -0.525407086846949
0.78 -0.522531408363918
0.79 -0.519726350729791
0.8 -0.516990424931821
0.81 -0.514322204423521
0.82 -0.511720322162441
0.83 -0.509183467833983
0.84 -0.50671038524798
0.85 -0.504299869895866
0.86 -0.50195076665731
0.87 -0.499661967646063
0.88 -0.497432410185618
0.89 -0.495261074906072
0.9 -0.493146983954204
0.91 -0.491089199309484
0.92 -0.489086821199263
0.93 -0.487138986606925
0.94 -0.485244867867275
0.95 -0.483403671343883
0.96 -0.481614636183483
0.97 -0.479877033142939
0.98 -0.47819016348459
0.99 -0.476553357936144
1 -0.474965975711556
1.01 -0.473427403589591
1.02 -0.471937055047057
1.03 -0.470494369443885
1.04 -0.469098811257474
1.05 -0.467749869363897
1.06 -0.466447056363794
1.07 -0.465189907950891
1.08 -0.463977982321297
1.09 -0.462810859621865
1.1 -0.461688141436057
1.11 -0.460609450305865
1.12 -0.459574429288504
1.13 -0.458582741546695
1.14 -0.457634069971473
1.15 -0.456728116836581
1.16 -0.455864603483596
1.17 -0.455043270037076
1.18 -0.454263875149065
1.19 -0.453526195772449
1.2 -0.452830026962712
1.21 -0.452175181707752
1.22 -0.451561490785522
1.23 -0.450988802649337
1.24 -0.450456983340802
1.25 -0.449965916430402
1.26 -0.4495155029859
1.27 -0.449105661568786
1.28 -0.448736328259136
1.29 -0.448407456709335
1.3 -0.448119018227229
1.31 -0.447871001889428
1.32 -0.447663414685535
1.33 -0.447496281694302
1.34 -0.447369646292769
1.35 -0.447283570399653
1.36 -0.447238134754398
1.37 -0.447233439233467
1.38 -0.44726960320566
1.39 -0.447346765928453
1.4 -0.447465086987566
1.41 -0.44762474678223
1.42 -0.447825947058896
1.43 -0.448068911496432
1.44 -0.448353886346171
1.45 -0.448681141130576
1.46 -0.449050969404664
1.47 -0.449463689584812
1.48 -0.449919645850068
1.49 -0.450419209121676
1.5 -0.45096277812714
1.51 -0.451550780555916
1.52 -0.452183674314594
1.53 -0.452861948890399
1.54 -0.453586126832859
1.55 -0.454356765364699
1.56 -0.455174458134381
1.57 -0.456039837124272
1.58 -0.456953574730224
1.59 -0.457916386030401
1.6 -0.458929031263601
1.61 -0.459992318540076
1.62 -0.461107106811097
1.63 -0.462274309127279
1.64 -0.463494896220097
1.65 -0.464769900446265
1.66 -0.466100420140797
1.67 -0.467487624431913
1.68 -0.468932758579709
1.69 -0.470437149911003
1.7 -0.472002214435458
1.71 -0.47362946424344
1.72 -0.475320515804826
1.73 -0.477077099310993
1.74 -0.478901069230672
1.75 -0.480794416285736
1.76 -0.482759281097387
1.77 -0.484797969809298
1.78 -0.486912972065763
1.79 -0.489106981814933
1.8 -0.491382921526782
1.81 -0.493743970572575
1.82 -0.496193598721454
1.83 -0.498735605991073
1.84 -0.501374170473489
1.85 -0.504113906290608
1.86 -0.506959934586245
1.87 -0.509917971545119
1.88 -0.512994439022416
1.89 -0.516196605768909
1.9 -0.519532770957882
1.91 -0.523012507674251
1.92 -0.526646993913237
1.93 -0.530449475786747
1.94 -0.534435938980857
1.95 -0.538626125556636
1.96 -0.543045161910974
1.97 -0.547726364398014
1.98 -0.552716595555137
1.99 -0.558088185520488
2 -0.56397344078505
2.01 -0.570752808284053
2.02 -0.581199819473072
2.03 -0.592562867777957
2.04 -0.604605908637172
2.05 -0.617250903228057
2.06 -0.630465356524797
2.07 -0.644237659406337
2.08 -0.658568241092106
2.09 -0.673465601566104
2.1 -0.688944293929517
2.11 -0.70502382489188
2.12 -0.721728046990065
2.13 -0.739084851015604
2.14 -0.757126073351924
2.15 -0.775887587407442
2.16 -0.79540957875221
2.17 -0.815737020794153
2.18 -0.836920376889607
2.19 -0.859016559164047
2.2 -0.882090177662248
2.21 -0.906215119885348
2.22 -0.931476514588067
2.23 -0.957973158944975
2.24 -0.985820528762911
2.25 -1.01515455204358
2.26 -1.04613641431164
2.27 -1.07895879280945
2.28 -1.11385410930929
2.29 -1.1511056895363
2.3 -1.19106319579235
2.31 -1.23416449611275
2.32 -1.28096751195758
2.33 -1.33219807643353
2.34 -1.38882455959882
2.35 -1.45217952131923
2.36 -1.52416917315356
2.37 -1.60765981358644
2.38 -1.70725808548365
2.39 -1.83109378776672
2.4 -1.99570728227448
2.41 -2.24421769655793
2.42 -2.78231507001993
2.43 -2.79886706595837
2.44 -2.34195260952364
2.45 -2.13504753969171
2.46 -2.00177299775842
2.47 -1.90418190535202
2.48 -1.8276649093359
2.49 -1.7650511667639
2.5 -1.71229073785575
2.51 -1.66687409946779
2.52 -1.62713979197154
2.53 -1.5919316590666
2.54 -1.56041327322143
2.55 -1.53196032508109
2.56 -1.50609472969212
2.57 -1.48244244784445
2.58 -1.46070547401358
2.59 -1.44064264447744
2.6 -1.42205613308949
2.61 -1.40478172706886
2.62 -1.38868168175535
2.63 -1.37363937587209
2.64 -1.35955524968805
2.65 -1.34634367403002
2.66 -1.33393050580199
2.67 -1.32225115731262
2.68 -1.31124905532455
2.69 -1.30087439932733
2.7 -1.29108315212617
2.71 -1.28183621265919
2.72 -1.2730987331145
2.73 -1.26483955132063
2.74 -1.25703071597814
2.75 -1.24964708723853
2.76 -1.24266599887242
2.77 -1.23606697112239
2.78 -1.22983146553389
2.79 -1.22394267476587
2.8 -1.21838534172023
2.81 -1.21314560338443
2.82 -1.20821085562017
2.83 -1.20356963580188
2.84 -1.19921152074903
2.85 -1.19512703783455
2.86 -1.19130758750865
2.87 -1.18774537577161
2.88 -1.18443335537131
2.89 -1.18136517470447
2.9 -1.17853513357016
2.91 -1.17593814506866
2.92 -1.17356970306289
2.93 -1.17142585472665
2.94 -1.16950317779939
2.95 -1.16779876225242
2.96 -1.16631019614999
2.97 -1.16503555556296
2.98 -1.16397339846444
2.99 -1.16312276260917
3 -1.16248316747286
3.01 -1.16205462040744
3.02 -1.16183762725595
3.03 -1.16183320777041
3.04 -1.16204291629146
3.05 -1.16246886828548
3.06 -1.16311377350123
3.07 -1.16398097671155
3.08 -1.16507450726146
3.09 -1.1663991389669
3.1 -1.16796046232472
3.11 -1.16976497153661
3.12 -1.17182016956727
3.13 -1.17413469542042
3.14 -1.1767184791313
3.15 -1.17958293180001
3.16 -1.18274118057524
3.17 -1.18620836223435
3.18 -1.19000199453484
3.19 -1.19414245290804
3.2 -1.19865359319374
3.21 -1.20356358233356
3.22 -1.20890603457924
3.23 -1.21472161337542
3.24 -1.22106037512394
3.25 -1.22798536080589
3.26 -1.23557843581126
3.27 -1.24395056421159
3.28 -1.2532620068507
3.29 -1.26376941924911
3.3 -1.27597537165358
3.31 -1.29185303079558
3.32 -1.31861731416128
3.33 -1.34867556815511
3.34 -1.38179567220357
3.35 -1.4181012442656
3.36 -1.45786720519502
3.37 -1.50146681816653
3.38 -1.54935970491338
3.39 -1.60211503443964
3.4 -1.6604790641879
3.41 -1.7254953784429
3.42 -1.79869401958241
3.43 -1.88240153133082
3.44 -1.98031343544196
3.45 -2.09870691557631
3.46 -2.24945541360631
3.47 -2.4594904765764
3.48 -2.81669290664793
3.49 -5.54990347785264
3.5 -2.92230102152867
3.51 -2.64687287285574
3.52 -2.49270428984944
3.53 -2.38753092038474
3.54 -2.30893795446003
3.55 -2.24703780384868
3.56 -2.19660466163924
3.57 -2.1545472252671
3.58 -2.1188896079005
3.59 -2.0882954391264
3.6 -2.06182127659597
3.61 -2.03877847992999
3.62 -2.01865099937244
3.63 -2.00104401313281
3.64 -1.98565054893038
3.65 -1.97222908406137
3.66 -1.96058811781652
3.67 -1.95057533065751
3.68 -1.94206986178478
3.69 -1.93497677801862
3.7 -1.92922313950192
3.71 -1.92475528126406
3.72 -1.92153707470904
3.73 -1.91954904003951
3.74 -1.91878826990858
3.75 -1.91926921270484
3.76 -1.92102546832391
3.77 -1.92411289450042
3.78 -1.92861454960093
3.79 -1.93464839041764
3.8 -1.94237938179717
3.81 -1.9520392054833
3.82 -1.96396033857281
3.83 -1.97864109913407
3.84 -1.99689177683916
3.85 -2.02027421949025
3.86 -2.05395216597995
3.87 -2.12831379142412
3.88 -2.22381619930321
3.89 -2.34553523533742
3.9 -2.50026636397068
3.91 -2.70263308353134
3.92 -2.99847095803901
3.93 -3.62224867717609
3.94 -3.69217620950154
3.95 -3.26408055984734
3.96 -3.08525090652447
3.97 -2.9791602149358
3.98 -2.90873420562318
3.99 -2.86006714592658
4 -2.82662243608155
4.01 -2.80501898741964
4.02 -2.793560888221
4.03 -2.79170938047879
4.04 -2.80005942842988
4.05 -2.82091224374295
4.06 -2.86041486361179
4.07 -2.94059997199689
4.08 -3.29314349712987
4.09 -4.00752108196165
4.1 -4.74372570009511
4.11 -4.16545936698972
4.12 -4.09623950694844
4.13 -4.56510025517454
4.14 -5.3689575515775
4.15 -6.47909771874189
4.16 -8.66453916430823
4.17 -9.88544820621622
4.18 -13.8073894375966
4.19 -25.699037347556
4.2 -47.4592283903074
4.21 -57.5991992071007
4.22 -92.1597453825988
4.23 -179.069745159379
4.24 -inf
4.25 -inf
4.26 -inf
4.27 -inf
4.28 -inf
4.29 -inf
};
\addlegendentry{$\log\lvert\dot\phi(t)\rvert$}
\end{axis}

\end{tikzpicture}

%% file: pics/simalphabetaarticle.tex
\begin{tikzpicture}

\begin{axis}[
width=\plottimewidth,
height=\plottimeheight,
legend cell align={left},
legend style={fill opacity=0.8, draw opacity=1, text opacity=1, at={(0.97,0.03)},
legend image post style={dash phase=0pt}, 
anchor=south east, draw=white!80!black},
tick align= inside,
tick pos=both,
x grid style={white!69.0196078431373!black},
xlabel={$t$},
xmajorgrids,
xmin=0, xmax=6,
xtick style={color=black},
y grid style={white!69.0196078431373!black},
ymajorgrids,
ymin=-0.2, ymax=0.05,
ytick style={color=black},
ytick={-0.2,-0.15,-0.1,-0.05,0,0.05},
yticklabels={-0.20,-0.15,-0.10,-0.05,0,0.05}
]
\addplot [black]
table {%
0 0
0.01 -0.00310234245109961
0.02 -0.00606050439958249
0.03 -0.00891071784167689
0.04 -0.0116717318321778
0.05 -0.0143554995521099
0.06 -0.0169704829888754
0.07 -0.0195230685777122
0.08 -0.0220182900940923
0.09 -0.0244602418801466
0.1 -0.0268523345659767
0.11 -0.0291974629620598
0.12 -0.031498121465185
0.13 -0.0337564863523963
0.14 -0.0359744762539783
0.15 -0.0381537977200498
0.16 -0.0402959802926041
0.17 -0.0424024039969934
0.18 -0.0444743212352622
0.19 -0.0465128744648771
0.2 -0.0485191106502555
0.21 -0.050493993205793
0.22 -0.0524384119627334
0.23 -0.0543531915604155
0.24 -0.056239098567524
0.25 -0.0580968475695485
0.26 -0.0599271064071227
0.27 -0.0617305007111661
0.28 -0.0635076178512514
0.29 -0.0652590103909203
0.3 -0.0669851991260185
0.31 -0.0686866757682676
0.32 -0.0703639053253222
0.33 -0.0720173282198064
0.34 -0.0736473621827785
0.35 -0.0752544039513743
0.36 -0.0768388307957217
0.37 -0.0784010018964095
0.38 -0.0799412595906408
0.39 -0.0814599305025876
0.4 -0.082957326571285
0.41 -0.0844337459875755
0.42 -0.0858894740500726
0.43 -0.0873247839488151
0.44 -0.0887399374841723
0.45 -0.0901351857276215
0.46 -0.091510769630209
0.47 -0.0928669205838122
0.48 -0.0942038609397241
0.49 -0.0955218044885604
0.5 -0.0968209569050433
0.51 -0.0981015161608241
0.52 -0.0993636729081653
0.53 -0.100607610837005
0.54 -0.10183350700766
0.55 -0.103041532161205
0.56 -0.104231851009345
0.57 -0.105404622505426
0.58 -0.106560000098085
0.59 -0.107698131968859
0.6 -0.108819161255002
0.61 -0.109923226258581
0.62 -0.111010460642892
0.63 -0.112080993617081
0.64 -0.113134950109828
0.65 -0.114172450932836
0.66 -0.115193612934842
0.67 -0.116198549146769
0.68 -0.117187368918614
0.69 -0.118160178048606
0.7 -0.119117078905128
0.71 -0.12005817054185
0.72 -0.120983548806493
0.73 -0.121893306443614
0.74 -0.12278753319175
0.75 -0.12366631587527
0.76 -0.124529738491209
0.77 -0.125377882291388
0.78 -0.126210825860057
0.79 -0.127028645187318
0.8 -0.127831413738522
0.81 -0.128619202519876
0.82 -0.129392080140417
0.83 -0.130150112870548
0.84 -0.13089336469729
0.85 -0.131621897376395
0.86 -0.132335770481467
0.87 -0.133035041450207
0.88 -0.133719765627914
0.89 -0.13438999630833
0.9 -0.135045784771949
0.91 -0.13568718032187
0.92 -0.136314230317277
0.93 -0.13692698020463
0.94 -0.13752547354663
0.95 -0.138109752049022
0.96 -0.138679855585298
0.97 -0.139235822219345
0.98 -0.139777688226085
0.99 -0.140305488110145
1 -0.140819254622601
1.01 -0.141319018775807
1.02 -0.141804809856348
1.03 -0.142276655436132
1.04 -0.142734581381631
1.05 -0.143178611861282
1.06 -0.143608769351056
1.07 -0.144025074638186
1.08 -0.144427546823054
1.09 -0.144816203319228
1.1 -0.145191059851621
1.11 -0.145552130452771
1.12 -0.145899427457189
1.13 -0.146232961493771
1.14 -0.146552741476211
1.15 -0.146858774591388
1.16 -0.147151066285669
1.17 -0.147429620249073
1.18 -0.147694438397228
1.19 -0.147945520851063
1.2 -0.148182865914137
1.21 -0.148406470047532
1.22 -0.148616327842214
1.23 -0.148812431988744
1.24 -0.148994773244242
1.25 -0.149163340396459
1.26 -0.149318120224836
1.27 -0.149459097458386
1.28 -0.149586254730243
1.29 -0.149699572528698
1.3 -0.149799029144524
1.31 -0.149884600614378
1.32 -0.149956260660053
1.33 -0.150013980623324
1.34 -0.150057729396105
1.35 -0.150087473345639
1.36 -0.150103176234369
1.37 -0.150104799134142
1.38 -0.15009230033436
1.39 -0.150065635243643
1.4 -0.150024756284532
1.41 -0.149969612780733
1.42 -0.149900150836316
1.43 -0.149816313206273
1.44 -0.149718039157729
1.45 -0.149605264321076
1.46 -0.149477920530184
1.47 -0.149335935650786
1.48 -0.149179233396014
1.49 -0.149007733127959
1.5 -0.148821349644018
1.51 -0.148619992946625
1.52 -0.148403567994818
1.53 -0.14817197443592
1.54 -0.147925106315373
1.55 -0.147662851762571
1.56 -0.14738509265022
1.57 -0.147091704224478
1.58 -0.146782554702739
1.59 -0.146457504835523
1.6 -0.146116407428434
1.61 -0.145759106819605
1.62 -0.145385438307349
1.63 -0.144995227521993
1.64 -0.144588289734937
1.65 -0.144164429096892
1.66 -0.143723437795964
1.67 -0.143265095124717
1.68 -0.142789166443484
1.69 -0.142295402024972
1.7 -0.141783535762509
1.71 -0.141253283720964
1.72 -0.140704342505364
1.73 -0.140136387417213
1.74 -0.139549070362343
1.75 -0.138942017466384
1.76 -0.138314826344173
1.77 -0.137667062957048
1.78 -0.136998257976053
1.79 -0.136307902548595
1.8 -0.135595443339224
1.81 -0.134860276679806
1.82 -0.134101741617003
1.83 -0.133319111580896
1.84 -0.132511584310562
1.85 -0.131678269549654
1.86 -0.130818173850773
1.87 -0.129930181575207
1.88 -0.129013030801628
1.89 -0.128065282291939
1.9 -0.127085278781199
1.91 -0.126071090440213
1.92 -0.125020439989869
1.93 -0.123930596810918
1.94 -0.122798221787512
1.95 -0.121619129716209
1.96 -0.120387904476476
1.97 -0.11909722777064
1.98 -0.117736581366566
1.99 -0.116289320106403
2 -0.114724078008322
2.01 -0.112947132264675
2.02 -0.110262591997959
2.03 -0.107415053511345
2.04 -0.104477340329944
2.05 -0.101479218420993
2.06 -0.0984379725705276
2.07 -0.0953652952109468
2.08 -0.0922698417625865
2.09 -0.0891584306177113
2.1 -0.0860366957663499
2.11 -0.0829094801375453
2.12 -0.0797810908103399
2.13 -0.0766554726668731
2.14 -0.0735363286378634
2.15 -0.070427201163699
2.16 -0.0673315230145112
2.17 -0.0642526428912879
2.18 -0.061193830612958
2.19 -0.0581582671088721
2.2 -0.0551490250584798
2.21 -0.0521690462461666
2.22 -0.0492211212138849
2.23 -0.0463078756024994
2.24 -0.0434317659372411
2.25 -0.0405950859142556
2.26 -0.0377999828331451
2.27 -0.0350484829128461
2.28 -0.0323425238892282
2.29 -0.0296839934699149
2.3 -0.0270747728139707
2.31 -0.0245167851284708
2.32 -0.0220120507259943
2.33 -0.0195627516007312
2.34 -0.0171713111286347
2.35 -0.0148404987129306
2.36 -0.0125735769428509
2.37 -0.0103745245642928
2.38 -0.00824840413769462
2.39 -0.00620203523393599
2.4 -0.00424542071477934
2.41 -0.00239557802375798
2.42 -0.00069392565960863
2.43 0.000667976141741175
2.44 0.00191282574064335
2.45 0.00308021106914395
2.46 0.00418653782289236
2.47 0.00524138560438402
2.48 0.00625119587898679
2.49 0.00722065924093927
2.5 0.00815337240304848
2.51 0.00905219534611781
2.52 0.00991946444159515
2.53 0.0107571285321397
2.54 0.0115668404131115
2.55 0.0123500208794718
2.56 0.0131079050623181
2.57 0.0138415768739705
2.58 0.0145519952032802
2.59 0.0152400142272219
2.6 0.015906399425679
2.61 0.0165518403932375
2.62 0.0171769612200361
2.63 0.017782328998075
2.64 0.0183684608614255
2.65 0.0189358298651224
2.66 0.0194848699335332
2.67 0.020015980055294
2.68 0.0205295278623276
2.69 0.0210258527008913
2.7 0.0215052682802236
2.71 0.0219680649672368
2.72 0.022414511782445
2.73 0.0228448581419614
2.74 0.0232593353822249
2.75 0.0236581580976131
2.76 0.0240415253158758
2.77 0.0244096215320963
2.78 0.0247626176184414
2.79 0.0251006716241227
2.8 0.0254239294776402
2.81 0.0257325256014118
2.82 0.0260265834472294
2.83 0.0263062159595677
2.84 0.0265715259725464
2.85 0.0268226065452798
2.86 0.0270595412394
2.87 0.0272824043416838
2.88 0.0274912610339293
2.89 0.0276861675114853
2.9 0.0278671710511301
2.91 0.0280343100282916
2.92 0.0281876138828914
2.93 0.0283271030323601
2.94 0.0284527887295855
2.95 0.0285646728627057
2.96 0.0286627476927107
2.97 0.0287469955237503
2.98 0.0288173882998175
2.99 0.0288738871200473
3 0.0289164416631909
3.01 0.0289449895098149
3.02 0.0289594553483617
3.03 0.0289597500482729
3.04 0.0289457695797758
3.05 0.0289173937554737
3.06 0.0288744847633163
3.07 0.0288168854534919
3.08 0.028744417332832
3.09 0.0286568782087701
3.1 0.0285540394098498
3.11 0.0284356424898988
3.12 0.0283013952963778
3.13 0.0281509672472452
3.14 0.0279839836107127
3.15 0.0278000185119862
3.16 0.0275985862903013
3.17 0.0273791306818821
3.18 0.0271410110829553
3.19 0.0268834848059542
3.2 0.0266056837016294
3.21 0.0263065826346265
3.22 0.0259849557942246
3.23 0.0256393141432676
3.24 0.0252678122813269
3.25 0.0248681029212661
3.26 0.0244370952215454
3.27 0.0239705198498981
3.28 0.0234620529217521
3.29 0.0229012188850092
3.3 0.0222665352856274
3.31 0.0214671794216576
3.32 0.0201841610491358
3.33 0.0188344285857483
3.34 0.0174514832481205
3.35 0.0160519120914807
3.36 0.0146474091070417
3.37 0.0132483318735745
3.38 0.0118650155574117
3.39 0.0105078286571338
3.4 0.00918647616843695
3.41 0.00790920147119671
3.42 0.00668242352543404
3.43 0.00551095018814492
3.44 0.00439860113012422
3.45 0.00334904591217928
3.46 0.00236686018082786
3.47 0.00145927633748309
3.48 0.000641113777968794
3.49 -1.18501688642225e-06
3.5 -0.00050272096311681
3.51 -0.000947883894124571
3.52 -0.00135183512640394
3.53 -0.00172225358138238
3.54 -0.00206390670745277
3.55 -0.00238007264605964
3.56 -0.00267315016015174
3.57 -0.00294496943445364
3.58 -0.00319696888176408
3.59 -0.00343030374759857
3.6 -0.00364591668673004
3.61 -0.00384458561802264
3.62 -0.00402695723480085
3.63 -0.00419357103346428
3.64 -0.00434487681547014
3.65 -0.00448124752192886
3.66 -0.00460298859786485
3.67 -0.00471034466283694
3.68 -0.00480350398278419
3.69 -0.00488260103551649
3.7 -0.0049477173026193
3.71 -0.00499888027826532
3.72 -0.00503606053945725
3.73 -0.00505916655074966
3.74 -0.00506803665163413
3.75 -0.00506242735495657
3.76 -0.00504199660127866
3.77 -0.00500627984414371
3.78 -0.00495465554772427
3.79 -0.00488629436355221
3.8 -0.00480008176355597
3.81 -0.00469449424973779
3.82 -0.00456738560021365
3.83 -0.00441557146160376
3.84 -0.00423385753220562
3.85 -0.00401193416221583
3.86 -0.00371258025795609
3.87 -0.00312834200067722
3.88 -0.00251079825191281
3.89 -0.00189711117196278
3.9 -0.00132850028814774
3.91 -0.000833671458687236
3.92 -0.000421849033931863
3.93 -0.000100318109610056
3.94 8.53988077210035e-05
3.95 0.000228848176931957
3.96 0.000345443567727616
3.97 0.000441029664550114
3.98 0.000518673203396514
3.99 0.000580177484865147
4 0.000626621916926517
4.01 0.000658580745682513
4.02 0.000676187485295311
4.03 0.000679076398091601
4.04 0.000666144717773806
4.05 0.000634915271450438
4.06 0.00057971315192681
4.07 0.000481978920798681
4.08 0.000214034863372524
4.09 4.13149217888288e-05
4.1 -7.58406972037936e-06
4.11 -2.8718956110364e-05
4.12 -3.36812739661638e-05
4.13 -1.14426802340258e-05
4.14 -1.79750648269001e-06
4.15 -1.39485896535364e-07
4.16 9.1009945145738e-10
4.17 5.47242490494045e-11
4.18 6.54995896710325e-15
4.19 8.40602484651393e-27
4.2 1.46015724141229e-48
4.21 -1.0578614955798e-58
4.22 -2.90993061912274e-93
4.23 -3.5799986059973e-180
4.24 0
4.25 0
4.26 0
4.27 0
4.28 0
4.29 0
4.3 0
4.31 0
4.32 0
4.33 0
4.34 0
4.35 0
4.36 0
4.37 0
4.38 0
4.39 0
4.4 0
4.41 0
4.42 0
4.43 0
4.44 0
4.45 0
4.46 0
4.47 0
4.48 0
4.49 0
4.5 0
4.51 0
4.52 0
4.53 0
4.54 0
4.55 0
4.56 0
4.57 0
4.58 0
4.59 0
4.6 0
4.61 0
4.62 0
4.63 0
4.64 0
4.65 0
4.66 0
4.67 0
4.68 0
4.69 0
4.7 0
4.71 0
4.72 0
4.73 0
4.74 0
4.75 0
4.76 0
4.77 0
4.78 0
4.79 0
4.8 0
4.81 0
4.82 0
4.83 0
4.84 0
4.85 0
4.86 0
4.87 0
4.88 0
4.89 0
4.9 0
4.91 0
4.92 0
4.93 0
4.94 0
4.95 0
4.96 0
4.97 0
4.98 0
4.99 0
5 0
5.01 0
5.02 0
5.03 0
5.04 0
5.05 0
5.06 0
5.07 0
5.08 0
5.09 0
5.1 0
5.11 0
5.12 0
5.13 0
5.14 0
5.15 0
5.16 0
5.17 0
5.18 0
5.19 0
5.2 0
5.21 0
5.22 0
5.23 0
5.24 0
5.25 0
5.26 0
5.27 0
5.28 0
5.29 0
5.3 0
5.31 0
5.32 0
5.33 0
5.34 0
5.35 0
5.36 0
5.37 0
5.38 0
5.39 0
5.4 0
5.41 0
5.42 0
5.43 0
5.44 0
5.45 0
5.46 0
5.47 0
5.48 0
5.49 0
5.5 0
5.51 0
5.52 0
5.53 0
5.54 0
5.55 0
5.56 0
5.57 0
5.58 0
5.59 0
5.6 0
5.61 0
5.62 0
5.63 0
5.64 0
5.65 0
5.66 0
5.67 0
5.68 0
5.69 0
5.7 0
5.71 0
5.72 0
5.73 0
5.74 0
5.75 0
5.76 0
5.77 0
5.78 0
5.79 0
5.8 0
5.81 0
5.82 0
5.83 0
5.84 0
5.85 0
5.86 0
5.87 0
5.88 0
5.89 0
5.9 0
5.91 0
5.92 0
5.93 0
5.94 0
5.95 0
5.96 0
5.97 0
5.98 0
5.99 0
6 0
};
\addlegendentry{$\beta(1, t)$}
\addplot [black, dashed, dash phase = 8pt] 
table {%
0 0
0.01 0
0.02 0
0.03 0
0.04 0
0.05 0
0.06 0
0.07 0
0.08 0
0.09 0
0.1 0
0.11 0
0.12 0
0.13 0
0.14 0
0.15 0
0.16 0
0.17 0
0.18 0
0.19 0
0.2 0
0.21 -1.52257264322914e-05
0.22 -9.96635876289693e-05
0.23 -0.000349792830926801
0.24 -0.000877558560946674
0.25 -0.00177039123473539
0.26 -0.00306343392399599
0.27 -0.00473556575874875
0.28 -0.00672450626744906
0.29 -0.00894908489317061
0.3 -0.0113286128122274
0.31 -0.0137947778504276
0.32 -0.0162961346042341
0.33 -0.0187975117287615
0.34 -0.021276928147185
0.35 -0.0237219133863527
0.36 -0.0261262585733715
0.37 -0.0284875641932649
0.38 -0.0308055751785292
0.39 -0.0330811396688176
0.4 -0.0353156023599415
0.41 -0.0375104744529257
0.42 -0.0396672673913748
0.43 -0.0417874176184877
0.44 -0.0438722588517549
0.45 -0.0459230174567734
0.46 -0.0479408179781279
0.47 -0.0499266923589168
0.48 -0.0518815898349367
0.49 -0.0538063862376888
0.5 -0.0557018922757562
0.51 -0.0575688607350503
0.52 -0.0594079926847674
0.53 -0.0612199428198306
0.54 -0.0630053240715552
0.55 -0.0647647116037198
0.56 -0.0664986462930418
0.57 -0.0682076377758849
0.58 -0.0698921671283264
0.59 -0.0715526892346364
0.6 -0.0731896348894724
0.61 -0.0748034126712652
0.62 -0.0763944106179726
0.63 -0.077962997731299
0.64 -0.0795095253313473
0.65 -0.0810343282803021
0.66 -0.0825377260909672
0.67 -0.0840200239336917
0.68 -0.0854815135533084
0.69 -0.086922474106116
0.7 -0.0883431729255939
0.71 -0.0897438662244058
0.72 -0.0911247997392871
0.73 -0.0924862093245918
0.74 -0.0938283214995723
0.75 -0.0951513539538645
0.76 -0.0964555160151303
0.77 -0.0977410090823599
0.78 -0.099008027027947
0.79 -0.100256756571309
0.8 -0.101487377626532
0.81 -0.102700063626247
0.82 -0.103894981823742
0.83 -0.105072293575078
0.84 -0.106232154602831
0.85 -0.107374715242914
0.86 -0.108500120675777
0.87 -0.109608511143188
0.88 -0.110700022151675
0.89 -0.111774784663593
0.9 -0.112832925276727
0.91 -0.113874566393232
0.92 -0.114899826378651
0.93 -0.115908819711704
0.94 -0.116901657125442
0.95 -0.11787844574036
0.96 -0.11883928918998
0.97 -0.119784287739371
0.98 -0.120713538397068
0.99 -0.12162713502078
1 -0.122525168417259
1.01 -0.123407726436678
1.02 -0.124274894061841
1.03 -0.125126753492499
1.04 -0.125963384225058
1.05 -0.126784863127924
1.06 -0.12759126451271
1.07 -0.128382660201526
1.08 -0.129159119590537
1.09 -0.129920709709992
1.1 -0.130667495280865
1.11 -0.131399538768284
1.12 -0.13211690043189
1.13 -0.13281963837324
1.14 -0.133507808580406
1.15 -0.134181464969857
1.16 -0.13484065942574
1.17 -0.135485441836655
1.18 -0.13611586013001
1.19 -0.136731960304031
1.2 -0.137333786457509
1.21 -0.137921380817345
1.22 -0.138494783763945
1.23 -0.13905403385454
1.24 -0.139599167844449
1.25 -0.140130220706354
1.26 -0.140647225647601
1.27 -0.141150214125578
1.28 -0.141639215861175
1.29 -0.142114258850364
1.3 -0.142575369373904
1.31 -0.143022572005179
1.32 -0.14345588961619
1.33 -0.143875343381676
1.34 -0.144280952781387
1.35 -0.144672735600469
1.36 -0.145050707927976
1.37 -0.145414884153466
1.38 -0.145765276961661
1.39 -0.146101897325144
1.4 -0.146424754495054
1.41 -0.146733855989724
1.42 -0.147029207581233
1.43 -0.147310813279796
1.44 -0.147578675315933
1.45 -0.147832794120357
1.46 -0.148073168301488
1.47 -0.148299794620509
1.48 -0.148512667963874
1.49 -0.148711781313156
1.5 -0.148897125712118
1.51 -0.14906869023089
1.52 -0.149226461927093
1.53 -0.149370425803786
1.54 -0.149500564764042
1.55 -0.149616859561994
1.56 -0.149719288750141
1.57 -0.149807828622706
1.58 -0.14988245315481
1.59 -0.149943133937204
1.6 -0.149989840106281
1.61 -0.15002253826907
1.62 -0.150041192422868
1.63 -0.150045763869162
1.64 -0.150036211121419
1.65 -0.150012489806346
1.66 -0.149974552558097
1.67 -0.149922348904952
1.68 -0.149855825147844
1.69 -0.149774924230152
1.7 -0.149679585598011
1.71 -0.149569745050423
1.72 -0.149445334578271
1.73 -0.149306282191337
1.74 -0.149152511732258
1.75 -0.148983942676273
1.76 -0.148800489915479
1.77 -0.148602063526148
1.78 -0.148388568517529
1.79 -0.148159904560317
1.8 -0.147915965692796
1.81 -0.147656640002388
1.82 -0.147381809280064
1.83 -0.147091348644728
1.84 -0.146785126134317
1.85 -0.146463002259892
1.86 -0.146124829518502
1.87 -0.145770451859958
1.88 -0.145399704101974
1.89 -0.145012411287253
1.9 -0.14460838797513
1.91 -0.144187437459166
1.92 -0.143749350900671
1.93 -0.143293906366442
1.94 -0.14282086775691
1.95 -0.142329983608405
1.96 -0.141820985750181
1.97 -0.141293587793067
1.98 -0.140747483421975
1.99 -0.140182344458671
2 -0.139597818653945
2.01 -0.138993527159077
2.02 -0.138369061614714
2.03 -0.137723980780074
2.04 -0.137057806605584
2.05 -0.136370019625928
2.06 -0.135660053515646
2.07 -0.13492728860219
2.08 -0.134171044066478
2.09 -0.133390568470214
2.1 -0.132585028119703
2.11 -0.13175349258704
2.12 -0.130894916426708
2.13 -0.130008115689577
2.14 -0.129091737139279
2.15 -0.128144216913288
2.16 -0.127163723325311
2.17 -0.126148074651676
2.18 -0.125094614782291
2.19 -0.124000010800473
2.2 -0.122859880846821
2.21 -0.121667875989532
2.22 -0.120410848269444
2.23 -0.119062578448354
2.24 -0.117583719592638
2.25 -0.115930863571774
2.26 -0.11406963080979
2.27 -0.111984286611436
2.28 -0.109679986498163
2.29 -0.107178579288539
2.3 -0.104511405349332
2.31 -0.101712359294165
2.32 -0.0988130074342417
2.33 -0.0958400906251268
2.34 -0.0928148975483673
2.35 -0.0897537528565766
2.36 -0.0866689786307045
2.37 -0.083569917778005
2.38 -0.0804638164120148
2.39 -0.077356501410568
2.4 -0.0742528625616373
2.41 -0.0711571769287758
2.42 -0.0680733169540109
2.43 -0.0650048776036139
2.44 -0.0619552493369201
2.45 -0.0589276561609548
2.46 -0.0559251723874567
2.47 -0.052950727733191
2.48 -0.0500071076094267
2.49 -0.0470969533986172
2.5 -0.0442227659396124
2.51 -0.0413869142057551
2.52 -0.0385916502261061
2.53 -0.0358391306781974
2.54 -0.033131445297027
2.55 -0.0304706523360425
2.56 -0.0278588218500439
2.57 -0.0252980887079583
2.58 -0.0227907193839206
2.59 -0.0203392007498345
2.6 -0.0179463682772143
2.61 -0.0156156160004412
2.62 -0.0133513372498387
2.63 -0.0111604825092085
2.64 -0.00905430942240314
2.65 -0.00704793701266023
2.66 -0.00515692991005355
2.67 -0.00339270780721649
2.68 -0.00175930430091183
2.69 -0.000252765355998649
2.7 0.00113716537384681
2.71 0.00242429417429448
2.72 0.00362332813309802
2.73 0.0047479849771518
2.74 0.00580999625274736
2.75 0.00681882582425502
2.76 0.00778182835485512
2.77 0.00870460660516968
2.78 0.00959140763854833
2.79 0.0104454750876313
2.8 0.0112693281597568
2.81 0.0120649676162554
2.82 0.0128340213672785
2.83 0.0135778449308846
2.84 0.014297590168152
2.85 0.0149942525396116
2.86 0.0156687041022447
2.87 0.0163217170875407
2.88 0.0169539812126348
2.89 0.0175661167494485
2.9 0.0181586846521241
2.91 0.0187321945866492
2.92 0.0192871114210816
2.93 0.0198238605553394
2.94 0.0203428323550891
2.95 0.0208443858796584
2.96 0.0213288520439605
2.97 0.0217965363200062
2.98 0.0222477210592118
2.99 0.0226826674989977
3 0.0231016175039977
3.01 0.02350479508222
3.02 0.0238924077087928
3.03 0.0242646474839094
3.04 0.0246216921468034
3.05 0.0249637059637613
3.06 0.0252908405050654
3.07 0.0256032353232167
3.08 0.0259010185426761
3.09 0.0261843073695884
3.1 0.0264532085284538
3.11 0.0267078186314182
3.12 0.0269482244847214
3.13 0.0271745033358381
3.14 0.0273867230639344
3.15 0.0275849423154134
3.16 0.0277692105855198
3.17 0.0279395682461826
3.18 0.0280960465194877
3.19 0.0282386673953517
3.2 0.0283674434911056
3.21 0.0284823778497553
3.22 0.0285834636726433
3.23 0.0286706839810522
3.24 0.0287440111999277
3.25 0.0288034066553039
3.26 0.0288488199751197
3.27 0.0288801883808411
3.28 0.0288974358545421
3.29 0.028900472162704
3.3 0.0288891917137825
3.31 0.0288634722213081
3.32 0.0288231731375852
3.33 0.0287681338144513
3.34 0.0286981713363806
3.35 0.0286130779564861
3.36 0.0285126180462814
3.37 0.0283965244433031
3.38 0.0282644940436724
3.39 0.0281161824344348
3.4 0.0279511972851695
3.41 0.0277690901069727
3.42 0.0275693458176186
3.43 0.0273513692860374
3.44 0.027114467596526
3.45 0.0268578260359968
3.46 0.0265804744801273
3.47 0.0262812382848912
3.48 0.0259586622976459
3.49 0.0256108829700795
3.5 0.0252353794138686
3.51 0.0248282049177635
3.52 0.0243810845346417
3.53 0.0238782432517526
3.54 0.0232971779593694
3.55 0.0226144289064114
3.56 0.0218130539768653
3.57 0.020887634148331
3.58 0.0198449628694402
3.59 0.0187012666810649
3.6 0.0174780190205056
3.61 0.0161981258133181
3.62 0.0148833406300962
3.63 0.0135529582373588
3.64 0.0122234298810992
3.65 0.0109084881289137
3.66 0.00961948991216089
3.67 0.00836584056948268
3.68 0.00715549324653897
3.69 0.00599590806792423
3.7 0.00489540325655775
3.71 0.00386371893220712
3.72 0.00291082073121183
3.73 0.00204444111626869
3.74 0.00126784812835137
3.75 0.000579016358879906
3.76 -2.86308940217638e-05
3.77 -0.000564520813240939
3.78 -0.00103907273780243
3.79 -0.00146224093727336
3.8 -0.00184265965673781
3.81 -0.00218732862043022
3.82 -0.00250166299193189
3.83 -0.00278972864169112
3.84 -0.0030545320722186
3.85 -0.00329829069895016
3.86 -0.00352265258788922
3.87 -0.00372886086794185
3.88 -0.00391787029115285
3.89 -0.00409042705584347
3.9 -0.00424712236133064
3.91 -0.0043884279623054
3.92 -0.00451471963792425
3.93 -0.00462629253896969
3.94 -0.00472337093751539
3.95 -0.00480611390411715
3.96 -0.00487461775375373
3.97 -0.00492891561480395
3.98 -0.00496897408534472
3.99 -0.00499468655984785
4 -0.00500586234017943
4.01 -0.00500220994093591
4.02 -0.00498331176311708
4.03 -0.00494858478825098
4.04 -0.00489721577284166
4.05 -0.0048280395718404
4.06 -0.00473920310910489
4.07 -0.00462659302116557
4.08 -0.00448197639885349
4.09 -0.00429339390601551
4.1 -0.00404873261354278
4.11 -0.00374073849613213
4.12 -0.00337076687540189
4.13 -0.00294965143319744
4.14 -0.00249602406533431
4.15 -0.00203299094063877
4.16 -0.00158395814227683
4.17 -0.00116862257370337
4.18 -0.000800267268388832
4.19 -0.0004850022357042
4.2 -0.000222760386339894
4.21 -9.27911401398841e-06
4.22 0.000161778236301534
4.23 0.000297126219158986
4.24 0.000402786397028366
4.25 0.000483573674917611
4.26 0.000542943040607302
4.27 0.000582935227342034
4.28 0.000603465116128717
4.29 0.000601974342361609
4.3 0.000575288091708814
4.31 0.000522975284975664
4.32 0.000449533877105367
4.33 0.0003638887448251
4.34 0.000276770207559032
4.35 0.000197607400948376
4.36 0.000132412876156298
4.37 8.32666232837481e-05
4.38 4.91094272188321e-05
4.39 2.71106609289791e-05
4.4 1.39413615726947e-05
4.41 6.60918895361788e-06
4.42 2.82365352803796e-06
4.43 1.02758751240948e-06
4.44 2.61288121196303e-07
4.45 -1.68425299446642e-08
4.46 -8.737913887636e-08
4.47 -8.29810219243484e-08
4.48 -5.97627181770654e-08
4.49 -3.78110421442576e-08
4.5 -2.20852847266273e-08
4.51 -1.2190646791431e-08
4.52 -6.44320527036814e-09
4.53 -3.28806665629885e-09
4.54 -1.62937213546471e-09
4.55 -7.87309614136738e-10
4.56 -3.72127800441438e-10
4.57 -1.72483083759301e-10
4.58 -7.85581421976902e-11
4.59 -3.52175489229059e-11
4.6 -1.55621677240674e-11
4.61 -6.78665971692594e-12
4.62 -2.92401281308481e-12
4.63 -1.2457915367945e-12
4.64 -5.25307454584175e-13
4.65 -2.19383602657171e-13
4.66 -9.08039341824997e-14
4.67 -3.72712957380404e-14
4.68 -1.51791991035352e-14
4.69 -6.13679935795949e-15
4.7 -2.46404937147605e-15
4.71 -9.829937352839e-16
4.72 -3.89772310175263e-16
4.73 -1.53667587174339e-16
4.74 -6.0256689871548e-17
4.75 -2.35077074504306e-17
4.76 -9.12679394311429e-18
4.77 -3.52729402184857e-18
4.78 -1.35732960897151e-18
4.79 -5.20172050364575e-19
4.8 -1.98572047722038e-19
4.81 -7.55238536462722e-20
4.82 -2.86236640720231e-20
4.83 -1.08122615784316e-20
4.84 -4.07125874597649e-21
4.85 -1.52836857979795e-21
4.86 -5.72106951090014e-22
4.87 -2.1356739577289e-22
4.88 -7.95164174877368e-23
4.89 -2.95321583659068e-23
4.9 -1.09420729279728e-23
4.91 -4.044984421059e-24
4.92 -1.49207691141227e-24
4.93 -5.49243082096781e-25
4.94 -2.01779353545224e-25
4.95 -7.39886760854365e-26
4.96 -2.70810754130321e-26
4.97 -9.89492677676465e-27
4.98 -3.60942063067828e-27
4.99 -1.31453225552277e-27
5 -4.78016555773906e-28
5.01 -1.73572234741565e-28
5.02 -6.29375087190035e-29
5.03 -2.27906079618673e-29
5.04 -8.24220190610632e-30
5.05 -2.97710635486478e-30
5.06 -1.07406592929967e-30
5.07 -3.87055950235689e-31
5.08 -1.3932933166704e-31
5.09 -5.01021360520828e-32
5.1 -1.79983564467071e-32
5.11 -6.45936172037032e-33
5.12 -2.31602360895078e-33
5.13 -8.29676376659266e-34
5.14 -2.96962692456244e-34
5.15 -1.06203026223326e-34
5.16 -3.79513747049518e-35
5.17 -1.35514869268488e-35
5.18 -4.83534923206773e-36
5.19 -1.72409941819724e-36
5.2 -6.14330285174962e-37
5.21 -2.18754916610772e-37
5.22 -7.78467797144738e-38
5.23 -2.76860264488919e-38
5.24 -9.84073538948657e-39
5.25 -3.49583418894317e-39
5.26 -1.24119341432303e-39
5.27 -4.40455585380729e-40
5.28 -1.56223769182662e-40
5.29 -5.53837567384527e-41
5.3 -1.96252718463358e-41
5.31 -6.95111192657201e-42
5.32 -2.46096412178529e-42
5.33 -8.70914411023114e-43
5.34 -3.08085596530874e-43
5.35 -1.08942995471376e-43
5.36 -3.8509267993857e-44
5.37 -1.36073968251153e-44
5.38 -4.8065587220687e-45
5.39 -1.69725946018152e-45
5.4 -5.99131522832266e-46
5.41 -2.11427277463993e-46
5.42 -7.45881027889141e-47
5.43 -2.63058537566462e-47
5.44 -9.27500108427658e-48
5.45 -3.26932825519993e-48
5.46 -1.15210001864999e-48
5.47 -4.05894241669683e-49
5.48 -1.42965239148841e-49
5.49 -5.03438640790927e-50
5.5 -1.77241206090317e-50
5.51 -6.23861679912698e-51
5.52 -2.19543542234967e-51
5.53 -7.72440306314464e-52
5.54 -2.71721583773789e-52
5.55 -9.55655211620882e-53
5.56 -3.36046242075141e-53
5.57 -1.18146346772382e-53
5.58 -4.15305445078828e-54
5.59 -1.45963262890895e-54
5.6 -5.12921090157777e-55
5.61 -1.80215011690531e-55
5.62 -6.330923430553e-56
5.63 -2.22372480545373e-56
5.64 -7.80971202059492e-57
5.65 -2.74240151548942e-57
5.66 -9.62877574139063e-58
5.67 -3.38031377668202e-58
5.68 -1.18656273474413e-58
5.69 -4.16460589620388e-59
5.7 -1.46153203195404e-59
5.71 -5.12856233412739e-60
5.72 -1.79944013011354e-60
5.73 -6.31299090182449e-61
5.74 -2.21457506775201e-61
5.75 -7.76791672397431e-62
5.76 -2.72445105145117e-62
5.77 -9.55465666313749e-63
5.78 -3.35053437073798e-63
5.79 -1.17483595206894e-63
5.8 -4.11913228176212e-64
5.81 -1.44411164384367e-64
5.82 -5.06248139307906e-65
5.83 -1.77457693136903e-65
5.84 -6.22008026893251e-66
5.85 -2.18005736641282e-66
5.86 -7.64032097885126e-67
5.87 -2.67749074087656e-67
5.88 -9.3824865768858e-68
5.89 -3.2876259901063e-68
5.9 -1.1519196490795e-68
5.91 -4.03587880135647e-69
5.92 -1.41394001544882e-69
5.93 -4.95337902766454e-70
5.94 -1.7352041676792e-70
5.95 -6.07825286432024e-71
5.96 -2.12905458248139e-71
5.97 -7.45719236000231e-72
5.98 -2.61183092349633e-72
5.99 -9.14737660074202e-73
6 -3.20354219245977e-73
};
\addlegendentry{$\beta(0, t)=\alpha(0, t)$}
\addplot [black, dash pattern=on 1pt off 3pt on 3pt off 3pt, dash phase = 12pt] 
table {%
0 0
0.01 0
0.02 0
0.03 0
0.04 0
0.05 0
0.06 0
0.07 0
0.08 0
0.09 0
0.1 0
0.11 0
0.12 0
0.13 0
0.14 0
0.15 0
0.16 0
0.17 0
0.18 0
0.19 0
0.2 0
0.21 0
0.22 0
0.23 0
0.24 0
0.25 0
0.26 0
0.27 0
0.28 0
0.29 0
0.3 0
0.31 0
0.32 0
0.33 0
0.34 0
0.35 0
0.36 0
0.37 0
0.38 0
0.39 0
0.4 0
0.41 -5.28386011161147e-08
0.42 -6.02499735752073e-07
0.43 -3.55114889575731e-06
0.44 -1.4425882363517e-05
0.45 -4.54441255562664e-05
0.46 -0.000118435605035775
0.47 -0.000266067934793598
0.48 -0.000530124140087158
0.49 -0.000956643726187997
0.5 -0.00158900981088997
0.51 -0.00246088353619269
0.52 -0.00359090000518518
0.53 -0.00498034344999758
0.54 -0.00661401295286434
0.55 -0.00846361101207019
0.56 -0.0104925063155539
0.57 -0.0126606835217078
0.58 -0.0149289736341158
0.59 -0.0172620737107193
0.6 -0.0196302550193568
0.61 -0.0220099333067892
0.62 -0.0243834111702141
0.63 -0.0267381242465441
0.64 -0.0290656727423327
0.65 -0.0313608391395865
0.66 -0.0336207115928304
0.67 -0.0358439668072439
0.68 -0.0380303212761482
0.69 -0.040180134048214
0.7 -0.042294133017439
0.71 -0.044373234914123
0.72 -0.0464184324024807
0.73 -0.0484307269288266
0.74 -0.050411091393178
0.75 -0.0523604514335273
0.76 -0.054279677806366
0.77 -0.0561695850381334
0.78 -0.0580309333749679
0.79 -0.05986443227486
0.8 -0.0616707444526368
0.81 -0.0634504899522378
0.82 -0.0652042499908323
0.83 -0.066932570469844
0.84 -0.0686359651274849
0.85 -0.0703149183463462
0.86 -0.0719698876462756
0.87 -0.0736013058978273
0.88 -0.075209583290985
0.89 -0.0767951090908286
0.9 -0.0783582532080019
0.91 -0.0798993676080303
0.92 -0.0814187875800647
0.93 -0.0829168328826021
0.94 -0.0843938087811623
0.95 -0.0858500069907264
0.96 -0.0872857065339258
0.97 -0.0887011745244442
0.98 -0.0900966668838113
0.99 -0.0914724289986892
1 -0.0928286963248332
1.01 -0.0941656949431374
1.02 -0.0954836420725054
1.03 -0.0967827465437259
1.04 -0.0980632092380367
1.05 -0.0993252234936454
1.06 -0.100568975483103
1.07 -0.101794644564118
1.08 -0.103002403606106
1.09 -0.104192419294548
1.1 -0.105364852415
1.11 -0.106519858118414
1.12 -0.107657586169277
1.13 -0.108778181177906
1.14 -0.109881782818127
1.15 -0.110968526031447
1.16 -0.112038541218715
1.17 -0.113091954420188
1.18 -0.114128887484833
1.19 -0.115149458229611
1.2 -0.116153780589455
1.21 -0.117141964758548
1.22 -0.118114117323492
1.23 -0.119070341388902
1.24 -0.120010736695893
1.25 -0.120935399733925
1.26 -0.121844423846395
1.27 -0.122737899330379
1.28 -0.123615913530843
1.29 -0.124478550929662
1.3 -0.125325893229735
1.31 -0.126158019434468
1.32 -0.126975005922877
1.33 -0.127776926520542
1.34 -0.128563852566627
1.35 -0.129335852977168
1.36 -0.130092994304798
1.37 -0.130835340795102
1.38 -0.131562954439725
1.39 -0.132275895026417
1.4 -0.132974220186109
1.41 -0.133657985437171
1.42 -0.13432724422695
1.43 -0.134982047970694
1.44 -0.135622446087966
1.45 -0.13624848603662
1.46 -0.136860213344429
1.47 -0.137457671638441
1.48 -0.138040902672111
1.49 -0.138609946350286
1.5 -0.13916484075209
1.51 -0.139705622151747
1.52 -0.140232325037399
1.53 -0.140744982127936
1.54 -0.14124362438788
1.55 -0.141728281040348
1.56 -0.142198979578104
1.57 -0.142655745772718
1.58 -0.143098603681844
1.59 -0.143527575654617
1.6 -0.143942682335166
1.61 -0.144343942664244
1.62 -0.144731373878947
1.63 -0.145104991510526
1.64 -0.145464809380246
1.65 -0.145810839593284
1.66 -0.146143092530615
1.67 -0.146461576838856
1.68 -0.146766299418019
1.69 -0.147057265407112
1.7 -0.147334478167536
1.71 -0.147597939264212
1.72 -0.147847648444348
1.73 -0.148083603613781
1.74 -0.148305800810783
1.75 -0.148514234177247
1.76 -0.148708895927121
1.77 -0.148889776311989
1.78 -0.149056863583639
1.79 -0.149210143953499
1.8 -0.149349601548751
1.81 -0.149475218364975
1.82 -0.149586974215103
1.83 -0.149684846674504
1.84 -0.149768811021944
1.85 -0.14983884017619
1.86 -0.149894904627977
1.87 -0.149936972367044
1.88 -0.149965008803915
1.89 -0.149978976686068
1.9 -0.149978836008097
1.91 -0.149964543915461
1.92 -0.149936054601327
1.93 -0.149893319196012
1.94 -0.149836285648456
1.95 -0.149764898599102
1.96 -0.149679099243508
1.97 -0.149578825185921
1.98 -0.149464010282005
1.99 -0.149334584469768
2 -0.149190473587693
2.01 -0.149031599178919
2.02 -0.148857878280206
2.03 -0.148669223194283
2.04 -0.148465541243985
2.05 -0.148246734506418
2.06 -0.148012699525159
2.07 -0.147763326998266
2.08 -0.147498501439544
2.09 -0.147218100810245
2.1 -0.146921996117913
2.11 -0.146610050978711
2.12 -0.146282121138987
2.13 -0.14593805395124
2.14 -0.145577687798912
2.15 -0.145200851463585
2.16 -0.144807363427107
2.17 -0.144397031099977
2.18 -0.14396964996584
2.19 -0.143525002630148
2.2 -0.143062857758914
2.21 -0.142582968890845
2.22 -0.142085073102882
2.23 -0.141568889505178
2.24 -0.141034117536532
2.25 -0.140480435024982
2.26 -0.139907495970269
2.27 -0.139314927994589
2.28 -0.138702329394774
2.29 -0.138069265711566
2.3 -0.137415265708487
2.31 -0.136739816621493
2.32 -0.136042358497774
2.33 -0.135322277382093
2.34 -0.134578897023446
2.35 -0.13381146864922
2.36 -0.13301915816395
2.37 -0.132201029830617
2.38 -0.131356024997342
2.39 -0.130482933553794
2.4 -0.129580354063305
2.41 -0.128646634135185
2.42 -0.127679761618238
2.43 -0.126677126325207
2.44 -0.125635026675886
2.45 -0.124547843843933
2.46 -0.12340699829987
2.47 -0.122200060758682
2.48 -0.120910530887821
2.49 -0.119518681503126
2.5 -0.118003511778339
2.51 -0.116345431654177
2.52 -0.114529019061112
2.53 -0.112545170537427
2.54 -0.110392184222336
2.55 -0.108075651003208
2.56 -0.105607342554767
2.57 -0.103003477205169
2.58 -0.100282789400367
2.59 -0.0974647566226967
2.6 -0.0945682058302189
2.61 -0.0916103841448937
2.62 -0.0886064711220888
2.63 -0.0855694460530377
2.64 -0.082510200352085
2.65 -0.0794377905747708
2.66 -0.0763597487384359
2.67 -0.0732823925550731
2.68 -0.0702111020377749
2.69 -0.0671505474568035
2.7 -0.0641048661413669
2.71 -0.0610777929342641
2.72 -0.0580727525432501
2.73 -0.0550929229300913
2.74 -0.0521412783160592
2.75 -0.0492206191255884
2.76 -0.0463335947262641
2.77 -0.0434827234399537
2.78 -0.040670413157317
2.79 -0.0378989850852944
2.8 -0.0351707027950896
2.81 -0.0324878090218774
2.82 -0.0298525742963889
2.83 -0.0272673683731117
2.84 -0.0247347800055012
2.85 -0.0222578199723384
2.86 -0.0198402192390871
2.87 -0.0174867666148454
2.88 -0.0152035484123514
2.89 -0.0129979165234675
2.9 -0.0108780641773119
2.91 -0.00885221894499254
2.92 -0.00692760882535197
2.93 -0.00510944844881686
2.94 -0.00340019007813396
2.95 -0.00179919804825456
2.96 -0.000302879864955096
2.97 0.00109480685900259
2.98 0.00240162098360847
2.99 0.00362624141555267
3 0.00477754243814175
3.01 0.00586402497430036
3.02 0.00689343553715331
3.03 0.0078725635324317
3.04 0.00880718294188396
3.05 0.00970209518992923
3.06 0.0105612318839339
3.07 0.0113877841266156
3.08 0.0121843350900967
3.09 0.0129529818274091
3.1 0.0136954395734733
3.11 0.014413126754145
3.12 0.0151072318564337
3.13 0.0157787647259532
3.14 0.0164285952733404
3.15 0.0170574824198725
3.16 0.0176660956982533
3.17 0.0182550314371082
3.18 0.0188248249993507
3.19 0.0193759601594005
3.2 0.0199088764023338
3.21 0.0204239747023242
3.22 0.0209216221744543
3.23 0.0214021558785179
3.24 0.0218658859729397
3.25 0.0223130983612159
3.26 0.0227440569347116
3.27 0.023159005488807
3.28 0.0235581693705013
3.29 0.0239417569020921
3.3 0.024309960615743
3.31 0.0246629583264803
3.32 0.0250009140656592
3.33 0.0253239788926951
3.34 0.0256322915995168
3.35 0.025925979319525
3.36 0.0262051580506593
3.37 0.026469933100365
3.38 0.0267203994587264
3.39 0.0269566421047102
3.4 0.0271787362493066
3.41 0.0273867475183066
3.42 0.0275807320764902
3.43 0.0277607366940767
3.44 0.0279267987553884
3.45 0.0280789462087674
3.46 0.0282171974558347
3.47 0.0283415611771677
3.48 0.0284520360903498
3.49 0.0285486106350935
3.5 0.0286312625786858
3.51 0.0286999585333197
3.52 0.0287546533748554
3.53 0.0287952895501294
3.54 0.0288217962569567
3.55 0.0288340884772846
3.56 0.0288320658393357
3.57 0.0288156112786957
3.58 0.0287845894607215
3.59 0.0287388449167315
3.6 0.0286781998332638
3.61 0.0286024514158719
3.62 0.0285113687243616
3.63 0.0284046888417449
3.64 0.0282821121891496
3.65 0.0281432967244775
3.66 0.0279878506479656
3.67 0.0278153230536635
3.68 0.0276251916534297
3.69 0.0274168461278599
3.7 0.0271895644700257
3.71 0.0269424761231166
3.72 0.0266744916836074
3.73 0.0263841500084904
3.74 0.0260693143670402
3.75 0.0257266888602439
3.76 0.025351244316154
3.77 0.0249357839450862
3.78 0.0244709394635927
3.79 0.0239457984815427
3.8 0.0233491487161746
3.81 0.0226710894879655
3.82 0.0219046207982792
3.83 0.0210468313970235
3.84 0.0200994468163911
3.85 0.0190686921836733
3.86 0.0179645939275434
3.87 0.0167999409554957
3.88 0.0155891417172295
3.89 0.0143471719619684
3.9 0.0130887436788852
3.91 0.0118277628827261
3.92 0.0105770851439219
3.93 0.00934851544821665
3.94 0.00815293808469715
3.95 0.00700042781985625
3.96 0.00590021225783158
3.97 0.00486042923865398
3.98 0.00388772227798074
3.99 0.00298679608144893
4 0.00216007956449247
4.01 0.00140761062812053
4.02 0.000727186624639472
4.03 0.000114749429392114
4.04 -0.000435078974213129
4.05 -0.000928413745897914
4.06 -0.00137156926944926
4.07 -0.00177062807682247
4.08 -0.00213114367161142
4.09 -0.00245797465225521
4.1 -0.0027552282333969
4.11 -0.0030262827623984
4.12 -0.00327385875571099
4.13 -0.00350011298306671
4.14 -0.00370673715392549
4.15 -0.00389504966489452
4.16 -0.00406607449179726
4.17 -0.00422060527996829
4.18 -0.00435925511199513
4.19 -0.00448249364911761
4.2 -0.00459067373059849
4.21 -0.00468404939959146
4.22 -0.00476278693568596
4.23 -0.00482696994495951
4.24 -0.00487659891810235
4.25 -0.00491158479775649
4.26 -0.00493173429623871
4.27 -0.00493671745129486
4.28 -0.00492599174352077
4.29 -0.00489864504956063
4.3 -0.00485314108466257
4.31 -0.00478702087923151
4.32 -0.00469670266583147
4.33 -0.00457756625066399
4.34 -0.00442445703743327
4.35 -0.00423260619927808
4.36 -0.00399879693266836
4.37 -0.00372248664245752
4.38 -0.0034065715643719
4.39 -0.00305756241085191
4.4 -0.00268509670365918
4.41 -0.00230089048097834
4.42 -0.00191737032875806
4.43 -0.00154628538589535
4.44 -0.00119756900258538
4.45 -0.000878623400817479
4.46 -0.000594078516445706
4.47 -0.000345968857861557
4.48 -0.000134208346575904
4.49 4.27702926630696e-05
4.5 0.000187329201074814
4.51 0.000301998290420379
4.52 0.000389016836834544
4.53 0.000450144057099453
4.54 0.000486782842331723
4.55 0.000500358744221063
4.56 0.000492794355939657
4.57 0.000466876139677801
4.58 0.000426354595051083
4.59 0.000375728194885078
4.6 0.000319782333331085
4.61 0.000263033238569346
4.62 0.000209238257198219
4.63 0.000161087326349098
4.64 0.000120116610298957
4.65 8.68171972317804e-05
4.66 6.08703703522286e-05
4.67 4.143095959393e-05
4.68 2.73940880791298e-05
4.69 1.7605777979424e-05
4.7 1.10032134398626e-05
4.71 6.68926364902879e-06
4.72 3.95610583956565e-06
4.73 2.27563913671731e-06
4.74 1.27243723224013e-06
4.75 6.90860319957036e-07
4.76 3.63543331337391e-07
4.77 1.84848321806386e-07
4.78 9.03682515087957e-08
4.79 4.2123726155029e-08
4.8 1.84415791471569e-08
4.81 7.35424788092613e-09
4.82 2.47366962011487e-09
4.83 5.10827984214235e-10
4.84 -1.61184759249823e-10
4.85 -3.1008574916929e-10
4.86 -2.77079231491413e-10
4.87 -2.0138975532273e-10
4.88 -1.32099456956149e-10
4.89 -8.12828186062906e-11
4.9 -4.78095898080834e-11
4.91 -2.71703439868394e-11
4.92 -1.50195448777184e-11
4.93 -8.11280740160753e-12
4.94 -4.29581997980203e-12
4.95 -2.23527922587845e-12
4.96 -1.14509127551696e-12
4.97 -5.78385125487725e-13
4.98 -2.88396293222286e-13
4.99 -1.42100752278087e-13
5 -6.92482500382439e-14
5.01 -3.34000041152843e-14
5.02 -1.59546745575988e-14
5.03 -7.55228542211232e-15
5.04 -3.54434585706213e-15
5.05 -1.6498971101153e-15
5.06 -7.62110078645072e-16
5.07 -3.49445678554184e-16
5.08 -1.59107315452942e-16
5.09 -7.1958927727415e-17
5.1 -3.23361659580187e-17
5.11 -1.44416847249479e-17
5.12 -6.4118131328596e-18
5.13 -2.83060029404427e-18
5.14 -1.24281237853185e-18
5.15 -5.42813708474624e-19
5.16 -2.35883700190693e-19
5.17 -1.02006039328838e-19
5.18 -4.39045256316678e-20
5.19 -1.88112999893315e-20
5.2 -8.02456455379306e-21
5.21 -3.40863864299577e-21
5.22 -1.44197539394352e-21
5.23 -6.07588905355153e-22
5.24 -2.55030397082894e-22
5.25 -1.0664898407627e-22
5.26 -4.44379865189671e-23
5.27 -1.84515532814376e-23
5.28 -7.63551602185472e-24
5.29 -3.14930580428287e-24
5.3 -1.29480547958077e-24
5.31 -5.30698737020547e-25
5.32 -2.16862461795407e-25
5.33 -8.8359037548806e-26
5.34 -3.58991248137603e-26
5.35 -1.45451333643127e-26
5.36 -5.877408834823e-27
5.37 -2.36875752282755e-27
5.38 -9.52254047230674e-28
5.39 -3.81867294440971e-28
5.4 -1.52766584368597e-28
5.41 -6.09717121367487e-29
5.42 -2.42794825351612e-29
5.43 -9.64689301476952e-30
5.44 -3.82470209899438e-30
5.45 -1.51319245818154e-30
5.46 -5.97448692327893e-31
5.47 -2.35418010556473e-31
5.48 -9.25834546279709e-32
5.49 -3.63415134503536e-32
5.5 -1.42386627907526e-32
5.51 -5.56867615719999e-33
5.52 -2.17405567529846e-33
5.53 -8.47313431283016e-34
5.54 -3.2967838184477e-34
5.55 -1.28064121241478e-34
5.56 -4.96674944865433e-35
5.57 -1.92327528286411e-35
5.58 -7.43620282847252e-36
5.59 -2.87089476782246e-36
5.6 -1.10676350861173e-36
5.61 -4.26067793097179e-37
5.62 -1.63795998405039e-37
5.63 -6.28843474812957e-38
5.64 -2.41107164566792e-38
5.65 -9.2324970612511e-39
5.66 -3.53087711385207e-39
5.67 -1.34869235997088e-39
5.68 -5.1454365842284e-40
5.69 -1.96075099052108e-40
5.7 -7.46319970051957e-41
5.71 -2.83753553877688e-41
5.72 -1.07766085892288e-41
5.73 -4.08844445625297e-42
5.74 -1.54945743708315e-42
5.75 -5.86619981545205e-43
5.76 -2.21870489771728e-43
5.77 -8.3833443895085e-44
5.78 -3.16460486326189e-44
5.79 -1.19348015215198e-44
5.8 -4.49690089578952e-45
5.81 -1.69286566379108e-45
5.82 -6.36723947852507e-46
5.83 -2.39280670535036e-46
5.84 -8.98462238066334e-47
5.85 -3.37082074506544e-47
5.86 -1.26363822308553e-47
5.87 -4.73334959408284e-48
5.88 -1.77166018838431e-48
5.89 -6.62621311147165e-49
5.9 -2.47645502415676e-49
5.91 -9.24873512514682e-50
5.92 -3.45165761811585e-50
5.93 -1.28727997873464e-50
5.94 -4.79760590696664e-51
5.95 -1.78685185018757e-51
5.96 -6.6507552788032e-52
5.97 -2.47387483877401e-52
5.98 -9.19633041589976e-53
5.99 -3.41654480119171e-53
6 -1.2685303979962e-53
};
\addlegendentry{$\alpha(1, t)$}
\addplot [semithick, black, mark=*, mark size=1.25, mark options={solid}, forget plot]
table {%
4.76 0
};
\draw (axis cs:4.76,0.00675) node[
  anchor=base west,
  text=black,
  rotate=0.0
]{$T_1$};
\end{axis}

\end{tikzpicture}

%% file: pics/simalphabetaarticle_log.tex
\begin{tikzpicture}

\begin{axis}[
width=\plottimewidth,
height=\plottimeheight,
legend cell align={left},
legend style={fill opacity=0.8, draw opacity=1, text opacity=1, at={(0.97,0.03)}, anchor=south east, draw=white!80!black},
legend pos=south west,
tick align=inside,
tick pos=both,
x grid style={white!69.0196078431373!black},
xlabel={$t$},
xmajorgrids,
xmin=0, xmax=6,
xtick style={color=black},
y grid style={white!69.0196078431373!black},
ymajorgrids,
ymin=-30, ymax=0,
ytick style={color=black},
ytick={-30,-20,-10,0},
yticklabels={-30,-20,-10,0}
]
\addplot [black]
table {%
0 -inf
0.01 -5.77559782344822
0.02 -5.10596224811016
0.03 -4.7205004748999
0.04 -4.45058544333592
0.05 -4.24362216614628
0.06 -4.07627973880481
0.07 -3.93615850862503
0.08 -3.81588180310542
0.09 -3.71070625996797
0.1 -3.61740251335355
0.11 -3.53367345834357
0.12 -3.45782737095472
0.13 -3.38858269191024
0.14 -3.32494558492543
0.15 -3.26612997915373
0.16 -3.21150355960472
0.17 -3.1605502203092
0.18 -3.11284330726087
0.19 -3.0680261344979
0.2 -3.02579752462383
0.21 -2.98590089619545
0.22 -2.94811590348313
0.23 -2.91225194477626
0.24 -2.8781430599239
0.25 -2.84564387529601
0.26 -2.81462634855568
0.27 -2.78497713129923
0.28 -2.75659541411742
0.29 -2.72939115203062
0.3 -2.70328359256195
0.31 -2.67820004662434
0.32 -2.65407485574519
0.33 -2.63084851920286
0.34 -2.60846695229156
0.35 -2.58688085279488
0.36 -2.56604515728654
0.37 -2.54591857241649
0.38 -2.52646316912256
0.39 -2.50764402990994
0.4 -2.489428941095
0.41 -2.47178812331652
0.42 -2.454693994752
0.43 -2.43812096239691
0.44 -2.42204523751469
0.45 -2.40644467198084
0.46 -2.39129861275027
0.47 -2.37658777209638
0.48 -2.36229411161847
0.49 -2.34840073830438
0.5 -2.33489181117876
0.51 -2.32175245727211
0.52 -2.30896869581829
0.53 -2.29652736973477
0.54 -2.28441608356385
0.55 -2.27262314715973
0.56 -2.26113752449629
0.57 -2.24994878704883
0.58 -2.2390470712696
0.59 -2.22842303973476
0.6 -2.21806784559051
0.61 -2.20797309996921
0.62 -2.19813084208406
0.63 -2.18853351174392
0.64 -2.17917392405806
0.65 -2.17004524612627
0.66 -2.16114097553132
0.67 -2.15245492047042
0.68 -2.1439811813792
0.69 -2.13571413391683
0.7 -2.12764841319443
0.71 -2.11977889914019
0.72 -2.11210070290575
0.73 -2.1046091542272
0.74 -2.09729978966254
0.75 -2.09016834163495
0.76 -2.08321072821748
0.77 -2.07642304360128
0.78 -2.069801549194
0.79 -2.06334266530068
0.8 -2.05704296334273
0.81 -2.0508991585755
0.82 -2.04490810326741
0.83 -2.03906678030761
0.84 -2.03337229721133
0.85 -2.02782188049512
0.86 -2.02241287039612
0.87 -2.01714271591196
0.88 -2.01200897013957
0.89 -2.00700928589299
0.9 -2.00214141158199
0.91 -1.99740318733446
0.92 -1.99279254134733
0.93 -1.98830748645139
0.94 -1.98394611687712
0.95 -1.97970660520908
0.96 -1.9755871995179
0.97 -1.9715862206592
0.98 -1.96770205973007
0.99 -1.96393317567408
1 -1.96027809302686
1.01 -1.95673539979433
1.02 -1.95330374545693
1.03 -1.94998183909319
1.04 -1.94676844761675
1.05 -1.94366239412128
1.06 -1.94066255632829
1.07 -1.93776786513306
1.08 -1.93497730324454
1.09 -1.93228990391513
1.1 -1.92970474975686
1.11 -1.92722097164053
1.12 -1.92483774767502
1.13 -1.92255430226382
1.14 -1.92036990523657
1.15 -1.91828387105322
1.16 -1.91629555807903
1.17 -1.9144043679287
1.18 -1.91260974487801
1.19 -1.910911175342
1.2 -1.9093081874185
1.21 -1.90780035049621
1.22 -1.90638727492697
1.23 -1.90506861176166
1.24 -1.90384405254966
1.25 -1.90271332920211
1.26 -1.90167621391912
1.27 -1.90073251918158
1.28 -1.89988209780843
1.29 -1.89912484308035
1.3 -1.89846068893121
1.31 -1.89788961020897
1.32 -1.89741162300779
1.33 -1.89702678507363
1.34 -1.89673519628581
1.35 -1.89653699921746
1.36 -1.896432379778
1.37 -1.8964215679415
1.38 -1.89650483856478
1.39 -1.89668251230002
1.4 -1.89695495660692
1.41 -1.89732258687006
1.42 -1.89778586762782
1.43 -1.89834531391982
1.44 -1.8990014927607
1.45 -1.89975502474889
1.46 -1.90060658581977
1.47 -1.90155690915415
1.48 -1.90260678725359
1.49 -1.9037570741958
1.5 -1.9050086880848
1.51 -1.90636261371194
1.52 -1.90781990544612
1.53 -1.90938169037333
1.54 -1.91104917170831
1.55 -1.91282363250381
1.56 -1.91470643968593
1.57 -1.91669904844784
1.58 -1.91880300703822
1.59 -1.92101996198537
1.6 -1.92335166380383
1.61 -1.92579997323621
1.62 -1.92836686809091
1.63 -1.93105445074466
1.64 -1.93386495638928
1.65 -1.93680076211396
1.66 -1.93986439692852
1.67 -1.94305855285018
1.68 -1.94638609719628
1.69 -1.94985008624974
1.7 -1.95345378049333
1.71 -1.95720066164376
1.72 -1.9610944517605
1.73 -1.9651391347564
1.74 -1.96933898070352
1.75 -1.97369857340837
1.76 -1.97822284183343
1.77 -1.98291709607073
1.78 -1.98778706873812
1.79 -1.99283896288044
1.8 -1.99807950773349
1.81 -2.00351602407017
1.82 -2.00915650132915
1.83 -2.01500968937446
1.84 -2.02108520861858
1.85 -2.02739368346981
1.86 -2.03394690579759
1.87 -2.04075803760362
1.88 -2.04784186575592
1.89 -2.05521512717168
1.9 -2.06289693140357
1.91 -2.07090932129423
1.92 -2.07927803512881
1.93 -2.08803357320713
1.94 -2.09721274393166
1.95 -2.10686100507791
1.96 -2.11703621231281
1.97 -2.12781507937675
1.98 -2.13930551124974
1.99 -2.15167405422964
2 -2.16522535527028
2.01 -2.18083542581341
2.02 -2.20489055804359
2.03 -2.23105494368139
2.04 -2.25878507003814
2.05 -2.2879012460841
2.06 -2.31832864925724
2.07 -2.35004054856847
2.08 -2.38303793233206
2.09 -2.41734037248444
2.1 -2.45298137877948
2.11 -2.49000586707575
2.12 -2.52846875986909
2.13 -2.56843427808032
2.14 -2.60997572769132
2.15 -2.65317571027756
2.16 -2.69812675653358
2.17 -2.74493242155805
2.18 -2.79370890152303
2.19 -2.84458724144024
2.2 -2.89771621143562
2.21 -2.95326594376788
2.22 -3.0114324546386
2.23 -3.07244323294918
2.24 -3.13656417157105
2.25 -3.20410825629466
2.26 -3.2754466305063
2.27 -3.35102294995196
2.28 -3.43137238553979
2.29 -3.51714731886097
2.3 -3.60915287711937
2.31 -3.70839728872579
2.32 -3.81616521531725
2.33 -3.93412794938524
2.34 -4.06451524539032
2.35 -4.21039543581494
2.36 -4.3761577349785
2.37 -4.56840203904378
2.38 -4.79773553520226
2.39 -5.08287777725972
2.4 -5.46191435581893
2.41 -6.03413073120179
2.42 -7.27314572191464
2.43 -7.31125810102178
2.44 -6.25917368486169
2.45 -5.78275715540181
2.46 -5.47588118187061
2.47 -5.25116938728663
2.48 -5.07498249289923
2.49 -4.93080902263637
2.5 -4.80932364553249
2.51 -4.70474797100647
2.52 -4.61325634688517
2.53 -4.53218662490609
2.54 -4.45961285950378
2.55 -4.39409752526409
2.56 -4.33453979090005
2.57 -4.2800783993023
2.58 -4.23002716739252
2.59 -4.18383079517933
2.6 -4.14103377112669
2.61 -4.10125798133324
2.62 -4.06418625699783
2.63 -4.02955006770679
2.64 -3.99712016870757
2.65 -3.96669939154242
2.66 -3.93811701542374
2.67 -3.91122432169623
2.68 -3.8858910456669
2.69 -3.86200251742275
2.7 -3.83945733757554
2.71 -3.81816547260307
2.72 -3.79804668245312
2.73 -3.77902921357421
2.74 -3.76104870572096
2.75 -3.74404727225293
2.76 -3.72797272224825
2.77 -3.71277789932277
2.78 -3.69842011710741
2.79 -3.68486067526921
2.8 -3.67206444304152
2.81 -3.65999949965832
2.82 -3.64863682301866
2.83 -3.63795001945174
2.84 -3.6279150886975
2.85 -3.61851021922602
2.86 -3.60971560984219
2.87 -3.60151331419838
2.88 -3.59388710539695
2.89 -3.58682235833089
2.9 -3.58030594780246
2.91 -3.57432616079224
2.92 -3.56887262153613
2.93 -3.56393622831546
2.94 -3.55950910108411
2.95 -3.55558453925338
2.96 -3.552156989136
2.97 -3.54922202072138
2.98 -3.5467763136199
2.99 -3.54481765217998
3 -3.54334492995358
3.01 -3.5423581638691
3.02 -3.54185851867319
3.03 -3.54184834243167
3.04 -3.5423312141461
3.05 -3.54331200485788
3.06 -3.54479695399405
3.07 -3.54679376317874
3.08 -3.54931171032169
3.09 -3.55236178754034
3.1 -3.5559568674294
3.11 -3.56011190344087
3.12 -3.5648441717894
3.13 -3.57017356451622
3.14 -3.57612294637242
3.15 -3.58271859238698
3.16 -3.58999072893679
3.17 -3.59797420973976
3.18 -3.60670937092316
3.19 -3.61624312865147
3.2 -3.62663041302572
3.21 -3.63793608082582
3.22 -3.65023753172677
3.23 -3.66362839676998
3.24 -3.67822393508018
3.25 -3.69416930388061
3.26 -3.71165300519796
3.27 -3.73093054324919
3.28 -3.75237093226429
3.29 -3.7765651434188
3.3 -3.80467038747113
3.31 -3.84123004872315
3.32 -3.90285708862568
3.33 -3.9720687761933
3.34 -4.04833063405369
3.35 -4.13192730307635
3.36 -4.22349181192101
3.37 -4.32388363080952
3.38 -4.43416107789327
3.39 -4.55563471323642
3.4 -4.69002285810183
3.41 -4.8397284541068
3.42 -5.00827455402206
3.43 -5.20101822274556
3.44 -5.42646871357855
3.45 -5.69907977604359
3.46 -6.0461910203987
3.47 -6.52981462539971
3.48 -7.35230361606365
3.49 -13.6457535333334
3.5 -7.59547528708725
3.51 -6.96127853776201
3.52 -6.6062922567237
3.53 -6.36412162401432
3.54 -6.1831546322619
3.55 -6.04062426821012
3.56 -5.92449766676581
3.57 -5.82765684072578
3.58 -5.74555214252576
3.59 -5.67510646557397
3.6 -5.61414745358191
3.61 -5.56108945347792
3.62 -5.51474421678658
3.63 -5.47420263273865
3.64 -5.43875787153661
3.65 -5.40785380660308
3.66 -5.38104949125963
3.67 -5.35799419680797
3.68 -5.3384096309727
3.69 -5.32207720202939
3.7 -5.30882895975036
3.71 -5.29854133597425
3.72 -5.29113114153447
3.73 -5.28655352254007
3.74 -5.28480178457739
3.75 -5.28590919629065
3.76 -5.28995312429862
3.77 -5.29706218578835
3.78 -5.30742762971659
3.79 -5.32132106163463
3.8 -5.33912232713921
3.81 -5.36136489315985
3.82 -5.38881431650336
3.83 -5.42261801692554
3.84 -5.46464175534631
3.85 -5.51848181923249
3.86 -5.59602815678232
3.87 -5.76725212702087
3.88 -5.98715454774832
3.89 -6.26742298565428
3.9 -6.62370457586737
3.91 -7.08967116764585
3.92 -7.77086304743026
3.93 -9.20716432485712
3.94 -9.36817841837958
3.95 -8.38245175711022
3.96 -7.97068126318367
3.97 -7.72639841822617
3.98 -7.56423653902389
3.99 -7.45217649284664
4 -7.37516720251796
4.01 -7.32542342367179
4.02 -7.29904017526295
4.03 -7.29477692113616
4.04 -7.31400361707413
4.05 -7.36201899876097
4.06 -7.45297714240517
4.07 -7.6376101776528
4.08 -8.44937164325558
4.09 -10.0942868208583
4.1 -11.7894605999658
4.11 -10.457953161577
4.12 -10.2985685437068
4.13 -11.3781603132623
4.14 -13.2291101408542
4.15 -15.7853023409009
4.16 -20.8174672350733
4.17 -23.6287141948723
4.18 -32.6593176098475
4.19 -60.0408488184536
4.2 -110.145540334196
4.21 -133.493685980338
4.22 -213.072284409773
4.23 -413.189954327902
4.24 -inf
4.25 -inf
4.26 -inf
4.27 -inf
4.28 -inf
4.29 -inf
4.3 -inf
4.31 -inf
4.32 -inf
4.33 -inf
4.34 -inf
4.35 -inf
4.36 -inf
4.37 -inf
4.38 -inf
4.39 -inf
4.4 -inf
4.41 -inf
4.42 -inf
4.43 -inf
4.44 -inf
4.45 -inf
4.46 -inf
4.47 -inf
4.48 -inf
4.49 -inf
4.5 -inf
4.51 -inf
4.52 -inf
4.53 -inf
4.54 -inf
4.55 -inf
4.56 -inf
4.57 -inf
4.58 -inf
4.59 -inf
4.6 -inf
4.61 -inf
4.62 -inf
4.63 -inf
4.64 -inf
4.65 -inf
4.66 -inf
4.67 -inf
4.68 -inf
4.69 -inf
4.7 -inf
4.71 -inf
4.72 -inf
4.73 -inf
4.74 -inf
4.75 -inf
4.76 -inf
4.77 -inf
4.78 -inf
4.79 -inf
4.8 -inf
4.81 -inf
4.82 -inf
4.83 -inf
4.84 -inf
4.85 -inf
4.86 -inf
4.87 -inf
4.88 -inf
4.89 -inf
4.9 -inf
4.91 -inf
4.92 -inf
4.93 -inf
4.94 -inf
4.95 -inf
4.96 -inf
4.97 -inf
4.98 -inf
4.99 -inf
5 -inf
5.01 -inf
5.02 -inf
5.03 -inf
5.04 -inf
5.05 -inf
5.06 -inf
5.07 -inf
5.08 -inf
5.09 -inf
5.1 -inf
5.11 -inf
5.12 -inf
5.13 -inf
5.14 -inf
5.15 -inf
5.16 -inf
5.17 -inf
5.18 -inf
5.19 -inf
5.2 -inf
5.21 -inf
5.22 -inf
5.23 -inf
5.24 -inf
5.25 -inf
5.26 -inf
5.27 -inf
5.28 -inf
5.29 -inf
5.3 -inf
5.31 -inf
5.32 -inf
5.33 -inf
5.34 -inf
5.35 -inf
5.36 -inf
5.37 -inf
5.38 -inf
5.39 -inf
5.4 -inf
5.41 -inf
5.42 -inf
5.43 -inf
5.44 -inf
5.45 -inf
5.46 -inf
5.47 -inf
5.48 -inf
5.49 -inf
5.5 -inf
5.51 -inf
5.52 -inf
5.53 -inf
5.54 -inf
5.55 -inf
5.56 -inf
5.57 -inf
5.58 -inf
5.59 -inf
5.6 -inf
5.61 -inf
5.62 -inf
5.63 -inf
5.64 -inf
5.65 -inf
5.66 -inf
5.67 -inf
5.68 -inf
5.69 -inf
5.7 -inf
5.71 -inf
5.72 -inf
5.73 -inf
5.74 -inf
5.75 -inf
5.76 -inf
5.77 -inf
5.78 -inf
5.79 -inf
5.8 -inf
5.81 -inf
5.82 -inf
5.83 -inf
5.84 -inf
5.85 -inf
5.86 -inf
5.87 -inf
5.88 -inf
5.89 -inf
5.9 -inf
5.91 -inf
5.92 -inf
5.93 -inf
5.94 -inf
5.95 -inf
5.96 -inf
5.97 -inf
5.98 -inf
5.99 -inf
6 -inf
};
\addlegendentry{$\log\lvert\beta(1, t)\rvert$}
\addplot [black, dashed]
table {%
0 -inf
0.01 -inf
0.02 -inf
0.03 -inf
0.04 -inf
0.05 -inf
0.06 -inf
0.07 -inf
0.08 -inf
0.09 -inf
0.1 -inf
0.11 -inf
0.12 -inf
0.13 -inf
0.14 -inf
0.15 -inf
0.16 -inf
0.17 -inf
0.18 -inf
0.19 -inf
0.2 -inf
0.21 -11.092524032384
0.22 -9.21371016707373
0.23 -7.95816949036738
0.24 -7.03836686872944
0.25 -6.33655472024772
0.26 -5.78821879482635
0.27 -5.35265407505756
0.28 -5.00199677366912
0.29 -4.71620399849603
0.3 -4.48042364636937
0.31 -4.2834651750723
0.32 -4.11682734013629
0.33 -3.97403077275902
0.34 -3.8501319789144
0.35 -3.7413560425767
0.36 -3.64481439480026
0.37 -3.55828763102929
0.38 -3.48005959309345
0.39 -3.40879195752114
0.4 -3.34343041949363
0.41 -3.2831350662633
0.42 -3.22722893033302
0.43 -3.17515999867698
0.44 -3.12647307546111
0.45 -3.08078881796304
0.46 -3.03778798758229
0.47 -2.99719950222941
0.48 -2.95879127553664
0.49 -2.92236311555631
0.5 -2.88774115999119
0.51 -2.85477346970249
0.52 -2.8233265046829
0.53 -2.79328227947657
0.54 -2.76453604708899
0.55 -2.73699439797486
0.56 -2.71057368802314
0.57 -2.68519872955448
0.58 -2.66080169427135
0.59 -2.63732118829624
0.6 -2.61470146791167
0.61 -2.59289177109505
0.62 -2.57184574493657
0.63 -2.55152095291043
0.64 -2.53187844901051
0.65 -2.51288240816041
0.66 -2.49449980421374
0.67 -2.47670012838386
0.68 -2.4594551421682
0.69 -2.44273865982372
0.7 -2.42652635625711
0.71 -2.41079559685279
0.72 -2.39552528630416
0.73 -2.38069573396133
0.74 -2.36628853357969
0.75 -2.35228645566319
0.76 -2.33867335085482
0.77 -2.32543406304405
0.78 -2.31255435104407
0.79 -2.3000208178465
0.8 -2.28782084659273
0.81 -2.27594254251279
0.82 -2.26437468017865
0.83 -2.25310665550022
0.84 -2.24212844196315
0.85 -2.23143055066806
0.86 -2.2210039937831
0.87 -2.21084025106737
0.88 -2.2009312391621
0.89 -2.19126928338067
0.9 -2.18184709175887
0.91 -2.17265773115256
0.92 -2.16369460519339
0.93 -2.15495143393288
0.94 -2.14642223502339
0.95 -2.13810130630001
0.96 -2.12998320964106
0.97 -2.12206275599764
0.98 -2.11433499149296
0.99 -2.10679518450235
1 -2.09943881363323
1.01 -2.09226155653191
1.02 -2.08525927945124
1.03 -2.07842802751888
1.04 -2.07176401565186
1.05 -2.06526362006763
1.06 -2.05892337034667
1.07 -2.05273994200515
1.08 -2.04671014954037
1.09 -2.04083093991442
1.1 -2.0350993864447
1.11 -2.02951268307261
1.12 -2.02406813898384
1.13 -2.01876317355631
1.14 -2.01359531161338
1.15 -2.00856217896202
1.16 -2.00366149819706
1.17 -1.99889108475439
1.18 -1.9942488431971
1.19 -1.98973276371989
1.2 -1.98534091885828
1.21 -1.98107146039018
1.22 -1.97692261641818
1.23 -1.97289268862198
1.24 -1.96898004967128
1.25 -1.96518314078981
1.26 -1.96150046946228
1.27 -1.95793060727651
1.28 -1.95447218789344
1.29 -1.95112390513863
1.3 -1.94788451120886
1.31 -1.94475281498857
1.32 -1.9417276804705
1.33 -1.9388080252762
1.34 -1.93599281927172
1.35 -1.93328108327445
1.36 -1.93067188784769
1.37 -1.92816435217917
1.38 -1.92575764304087
1.39 -1.92345097382696
1.4 -1.92124360366774
1.41 -1.919134836617
1.42 -1.91712402091105
1.43 -1.91521054829757
1.44 -1.91339385343277
1.45 -1.91167341334566
1.46 -1.9100487469683
1.47 -1.90851941473126
1.48 -1.90708501822371
1.49 -1.9057451999177
1.5 -1.90449964295654
1.51 -1.90334807100733
1.52 -1.902290248178
1.53 -1.9013259789993
1.54 -1.90045510847263
1.55 -1.89967752218475
1.56 -1.89899314649058
1.57 -1.89840194876585
1.58 -1.89790393773123
1.59 -1.89749916385045
1.6 -1.89718771980463
1.61 -1.89696974104592
1.62 -1.89684540643354
1.63 -1.89681493895605
1.64 -1.89687860654385
1.65 -1.89703672297661
1.66 -1.89728964889069
1.67 -1.8976377928923
1.68 -1.89808161278285
1.69 -1.89862161690344
1.7 -1.89925836560638
1.71 -1.89999247286258
1.72 -1.90082460801433
1.73 -1.90175549768434
1.74 -1.90278592785303
1.75 -1.90391674611732
1.76 -1.90514886414583
1.77 -1.90648326034706
1.78 -1.90792098276898
1.79 -1.9094631522508
1.8 -1.91111096585013
1.81 -1.9128657005715
1.82 -1.91472871742563
1.83 -1.91670146585256
1.84 -1.91878548854592
1.85 -1.92098242672078
1.86 -1.92329402587321
1.87 -1.92572214208653
1.88 -1.92826874894687
1.89 -1.93093594514011
1.9 -1.93372596281311
1.91 -1.93664117679483
1.92 -1.93968411478846
1.93 -1.94285746866387
1.94 -1.94616410700154
1.95 -1.94960708906574
1.96 -1.953189680417
1.97 -1.95691537041289
1.98 -1.96078789189501
1.99 -1.96481124341955
2 -1.96898971446402
2.01 -1.97332791413662
2.02 -1.97783080403437
2.03 -1.98250373604955
2.04 -1.98735249612204
2.05 -1.99238335519565
2.06 -1.99760312898067
2.07 -2.00301924858912
2.08 -2.00863984474142
2.09 -2.01447384912289
2.1 -2.02053111771302
2.11 -2.02682258271248
2.12 -2.03336044236885
2.13 -2.04015840209377
2.14 -2.04723198674657
2.15 -2.0545989546715
2.16 -2.06227986273506
2.17 -2.07029886638431
2.18 -2.07868490973966
2.19 -2.08747362627652
2.2 -2.09671075307573
2.21 -2.10646027447716
2.22 -2.11684564826067
2.23 -2.12810605477752
2.24 -2.14060469193927
2.25 -2.15476126924818
2.26 -2.17094622047085
2.27 -2.18939671564192
2.28 -2.21018836683017
2.29 -2.23325887038173
2.3 -2.25845907143397
2.31 -2.28560645632849
2.32 -2.31452602870124
2.33 -2.34507419902773
2.34 -2.3771481181478
2.35 -2.41068543793
2.36 -2.4456592604824
2.37 -2.48207165883903
2.38 -2.51994768117157
2.39 -2.5593306536456
2.4 -2.60027894902379
2.41 -2.64286408913594
2.42 -2.68716996413579
2.43 -2.73329297184627
2.44 -2.78134293943755
2.45 -2.83144475419744
2.46 -2.88374068913494
2.47 -2.93839346324841
2.48 -2.99559013146812
2.49 -3.05554696373301
2.5 -3.11851555599264
2.51 -3.1847905301019
2.52 -3.25471934129675
2.53 -3.32871494680023
2.54 -3.40727243859842
2.55 -3.49099127708857
2.56 -3.58060559974992
2.57 -3.67702643124233
2.58 -3.78140187045176
2.59 -3.89520518363552
2.6 -4.02036750890736
2.61 -4.15948383974942
2.62 -4.31613873062516
2.63 -4.49537608737144
2.64 -4.70451445533661
2.65 -4.95502032729792
2.66 -5.26741385523051
2.67 -5.68612691324061
2.68 -6.34283683165005
2.69 -8.28304894622715
2.7 -6.77921662725563
2.71 -6.0222148591061
2.72 -5.6203623017257
2.73 -5.35003496627683
2.74 -5.14817535308456
2.75 -4.98806798848065
2.76 -4.85596396135088
2.77 -4.74390289860744
2.78 -4.64688761894492
2.79 -4.56158640031377
2.8 -4.48567056586118
2.81 -4.41744926401521
2.82 -4.35565571475895
2.83 -4.29931586378326
2.84 -4.24766427562697
2.85 -4.20008831533191
2.86 -4.15608992531638
2.87 -4.1152587212834
2.88 -4.07725259293416
2.89 -4.0417834171784
2.9 -4.00860633948268
2.91 -3.97751159978776
2.92 -3.94831820816688
2.93 -3.92086898825349
2.94 -3.89502664756015
2.95 -3.87067062936145
2.96 -3.84769456679571
2.97 -3.82600420621395
2.98 -3.8055156999207
2.99 -3.78615419244439
3 -3.76785264208772
3.01 -3.75055083259467
3.02 -3.73419453960035
3.03 -3.71873482399025
3.04 -3.70412743001076
3.05 -3.69033227039016
3.06 -3.6773129841722
3.07 -3.66503655566873
3.08 -3.65347298507784
3.09 -3.64259500301972
3.1 -3.63237782261115
3.11 -3.62279892380565
3.12 -3.61383786562595
3.13 -3.60547612265198
3.14 -3.59769694273553
3.15 -3.59048522341768
3.16 -3.58382740494829
3.17 -3.57771137816525
3.18 -3.57212640579733
3.19 -3.56706305602038
3.2 -3.56251314733053
3.21 -3.55846970400805
3.22 -3.55492692163873
3.23 -3.55188014234136
3.24 -3.54932583952698
3.25 -3.54726161219305
3.26 -3.54568618894011
3.27 -3.54459944209647
3.28 -3.54400241255635
3.29 -3.54389734618738
3.3 -3.54428774295761
3.31 -3.54517842028475
3.32 -3.54657559254283
3.33 -3.5484869692003
3.34 -3.55092187474896
3.35 -3.55389139446714
3.36 -3.55740855122127
3.37 -3.5614885200608
3.38 -3.56614888946893
3.39 -3.57140998105549
3.4 -3.5772952436198
3.41 -3.5838317435235
3.42 -3.5910507822795
3.43 -3.59898868606436
3.44 -3.60768783389844
3.45 -3.61719802802051
3.46 -3.62757837483577
3.47 -3.63889996756821
3.48 -3.65124991784042
3.49 -3.6647379018141
3.5 -3.6795083240793
3.51 -3.69577497697021
3.52 -3.71394767131458
3.53 -3.73478755883583
3.54 -3.75942304336945
3.55 -3.78916712923825
3.56 -3.82524668212181
3.57 -3.86859796265966
3.58 -3.91980506377732
3.59 -3.9791640204603
3.6 -4.04681124352325
3.61 -4.12285973397238
3.62 -4.20751277006384
3.63 -4.30115043538451
3.64 -4.4044006869536
3.65 -4.51821406537464
3.66 -4.64396403939423
3.67 -4.78359846310252
3.68 -4.93987493102778
3.69 -5.11667803109344
3.7 -5.31945862507356
3.71 -5.55612510538959
3.72 -5.83932019934691
3.73 -6.1926308196168
3.74 -6.67043420273397
3.75 -7.45417982711181
3.76 -10.4610242125949
3.77 -7.47953330483115
3.78 -6.8694265618059
3.79 -6.52778513145596
3.8 -6.29654528542268
3.81 -6.12507428774781
3.82 -5.99079957148053
3.83 -5.88181094892164
3.84 -5.79112886593986
3.85 -5.71435091471761
3.86 -5.64854099706315
3.87 -5.5916524893043
3.88 -5.54220706589018
3.89 -5.4991058997466
3.9 -5.46151361682737
3.91 -5.42878421117327
3.92 -5.40041218973905
3.93 -5.37599947917629
3.94 -5.35523255261406
3.95 -5.3378664415498
3.96 -5.32371358709293
3.97 -5.31263627155072
3.98 -5.30454188164126
3.99 -5.29938061963174
4 -5.29714558531601
4.01 -5.29787547600886
4.02 -5.30166059629441
4.03 -5.30865364464256
4.04 -5.31908844498197
4.05 -5.3333147794589
4.06 -5.35188627785098
4.07 -5.37593453031379
4.08 -5.40769116937954
4.09 -5.45067773852958
4.1 -5.50935138176575
4.11 -5.58847222818927
4.12 -5.69261500103916
4.13 -5.8260682738498
4.14 -5.99305618696122
4.15 -6.19824720049993
4.16 -6.44782841126934
4.17 -6.75192951115273
4.18 -7.13056480060444
4.19 -7.63135705733803
4.2 -8.40941386491187
4.21 -11.5877444883545
4.22 -8.72928407226642
4.23 -8.12135352910686
4.24 -7.81710416872198
4.25 -7.63430687635955
4.26 -7.51850614112492
4.27 -7.44743448011538
4.28 -7.41282232170094
4.29 -7.41529573422552
4.3 -7.46063961358127
4.31 -7.55597635127664
4.32 -7.7072993406959
4.33 -7.91866238313171
4.34 -8.19232297166436
4.35 -8.52922831915971
4.36 -8.92958566725424
4.37 -9.39346276995287
4.38 -9.92145954120997
4.39 -10.5155835152238
4.4 -11.1806504836072
4.41 -11.9270496115355
4.42 -12.7774789342319
4.43 -13.7882967239746
4.44 -15.1576421258577
4.45 -17.8993586051725
4.46 -16.2530092684
4.47 -16.3046539068084
4.48 -16.6328838122656
4.49 -17.0906646567612
4.5 -17.6283542996728
4.51 -18.2225968356786
4.52 -18.860239707858
4.53 -19.532966087361
4.54 -20.2350710893
4.55 -20.9623995342837
4.56 -21.7117837710813
4.57 -22.4807219494293
4.58 -23.2671821003584
4.59 -24.0694766086125
4.6 -24.8861782929955
4.61 -25.7160622369435
4.62 -26.5580641920107
4.63 -27.4112500155039
4.64 -28.2747926759527
4.65 -29.1479545870678
4.66 -30.030073782241
4.67 -30.9205529157765
4.68 -31.818850384313
4.69 -32.7244730658691
4.7 -33.6369703123151
4.71 -34.5559289268241
4.72 -35.4809689253219
4.73 -36.4117399292937
4.74 -37.3479180724736
4.75 -38.2892033299012
4.76 -39.2353171972412
4.77 -40.1860006629582
4.78 -41.141012427198
4.79 -42.1001273299379
4.8 -43.0631349577857
4.81 -44.0298384041845
4.82 -45.0000531620626
4.83 -45.9736061314098
4.84 -46.9503347270402
4.85 -47.9300860740736
4.86 -48.9127162805253
4.87 -49.8980897779346
4.88 -50.8860787222372
4.89 -51.8765624481557
4.9 -52.8694269712795
4.91 -53.864564532761
4.92 -54.8618731822002
4.93 -55.8612563948349
4.94 -56.8626227196245
4.95 -57.8658854552172
4.96 -58.8709623511379
4.97 -59.8777753318359
4.98 -60.8862502414912
4.99 -61.8963166077077
5 -62.9079074224194
5.01 -63.9209589385112
5.02 -64.9354104808083
5.03 -65.9512042702251
5.04 -66.9682852599825
5.05 -67.9866009829086
5.06 -69.0061014089338
5.07 -70.0267388119711
5.08 -71.0484676454523
5.09 -72.0712444258514
5.1 -73.0950276235925
5.11 -74.1197775607868
5.12 -75.1454563152972
5.13 -76.1720276306682
5.14 -77.1994568314986
5.15 -78.2277107438723
5.16 -79.2567576204907
5.17 -80.2865670701817
5.18 -81.3171099914852
5.19 -82.3483585100391
5.2 -83.380285919512
5.21 -84.4128666258477
5.22 -85.4460760946038
5.23 -86.4798908011862
5.24 -87.5142881837922
5.25 -88.5492465988913
5.26 -89.5847452790848
5.27 -90.6207642931967
5.28 -91.6572845084573
5.29 -92.694287554655
5.3 -93.7317557901348
5.31 -94.7696722695353
5.32 -95.8080207131611
5.33 -96.8467854778938
5.34 -97.885951529554
5.35 -98.9255044166301
5.36 -99.9654302452953
5.37 -101.005715655642
5.38 -102.046347799063
5.39 -103.087314316719
5.4 -104.128603319029
5.41 -105.170203366133
5.42 -106.212103449264
5.43 -107.254292972999
5.44 -108.296761738318
5.45 -109.339499926452
5.46 -110.382498083464
5.47 -111.425747105526
5.48 -112.469238224861
5.49 -113.512962996319
5.5 -114.55691328454
5.51 -115.601081251693
5.52 -116.645459345751
5.53 -117.690040289272
5.54 -118.734817068674
5.55 -119.779782923974
5.56 -120.824931338955
5.57 -121.870256031768
5.58 -122.915750945922
5.59 -123.961410241656
5.6 -125.007228287677
5.61 -126.053199653237
5.62 -127.099319100546
5.63 -128.145581577489
5.64 -129.191982210652
5.65 -130.238516298625
5.66 -131.285179305585
5.67 -132.331966855131
5.68 -133.378874724375
5.69 -134.425898838267
5.7 -135.473035264143
5.71 -136.520280206495
5.72 -137.567630001949
5.73 -138.615081114433
5.74 -139.662630130545
5.75 -140.710273755094
5.76 -141.758008806818
5.77 -142.805832214269
5.78 -143.853741011856
5.79 -144.901732336036
5.8 -145.949803421662
5.81 -146.997951598456
5.82 -148.046174287628
5.83 -149.09446899862
5.84 -150.142833325966
5.85 -151.191264946284
5.86 -152.239761615366
5.87 -153.28832116539
5.88 -154.336941502223
5.89 -155.385620602839
5.9 -156.434356512817
5.91 -157.483147343945
5.92 -158.531991271904
5.93 -159.580886534038
5.94 -160.629831427212
5.95 -161.678824305736
5.96 -162.72786357938
5.97 -163.776947711446
5.98 -164.82607521692
5.99 -165.875244660689
6 -166.924454655817
};
\addlegendentry{$\log\lvert\beta(0, t)\rvert=\log\lvert\alpha(0, t)\rvert$}
\addplot [black, dash pattern=on 1pt off 3pt on 3pt off 3pt]
table {%
0 -inf
0.01 -inf
0.02 -inf
0.03 -inf
0.04 -inf
0.05 -inf
0.06 -inf
0.07 -inf
0.08 -inf
0.09 -inf
0.1 -inf
0.11 -inf
0.12 -inf
0.13 -inf
0.14 -inf
0.15 -inf
0.16 -inf
0.17 -inf
0.18 -inf
0.19 -inf
0.2 -inf
0.21 -inf
0.22 -inf
0.23 -inf
0.24 -inf
0.25 -inf
0.26 -inf
0.27 -inf
0.28 -inf
0.29 -inf
0.3 -inf
0.31 -inf
0.32 -inf
0.33 -inf
0.34 -inf
0.35 -inf
0.36 -inf
0.37 -inf
0.38 -inf
0.39 -inf
0.4 -inf
0.41 -16.7560238315997
0.42 -14.3221786101675
0.43 -12.548239374227
0.44 -11.1464865783929
0.45 -9.99902699638158
0.46 -9.04114116251347
0.47 -8.23175888782917
0.48 -7.54239935226533
0.49 -6.95207951776066
0.5 -6.44464421719239
0.51 -6.00723483246682
0.52 -5.62935241007823
0.53 -5.30225642446064
0.54 -5.01856470602344
0.55 -4.77197936285351
0.56 -4.55709396084545
0.57 -4.36925387306732
0.58 -4.20445141499733
0.59 -4.05924345507912
0.6 -3.93068327947113
0.61 -3.81626141358061
0.62 -3.71385224797264
0.63 -3.62166485761683
0.64 -3.53819743196448
0.65 -3.4621953255803
0.66 -3.39261298546628
0.67 -3.32858001599441
0.68 -3.26937150911162
0.69 -3.21438258340843
0.7 -3.16310690190087
0.71 -3.11511880856081
0.72 -3.07005864857938
0.73 -3.02762081279266
0.74 -2.98754406079112
0.75 -2.94960371625248
0.76 -2.91360537980379
0.77 -2.87937986007175
0.78 -2.8467790765566
0.79 -2.81567273528167
0.8 -2.78594561849919
0.81 -2.75749536298807
0.82 -2.73023062826544
0.83 -2.70406957726051
0.84 -2.67893860863796
0.85 -2.65477129291572
0.86 -2.6315074745825
0.87 -2.60909251022771
0.88 -2.58747661876568
0.89 -2.56661432456957
0.9 -2.54646397803786
0.91 -2.52698734103434
0.92 -2.50814922695327
0.93 -2.48991718699926
0.94 -2.47226123574432
0.95 -2.45515361020911
0.96 -2.43856855767522
0.97 -2.42248214821601
0.98 -2.40687210857253
0.99 -2.39171767452608
1 -2.3769994593527
1.01 -2.36269933630492
1.02 -2.34880033336568
1.03 -2.33528653877057
1.04 -2.32214301600479
1.05 -2.30935572715921
1.06 -2.29691146368012
1.07 -2.2847977836747
1.08 -2.27300295504317
1.09 -2.26151590380125
1.1 -2.25032616703627
1.11 -2.23942385000862
1.12 -2.22879958696936
1.13 -2.2184445053155
1.14 -2.20835019274902
1.15 -2.19850866714374
1.16 -2.18891234885783
1.17 -2.17955403525877
1.18 -2.17042687725321
1.19 -2.16152435763653
1.2 -2.15284027109651
1.21 -2.14436870572301
1.22 -2.1361040258907
1.23 -2.12804085639541
1.24 -2.12017406773674
1.25 -2.11249876245003
1.26 -2.10501026240044
1.27 -2.09770409696006
1.28 -2.09057599199661
1.29 -2.08362185960903
1.3 -2.07683778855108
1.31 -2.07022003528961
1.32 -2.06376501564882
1.33 -2.0574692969964
1.34 -2.05132959093104
1.35 -2.04534274643452
1.36 -2.0395057434547
1.37 -2.03381568688867
1.38 -2.02826980093772
1.39 -2.02286542380836
1.4 -2.01760000273561
1.41 -2.0124710893069
1.42 -2.00747633506629
1.43 -2.00261348738092
1.44 -1.99788038555252
1.45 -1.99327495715836
1.46 -1.9887952146074
1.47 -1.9844392518983
1.48 -1.98020524156688
1.49 -1.9760914318121
1.5 -1.97209614378967
1.51 -1.9682177690641
1.52 -1.96445476720988
1.53 -1.9608056635539
1.54 -1.95726904705125
1.55 -1.9538435682876
1.56 -1.95052793760144
1.57 -1.94732092332049
1.58 -1.94422135010644
1.59 -1.94122809740323
1.6 -1.93834009798395
1.61 -1.93555633659229
1.62 -1.93287584867439
1.63 -1.93029771919763
1.64 -1.927821081553
1.65 -1.92544511653813
1.66 -1.9231690514182
1.67 -1.92099215906239
1.68 -1.91891375715363
1.69 -1.91693320746983
1.7 -1.91504991523485
1.71 -1.91326332853781
1.72 -1.91157293781962
1.73 -1.90997827542561
1.74 -1.90847891522372
1.75 -1.90707447228756
1.76 -1.90576460264425
1.77 -1.90454900308676
1.78 -1.90342741105119
1.79 -1.90239960455915
1.8 -1.90146540222615
1.81 -1.90062466333666
1.82 -1.8998772879873
1.83 -1.89922321729933
1.84 -1.89866243370238
1.85 -1.89819496129144
1.86 -1.89782086625935
1.87 -1.89754025740771
1.88 -1.89735328673921
1.89 -1.89726015013479
1.9 -1.89726108811983
1.91 -1.89735638672352
1.92 -1.89754637843651
1.93 -1.89783144327231
1.94 -1.89821200993855
1.95 -1.89868855712499
1.96 -1.89926161491584
1.97 -1.89993176633475
1.98 -1.90069964903208
1.99 -1.90156595712459
2 -1.90253144319941
2.01 -1.90359692049511
2.02 -1.90476326527431
2.03 -1.90603141940394
2.04 -1.90740239316115
2.05 -1.90887726828506
2.06 -1.91045720129686
2.07 -1.91214342711382
2.08 -1.91393726298577
2.09 -1.91584011278634
2.1 -1.91785347169557
2.11 -1.91997893131556
2.12 -1.92221818526612
2.13 -1.92457303531463
2.14 -1.92704539810182
2.15 -1.92963731253436
2.16 -1.93235094792631
2.17 -1.93518861298424
2.18 -1.93815276574612
2.19 -1.94124602460307
2.2 -1.94447118055487
2.21 -1.94783121087762
2.22 -1.95132929441512
2.23 -1.95496882874594
2.24 -1.95875344952939
2.25 -1.9626870523961
2.26 -1.96677381782924
2.27 -1.97101823958351
2.28 -1.97542515731951
2.29 -1.97999979430124
2.3 -1.98474780122803
2.31 -1.98967530757174
2.32 -1.99478898219708
2.33 -2.00009610560666
2.34 -2.00560465694997
2.35 -2.01132342009421
2.36 -2.01726211478858
2.37 -2.02343156165074
2.38 -2.02984389410934
2.39 -2.03651283813705
2.4 -2.04345409557245
2.41 -2.05068590355673
2.42 -2.05823001230746
2.43 -2.06611374208741
2.44 -2.07437418901838
2.45 -2.08306534899398
2.46 -2.09226745680232
2.47 -2.10209573503788
2.48 -2.11270442103711
2.49 -2.12428258925883
2.5 -2.13704089412594
2.51 -2.15119165382783
2.52 -2.16692704651697
2.53 -2.18440062203536
2.54 -2.20371594273277
2.55 -2.22492382482217
2.56 -2.24802737835857
2.57 -2.27299253204885
2.58 -2.29976118996037
2.59 -2.328264436861
2.6 -2.35843395000506
2.61 -2.39021064970004
2.62 -2.42355038654285
2.63 -2.45842699818493
2.64 -2.49483335265304
2.65 -2.53278107213488
2.66 -2.5722995706167
2.67 -2.61343490965316
2.68 -2.6562488317683
2.69 -2.70081820330147
2.7 -2.74723500309498
2.71 -2.79560693330119
2.72 -2.84605870028672
2.73 -2.89873401154266
2.74 -2.95379835378347
2.75 -3.01144265535823
2.76 -3.07188799296142
2.77 -3.13539158207725
2.78 -3.20225440034114
2.79 -3.27283094600444
2.8 -3.34754184989551
2.81 -3.42689036703584
2.82 -3.51148420177702
2.83 -3.60206458941749
2.84 -3.69954492839897
2.85 -3.80506187268635
2.86 -3.92004412682925
2.87 -4.04631087754855
2.88 -4.18622643019368
2.89 -4.34296620179037
2.9 -4.52100697826427
2.91 -4.72708712314163
2.92 -4.97224057216944
2.93 -5.27666381624427
2.94 -5.68392394358915
2.95 -6.3204142421049
2.96 -8.1021743163727
2.97 -6.81717731574376
2.98 -6.03161135977795
2.99 -5.61955858983855
3 -5.34382899897109
3.01 -5.13891905548797
3.02 -4.97718569164701
3.03 -4.84437153484932
3.04 -4.73218764704755
3.05 -4.63541341783544
3.06 -4.55056535183221
3.07 -4.47521406596172
3.08 -4.40760416129728
3.09 -4.34642926039886
3.1 -4.29069237941285
3.11 -4.23961590748808
3.12 -4.1925817189224
3.13 -4.14909024762008
3.14 -4.10873184825453
3.15 -4.07116631991439
3.16 -4.03610797358719
3.17 -4.00331454176656
3.18 -3.97257880152262
3.19 -3.9437221483564
3.2 -3.91658959633405
3.21 -3.89104583766866
3.22 -3.8669721009917
3.23 -3.84426362004165
3.24 -3.82282757479458
3.25 -3.80258140249246
3.26 -3.78345140135168
3.27 -3.76537156809713
3.28 -3.74828262396034
3.29 -3.73213119382191
3.3 -3.71686911074003
3.31 -3.70245282386365
3.32 -3.68884289215596
3.33 -3.67600354979077
3.34 -3.66390233177526
3.35 -3.65250975047779
3.36 -3.6417990154293
3.37 -3.63174579012203
3.38 -3.62232798062281
3.39 -3.61352555170738
3.4 -3.6053203669506
3.41 -3.59769604980792
3.42 -3.59063786322415
3.43 -3.5841326057246
3.44 -3.57816852230051
3.45 -3.57273522870653
3.46 -3.56782364805483
3.47 -3.56342595882798
3.48 -3.55953555364725
3.49 -3.55614700833419
3.5 -3.55325606099517
3.51 -3.55085960104967
3.52 -3.54895566831942
3.53 -3.54754346250436
3.54 -3.54662336360224
3.55 -3.54619696409174
3.56 -3.54626711400805
3.57 -3.54683798041272
3.58 -3.54791512322039
3.59 -3.54950558992767
3.6 -3.55161803254062
3.61 -3.55426285098788
3.62 -3.55745236863885
3.63 -3.56120104738105
3.64 -3.56552575230046
3.65 -3.57044607976996
3.66 -3.57598476840093
3.67 -3.5821682211537
3.68 -3.58902718148757
3.69 -3.59659763232917
3.7 -3.60492203848643
3.71 -3.61405119971001
3.72 -3.62404753788564
3.73 -3.63499182759546
3.74 -3.64699635107416
3.75 -3.66022634887294
3.76 -3.67492746469974
3.77 -3.691451400932
3.78 -3.71026900980058
3.79 -3.73196239926024
3.8 -3.75719475308323
3.81 -3.78666475703938
3.82 -3.82105766906086
3.83 -3.86100525764146
3.84 -3.90706298586842
3.85 -3.95970744149743
3.86 -4.01935246220435
3.87 -4.08637990713292
3.88 -4.16118065084491
3.89 -4.24420243201952
3.9 -4.33600267930531
3.91 -4.43730572430765
3.92 -4.54906539675637
3.93 -4.67253772387952
3.94 -4.80937691550447
3.95 -4.96178401467208
3.96 -5.13276695281392
3.97 -5.32662852425908
3.98 -5.54993182549574
3.99 -5.81355401075989
4 -6.13761022255129
4.01 -6.56586160201926
4.02 -7.22632740823018
4.03 -9.07275968165372
4.04 -7.73998299343854
4.05 -6.98203307766551
4.06 -6.59179974251914
4.07 -6.33642094967682
4.08 -6.15109650833
4.09 -6.00841758023679
4.1 -5.89425499647032
4.11 -5.80042022369927
4.12 -5.72178594216387
4.13 -5.65496003013158
4.14 -5.59760326266181
4.15 -5.54804884880503
4.16 -5.50507724328535
4.17 -5.46777672994739
4.18 -5.43545408208969
4.19 -5.40757576930303
4.2 -5.38372848341568
4.21 -5.36359228684808
4.22 -5.34692229138128
4.23 -5.33353634876619
4.24 -5.32330724510658
4.25 -5.31615861985404
4.26 -5.31206456856218
4.27 -5.31105465219289
4.28 -5.31322965538735
4.29 -5.31879663260223
4.3 -5.32812913734294
4.31 -5.34184700692857
4.32 -5.36089457693621
4.33 -5.38658780840969
4.34 -5.42060771145583
4.35 -5.46493735288888
4.36 -5.52176172993457
4.37 -5.59336338173826
4.38 -5.68204889978924
4.39 -5.79013727823126
4.4 -5.92003853540948
4.41 -6.07445906533355
4.42 -6.25680065216832
4.43 -6.47189974954265
4.44 -6.72746160813856
4.45 -7.03715419256663
4.46 -7.42849906476845
4.47 -7.96916179313782
4.48 -8.91611714030928
4.49 -10.0596667930304
4.5 -8.58264305530148
4.51 -8.1050892014651
4.52 -7.85188793293147
4.53 -7.70594289954191
4.54 -7.62769244329842
4.55 -7.60018522837174
4.56 -7.61541859885508
4.57 -7.66944656100322
4.58 -7.76023917518428
4.59 -7.88664456192635
4.6 -8.04787000195751
4.61 -8.24323015137914
4.62 -8.47203696869642
4.63 -8.73356394034151
4.64 -9.02704753453672
4.65 -9.351705831132
4.66 -9.70676403116237
4.67 -10.0914821403298
4.68 -10.505183331386
4.69 -10.9472834155115
4.7 -11.4173231969312
4.71 -11.9150067572917
4.72 -12.4402503903742
4.73 -12.993249605801
4.74 -13.5745764160722
4.75 -14.1853281754978
4.76 -14.8273673411973
4.77 -15.5037302302829
4.78 -16.2193728313574
4.79 -16.9826546883743
4.8 -17.8086579854371
4.81 -18.7279877476832
4.82 -19.8175631126716
4.83 -21.3949882081999
4.84 -22.5484698359224
4.85 -21.8941722464754
4.86 -22.0067176163762
4.87 -22.325779004253
4.88 -22.747466015263
4.89 -23.2330864549648
4.9 -23.7637948729533
4.91 -24.3288950320316
4.92 -24.9216687711356
4.93 -25.5375771422763
4.94 -26.1733786637429
4.95 -26.8266549623876
4.96 -27.4955367654908
4.97 -28.1785364410501
4.98 -28.8744408424872
4.99 -29.5822400658168
5 -30.3010785203684
5.01 -31.0302203717159
5.02 -31.7690245329066
5.03 -32.5169261725401
5.04 -33.2734227774649
5.05 -34.0380634664484
5.06 -34.8104406684754
5.07 -35.5901835504671
5.08 -36.3769527593865
5.09 -37.1704361658746
5.1 -37.9703453809458
5.11 -38.7764128765309
5.12 -39.5883895829171
5.13 -40.4060428666947
5.14 -41.2291548152095
5.15 -42.0575207701215
5.16 -42.8909480651153
5.17 -43.7292549322372
5.18 -44.5722695485476
5.19 -45.4198292003429
5.2 -46.2717795465407
5.21 -47.1279739662189
5.22 -47.9882729779972
5.23 -48.8525437210863
5.24 -49.7206594895563
5.25 -50.5924993127585
5.26 -51.4679475759661
5.27 -52.3468936762184
5.28 -53.2292317091109
5.29 -54.1148601828994
5.3 -55.003681756806
5.31 -55.8956030008451
5.32 -56.7905341748565
5.33 -57.6883890247323
5.34 -58.5890845940889
5.35 -59.4925410498511
5.36 -60.398681520404
5.37 -61.3074319451318
5.38 -62.2187209342951
5.39 -63.1324796383224
5.4 -64.048641625692
5.41 -64.9671427686697
5.42 -65.8879211362476
5.43 -66.8109168936945
5.44 -67.7360722081917
5.45 -68.6633311600767
5.46 -69.5926396592659
5.47 -70.5239453664688
5.48 -71.4571976188396
5.49 -72.3923473597474
5.5 -73.3293470723728
5.51 -74.2681507168654
5.52 -75.2087136708195
5.53 -76.1509926728445
5.54 -77.0949457690273
5.55 -78.0405322620991
5.56 -78.9877126631329
5.57 -79.9364486456149
5.58 -80.8867030017422
5.59 -81.8384396008112
5.6 -82.7916233495723
5.61 -83.7462201544339
5.62 -84.7021968854089
5.63 -85.6595213417033
5.64 -86.6181622188537
5.65 -87.578089077328
5.66 -88.5392723125081
5.67 -89.5016831259789
5.68 -90.4652934980548
5.69 -91.4300761614767
5.7 -92.3960045762189
5.71 -93.3630529053477
5.72 -94.3311959918791
5.73 -95.300409336583
5.74 -96.2706690766888
5.75 -97.241951965446
5.76 -98.2142353524984
5.77 -99.1874971650332
5.78 -100.161715889667
5.79 -101.136870555032
5.8 -102.112940715036
5.81 -103.089906432756
5.82 -104.06774826494
5.83 -105.046447247094
5.84 -106.025984879118
5.85 -107.006343111478
5.86 -107.987504331876
5.87 -108.969451352411
5.88 -109.952167397198
5.89 -110.935636090433
5.9 -111.919841444878
5.91 -112.904767850756
5.92 -113.890400065028
5.93 -114.876723201053
5.94 -115.863722718591
5.95 -116.85138441416
5.96 -117.839694411712
5.97 -118.828639153633
5.98 -119.818205392035
5.99 -120.808380180348
6 -121.799150865181
};
\addlegendentry{$\log\lvert\alpha(1, t)\rvert$}
\end{axis}

\end{tikzpicture}

%% file: pics/simXparticle.tex
\begin{tikzpicture}

\begin{axis}[
width=\plottimewidth,
height=\plottimeheight,
legend cell align={left},
legend style={fill opacity=0.8, draw opacity=1, text opacity=1, draw=white!80!black},
tick align=inside,
tick pos=both,
x grid style={white!69.0196078431373!black},
xlabel={$t$},
xmajorgrids,
xmin=0, xmax=6,
xtick style={color=black},
y grid style={white!69.0196078431373!black},
ymajorgrids,
ymin=-0.1, ymax=0.6,
ytick style={color=black},
ytick={-0.1,0,0.1,0.2,0.3,0.4,0.5,0.6},
yticklabels={-0.1,0.0,0.1,0.2,0.3,0.4,0.5,0.6}
]
\addplot [black]
table {%
0 0.5
0.0499999523162842 0.499308347702026
0.100000023841858 0.498069286346436
0.149999976158142 0.496296167373657
0.200000047683716 0.494005799293518
0.25 0.491213083267212
0.299999952316284 0.4879310131073
0.350000023841858 0.484171390533447
0.399999976158142 0.479944944381714
0.450000047683716 0.475261211395264
0.5 0.470102906227112
0.539999961853027 0.465537071228027
0.579999923706055 0.460451602935791
0.620000004768372 0.454785585403442
0.660000085830688 0.448544144630432
0.700000047683716 0.441758036613464
0.740000009536743 0.434461355209351
0.779999971389771 0.426685810089111
0.819999933242798 0.418460845947266
0.860000014305115 0.409814238548279
0.910000085830688 0.398452162742615
0.960000038146973 0.38651978969574
1.00999999046326 0.374063014984131
1.05999994277954 0.361125588417053
1.11000001430511 0.347750067710876
1.16999995708466 0.331179141998291
1.23000001907349 0.314105749130249
1.29999995231628 0.293642044067383
1.37999999523163 0.269665718078613
1.47000002861023 0.242126822471619
1.5900000333786 0.204846739768982
1.78999996185303 0.142662525177002
1.87000000476837 0.118284583091736
1.94000005722046 0.0974357128143311
2 0.0800683498382568
2.04999995231628 0.0660967826843262
2.08999991416931 0.0554928779602051
2.13000011444092 0.0455037355422974
2.16000008583069 0.0384615659713745
2.19000005722046 0.0318293571472168
2.22000002861023 0.0256085395812988
2.25 0.0197813510894775
2.28999996185303 0.0125750303268433
2.32999992370605 0.0059514045715332
2.36999988555908 -0.000142693519592285
2.41000008583069 -0.00575411319732666
2.46000003814697 -0.0122324228286743
2.5 -0.0169479846954346
2.53999996185303 -0.0211842060089111
2.57999992370605 -0.0248680114746094
2.61999988555908 -0.0279453992843628
2.66000008583069 -0.0304011106491089
2.70000004768372 -0.0322470664978027
2.74000000953674 -0.0335085391998291
2.77999997138977 -0.0342191457748413
2.8199999332428 -0.0344183444976807
2.85999989509583 -0.0341500043869019
2.90000009536743 -0.0334641933441162
2.95000004768372 -0.0321155786514282
3.00999999046326 -0.0299614667892456
3.07999992370605 -0.0269467830657959
3.17000007629395 -0.0225554704666138
3.39000010490417 -0.0114951133728027
3.45000004768372 -0.00938355922698975
3.57999992370605 -0.00546836853027344
3.6800000667572 -0.00292015075683594
3.77999997138977 -0.000844597816467285
3.85999989509583 0.000371217727661133
3.96000003814697 0.00129401683807373
4.03999996185303 0.00142848491668701
4.15000009536743 0.00110697746276855
4.46999979019165 -5.00679016113281e-05
6 0
};
\addlegendentry{$X_\mathrm{p}(t)=y(1,t)$}
\addplot [black, dashed]
table {%
0 0.5
0.129999995231628 0.499287128448486
0.279999971389771 0.498401284217834
0.309999942779541 0.497646689414978
0.340000033378601 0.4964280128479
0.370000004768372 0.49474573135376
0.399999976158142 0.492621660232544
0.430000066757202 0.490078926086426
0.460000038146973 0.487138032913208
0.490000009536743 0.483817458152771
0.529999971389771 0.478828310966492
0.569999933242798 0.473230481147766
0.610000014305115 0.467057824134827
0.649999976158142 0.460341930389404
0.690000057220459 0.453112721443176
0.730000019073486 0.445398092269897
0.769999980926514 0.437224626541138
0.809999942779541 0.42861795425415
0.860000014305115 0.417287230491638
0.910000085830688 0.405363321304321
0.960000038146973 0.392890095710754
1.00999999046326 0.379910111427307
1.05999994277954 0.366464495658875
1.12000000476837 0.34977126121521
1.17999994754791 0.332532525062561
1.24000000953674 0.31481409072876
1.30999994277954 0.293622493743896
1.38999998569489 0.26884138584137
1.48000001907349 0.240423917770386
1.61000001430511 0.198785543441772
1.79999995231628 0.137909531593323
1.88999998569489 0.109610795974731
1.96000003814697 0.088067889213562
2.02999997138977 0.0670318603515625
2.08999991416931 0.0494798421859741
2.14000010490417 0.0354253053665161
2.1800000667572 0.024682879447937
2.22000002861023 0.0144593715667725
2.25999999046326 0.00480711460113525
2.28999996185303 -0.0019681453704834
2.3199999332428 -0.00825059413909912
2.34999990463257 -0.0139920711517334
2.38000011444092 -0.0191788673400879
2.41000008583069 -0.0238159894943237
2.44000005722046 -0.0279419422149658
2.47000002861023 -0.0315290689468384
2.5 -0.0345655679702759
2.52999997138977 -0.0370662212371826
2.55999994277954 -0.0390521287918091
2.58999991416931 -0.0405490398406982
2.61999988555908 -0.0415873527526855
2.65000009536743 -0.0422040224075317
2.69000005722046 -0.0424654483795166
2.73000001907349 -0.042232871055603
2.77999997138977 -0.0414227247238159
2.82999992370605 -0.0401459932327271
2.89000010490417 -0.0381019115447998
2.95000004768372 -0.0355924367904663
3.01999998092651 -0.0321879386901855
3.09999990463257 -0.0278100967407227
3.20000004768372 -0.0218400955200195
3.4300000667572 -0.00783979892730713
3.54999995231628 -0.0016028881072998
3.59999990463257 0.00051724910736084
3.64000010490417 0.00175929069519043
3.6800000667572 0.00256037712097168
3.73000001907349 0.0030134916305542
3.78999996185303 0.00300121307373047
3.88000011444092 0.0024104118347168
4.05999994277954 0.0006866455078125
4.15999984741211 -5.59091567993164e-05
4.25 -0.000107049942016602
4.57999992370605 -0
6 -0
};
\addlegendentry{$y(0, t)$}
\addplot [semithick, black, mark=*, mark size=1.25, mark options={solid}, forget plot]
table {%
4.76 0
};
\draw (axis cs:4.76,0.021) node[
  anchor=base west,
  text=black,
  rotate=0.0
]{$T_1$};
\end{axis}

\end{tikzpicture}

%% file: pics/simXparticle_log.tex
\begin{tikzpicture}

\begin{axis}[
width=\plottimewidth,
height=\plottimeheight,
legend cell align={left},
legend style={fill opacity=0.8, draw opacity=1, text opacity=1, at={(0.03,0.03)}, anchor=south west, draw=white!80!black},
tick align=inside,
tick pos=both,
x grid style={white!69.0196078431373!black},
xlabel={$t$},
xmajorgrids,
xmin=0, xmax=6,
xtick style={color=black},
y grid style={white!69.0196078431373!black},
ymajorgrids,
ymin=-30, ymax=0,
ytick style={color=black},
ytick={-30,-20,-10,0},
yticklabels={-30,-20,-10,0}
]
\addplot [black]
table {%
0 -0.301029995663981
0.01 -0.301112639926184
0.02 -0.301213573323474
0.03 -0.301333537773887
0.04 -0.301472736566344
0.05 -0.301631199687336
0.06 -0.30180888590039
0.07 -0.302005721980408
0.08 -0.302221620508
0.09 -0.302456488774472
0.1 -0.302710233525773
0.11 -0.302982763577326
0.12 -0.303273991274356
0.13 -0.30358383329613
0.14 -0.303912211075615
0.15 -0.304259050988985
0.16 -0.304624284405761
0.17 -0.305007847654599
0.18 -0.305409681938544
0.19 -0.305829733221192
0.2 -0.306267952097253
0.21 -0.306724293656353
0.22 -0.307198717362605
0.23 -0.307691186977227
0.24 -0.30820167052751
0.25 -0.308730140293766
0.26 -0.30927657277493
0.27 -0.309840948607773
0.28 -0.310423252438581
0.29 -0.311023472763807
0.3 -0.311641601761693
0.31 -0.312277635132926
0.32 -0.312931571960716
0.33 -0.313603414593583
0.34 -0.314293168549486
0.35 -0.315000842437785
0.36 -0.315726447895114
0.37 -0.316469999531753
0.38 -0.317231514885887
0.39 -0.318011014383958
0.4 -0.318808521305919
0.41 -0.319624062089665
0.42 -0.320457668880366
0.43 -0.321309389468802
0.44 -0.322179312940968
0.45 -0.323067618801235
0.46 -0.323974650399649
0.47 -0.324901001753223
0.48 -0.325847596019403
0.49 -0.32681573027978
0.5 -0.327807067685499
0.51 -0.328823571982028
0.52 -0.329867395072008
0.53 -0.330940739606409
0.54 -0.332045722305263
0.55 -0.333184259985051
0.56 -0.334357991831089
0.57 -0.335568241824549
0.58 -0.336816017263958
0.59 -0.338102034472173
0.6 -0.339426761166251
0.61 -0.340790465827265
0.62 -0.34219326665918
0.63 -0.343635175389563
0.64 -0.345116133567873
0.65 -0.346636040826007
0.66 -0.348194775713599
0.67 -0.349792210299967
0.68 -0.351428219903709
0.69 -0.353102689229242
0.7 -0.354815515985182
0.71 -0.356566612819702
0.72 -0.358355908183333
0.73 -0.360183346543205
0.74 -0.362048888230056
0.75 -0.363952509096731
0.76 -0.365894200096817
0.77 -0.367873966846209
0.78 -0.369891829201712
0.79 -0.371947820873459
0.8 -0.374041989078024
0.81 -0.376174394233659
0.82 -0.378345109696347
0.83 -0.380554221534165
0.84 -0.382801828337051
0.85 -0.385088041059103
0.86 -0.387412982890788
0.87 -0.389776789158705
0.88 -0.392179607250879
0.89 -0.394621596565846
0.9 -0.397102928484045
0.91 -0.399623786360273
0.92 -0.402184365536143
0.93 -0.404784873371681
0.94 -0.407425529295323
0.95 -0.410106564871739
0.96 -0.412828223886994
0.97 -0.415590762450683
0.98 -0.418394449114785
0.99 -0.421239565009013
1 -0.424126403992586
1.01 -0.427055272822369
1.02 -0.430026491337417
1.03 -0.433040392660019
1.04 -0.436097323413402
1.05 -0.439197643956311
1.06 -0.442341728634741
1.07 -0.445529966051147
1.08 -0.448762759351522
1.09 -0.452040526530792
1.1 -0.455363700757003
1.11 -0.45873273071489
1.12 -0.462148080969402
1.13 -0.465610232349883
1.14 -0.46911968235563
1.15 -0.472676945583614
1.16 -0.476282554179239
1.17 -0.479937058311063
1.18 -0.483641026670471
1.19 -0.487395046997387
1.2 -0.49119972663317
1.21 -0.495055693101939
1.22 -0.498963594721645
1.23 -0.502924101246309
1.24 -0.50693790454095
1.25 -0.511005719290824
1.26 -0.515128283746708
1.27 -0.519306360508094
1.28 -0.523540737346265
1.29 -0.527832228069396
1.3 -0.532181673431932
1.31 -0.536589942090681
1.32 -0.541057931610224
1.33 -0.545586569520413
1.34 -0.550176814428944
1.35 -0.554829657192182
1.36 -0.559546122147669
1.37 -0.564327268411959
1.38 -0.569174191247718
1.39 -0.574088023504292
1.4 -0.579069937136272
1.41 -0.584121144804906
1.42 -0.589242901567575
1.43 -0.594436506660952
1.44 -0.599703305383867
1.45 -0.6050446910864
1.46 -0.610462107272179
1.47 -0.615957049821459
1.48 -0.621531069343113
1.49 -0.627185773664339
1.5 -0.63292283046759
1.51 -0.638743970085001
1.52 -0.644650988461451
1.53 -0.650645750298296
1.54 -0.65673019239086
1.55 -0.662906327173847
1.56 -0.669176246490087
1.57 -0.675542125599345
1.58 -0.682006227445426
1.59 -0.688570907201412
1.6 -0.69523861711467
1.61 -0.702011911675248
1.62 -0.708893453133467
1.63 -0.715886017394929
1.64 -0.722992500323858
1.65 -0.730215924488647
1.66 -0.737559446386812
1.67 -0.745026364190206
1.68 -0.752620126055456
1.69 -0.760344339049132
1.7 -0.768202778742239
1.71 -0.776199399534316
1.72 -0.784338345773761
1.73 -0.792623963748079
1.74 -0.801060814625717
1.75 -0.809653688440008
1.76 -0.81840761921568
1.77 -0.827327901349497
1.78 -0.836420107369012
1.79 -0.845690107207239
1.8 -0.855144089146422
1.81 -0.864788582601141
1.82 -0.874630482929757
1.83 -0.884677078483677
1.84 -0.894936080126142
1.85 -0.905415653475733
1.86 -0.916124454154261
1.87 -0.927071666343089
1.88 -0.938267044974688
1.89 -0.949720961905031
1.9 -0.961444456423292
1.91 -0.973449290452847
1.92 -0.985748008774016
1.93 -0.998354004545475
1.94 -1.0112815903132
1.95 -1.02454607459516
1.96 -1.03816384415034
1.97 -1.05215245276753
1.98 -1.06653072139403
1.99 -1.08131887444149
2 -1.09653885970402
2.01 -1.11221647197167
2.02 -1.12841110597794
2.03 -1.14506653941159
2.04 -1.16219618134697
2.05 -1.17982003076221
2.06 -1.19796140855382
2.07 -1.21664626137547
2.08 -1.23590300016598
2.09 -1.25576256040899
2.1 -1.27625861791426
2.11 -1.29742796261683
2.12 -1.31931106174545
2.13 -1.34195285901647
2.14 -1.36540386461547
2.15 -1.389721593742
2.16 -1.41497241189574
2.17 -1.44123384750904
2.18 -1.46859744430473
2.19 -1.49717225644219
2.2 -1.52708915008616
2.21 -1.55850617786538
2.22 -1.59161545388815
2.23 -1.62665220262399
2.24 -1.66390703359856
2.25 -1.70374309892066
2.26 -1.74662080385149
2.27 -1.79313452360142
2.28 -1.84406908284833
2.29 -1.90049022816673
2.3 -1.96389686715261
2.31 -2.03649353073469
2.32 -2.12171837893998
2.33 -2.22538179797376
2.34 -2.35852503274823
2.35 -2.54653039374301
2.36 -2.87495669592653
2.37 -3.84565685654626
2.38 -2.79895642435782
2.39 -2.52212060243218
2.4 -2.35721304857322
2.41 -2.24001787674721
2.42 -2.14930708990716
2.43 -2.07502537901181
2.44 -2.01268460043552
2.45 -1.95918551572963
2.46 -1.91248625400852
2.47 -1.87118385987483
2.48 -1.83427694201682
2.49 -1.80102847451835
2.5 -1.77088231839586
2.51 -1.74340988038682
2.52 -1.71827453721841
2.53 -1.69520707487452
2.54 -1.6739883123886
2.55 -1.65443663437263
2.56 -1.63639899851122
2.57 -1.61974445211148
2.58 -1.60435946673819
2.59 -1.59014457778542
2.6 -1.57701194419002
2.61 -1.56488354365067
2.62 -1.55368979889879
2.63 -1.543368493708
2.64 -1.53386388505104
2.65 -1.52512595204653
2.66 -1.51710974550173
2.67 -1.50977481653638
2.68 -1.50308471140439
2.69 -1.49700652428356
2.7 -1.49151050206734
2.71 -1.48656969616782
2.72 -1.48215965671218
2.73 -1.47825816466697
2.74 -1.47484499754367
2.75 -1.47190172449726
2.76 -1.46941152684143
2.77 -1.46735904025772
2.78 -1.46573021524483
2.79 -1.464512192609
2.8 -1.46369319099834
2.81 -1.46326240357553
2.82 -1.46320990077933
2.83 -1.46352653525489
2.84 -1.46420384261836
2.85 -1.46523392798999
2.86 -1.4666093260741
2.87 -1.46832282667688
2.88 -1.47036727050476
2.89 -1.47273533869226
2.9 -1.47541937532955
2.91 -1.47841128641282
2.92 -1.48170254730992
2.93 -1.48528432751876
2.94 -1.48914771486673
2.95 -1.49328400061772
2.96 -1.49768497802102
2.97 -1.50234321039311
2.98 -1.50725223745853
2.99 -1.51240670511036
3 -1.51780241908967
3.01 -1.52343633419486
3.02 -1.52930649645185
3.03 -1.53541195676897
3.04 -1.54175267239266
3.05 -1.54832940857674
3.06 -1.55514364860011
3.07 -1.56219751647509
3.08 -1.56949371378512
3.09 -1.57703547014498
3.1 -1.58482650565478
3.11 -1.5928710032066
3.12 -1.60117358838145
3.13 -1.60973931476056
3.14 -1.61857365263898
3.15 -1.6276824792834
3.16 -1.63707206897361
3.17 -1.64674908108286
3.18 -1.65672054437956
3.19 -1.66699383557224
3.2 -1.67757664987576
3.21 -1.68847696105703
3.22 -1.69970296802976
3.23 -1.71126302462026
3.24 -1.72316554864657
3.25 -1.73541890605544
3.26 -1.74803126605633
3.27 -1.76101042622427
3.28 -1.77436362287951
3.29 -1.78809742883467
3.3 -1.80221844095005
3.31 -1.81674706482377
3.32 -1.83176316988639
3.33 -1.84702531289595
3.34 -1.86246935690995
3.35 -1.87802016759285
3.36 -1.89358670367518
3.37 -1.90907092836943
3.38 -1.9243875163475
3.39 -1.93948621018508
3.4 -1.95436462770345
3.41 -1.96906495809447
3.42 -1.98365878941546
3.43 -1.99822974959508
3.44 -2.01286112004973
3.45 -2.02763048845054
3.46 -2.04261057272518
3.47 -2.05787578132984
3.48 -2.07352064645216
3.49 -2.08985119235477
3.5 -2.10681942785368
3.51 -2.12428104281094
3.52 -2.14224235585096
3.53 -2.1607249652161
3.54 -2.17975760821073
3.55 -2.19937422376894
3.56 -2.21961353052952
3.57 -2.24051918192155
3.58 -2.26214023365109
3.59 -2.28453184437722
3.6 -2.30775619880085
3.61 -2.33188367697616
3.62 -2.35699431895798
3.63 -2.38317965953622
3.64 -2.41054503937517
3.65 -2.43921254165734
3.66 -2.4693247642258
3.67 -2.50104972664145
3.68 -2.5345873462509
3.69 -2.57017812518838
3.7 -2.6081150191465
3.71 -2.64875999196839
3.72 -2.69256763895803
3.73 -2.74011974508835
3.74 -2.79217729130632
3.75 -2.84976155034328
3.76 -2.91428667684049
3.77 -2.98779046674031
3.78 -3.07336940479245
3.79 -3.17608745671987
3.8 -3.30515742412833
3.81 -3.48038851129885
3.82 -3.75977642193607
3.83 -4.60466804181663
3.84 -3.93616697095538
3.85 -3.60568229508545
3.86 -3.43043156106934
3.87 -3.3154542121149
3.88 -3.22723074301809
3.89 -3.15316581550195
3.9 -3.08908570470172
3.91 -3.03493406556195
3.92 -2.99051135697559
3.93 -2.95493226021549
3.94 -2.92762488125325
3.95 -2.90567646016441
3.96 -2.88806792406802
3.97 -2.87411294348762
3.98 -2.8632884606014
3.99 -2.85517624822521
4 -2.84943305873101
4.01 -2.84577344098388
4.02 -2.8439592415992
4.03 -2.84379291612553
4.04 -2.8451136313804
4.05 -2.84779817852338
4.06 -2.8517781082707
4.07 -2.8570919031162
4.08 -2.86402343294943
4.09 -2.8744438626718
4.1 -2.88546247542404
4.11 -2.89745830928865
4.12 -2.91053815454675
4.13 -2.9247135803935
4.14 -2.93971045434856
4.15 -2.95583859139238
4.16 -2.97312221623853
4.17 -2.99161147213696
4.18 -3.01136670506347
4.19 -3.03246033164974
4.2 -3.0549783572543
4.21 -3.07902246690585
4.22 -3.10471275584012
4.23 -3.13219127633453
4.24 -3.16162665210853
4.25 -3.19322011773029
4.26 -3.22721349015962
4.27 -3.26389976540743
4.28 -3.30363720202008
4.29 -3.34686785633513
4.3 -3.39414168497822
4.31 -3.44614796730426
4.32 -3.50375774717594
4.33 -3.5680854997878
4.34 -3.64058746948344
4.35 -3.72323313732784
4.36 -3.81882900536468
4.37 -3.93168339299962
4.38 -4.06913056313092
4.39 -4.24566111830042
4.4 -4.49780531924893
4.41 -4.98549068480574
4.42 -5.12272461412178
4.43 -4.65828050431187
4.44 -4.47961117856479
4.45 -4.38355613600314
4.46 -4.3288716897179
4.47 -4.2999338711839
4.48 -4.28920280062375
4.49 -4.2925692578997
4.5 -4.30767510316988
4.51 -4.33318128190469
4.52 -4.36839772358187
4.53 -4.41306687447447
4.54 -4.46721775146727
4.55 -4.53105671983472
4.56 -4.60488243267376
4.57 -4.68902099366599
4.58 -4.78377985329159
4.59 -4.88941870945584
4.6 -5.00613478046808
4.61 -5.13405922960793
4.62 -5.27326150602025
4.63 -5.42375881080176
4.64 -5.5855285758822
4.65 -5.75852256977159
4.66 -5.94268191509057
4.67 -6.13795289232939
4.68 -6.34430394436495
4.69 -6.56174486356928
4.7 -6.79034986505754
4.71 -7.03028733906437
4.72 -7.28186092335979
4.73 -7.5455699269525
4.74 -7.82220377597015
4.75 -8.11299907341674
4.76 -8.41991950593748
4.77 -8.74619843246654
4.78 -9.09751281904137
4.79 -9.48494910937033
4.8 -9.93456052577791
4.81 -10.5366810757248
4.82 -11.6237456337839
4.83 -10.9738832089825
4.84 -10.9836271049025
4.85 -11.1059129425058
4.86 -11.2794645695132
4.87 -11.4836514771982
4.88 -11.7089167917248
4.89 -11.9500029173783
4.9 -12.2036787262294
4.91 -12.4677964047163
4.92 -12.7408414536128
4.93 -13.0216954289608
4.94 -13.3095010335051
4.95 -13.6035807115041
4.96 -13.9033851076646
4.97 -14.2084591072037
4.98 -14.5184186937751
4.99 -14.8329347170367
5 -15.1517212171786
5.01 -15.4745268395135
5.02 -15.801128396172
5.03 -16.1313259522874
5.04 -16.4649390157011
5.05 -16.8018035394642
5.06 -17.1417695325036
5.07 -17.4846991319146
5.08 -17.8304650302914
5.09 -18.1789491794427
5.1 -18.5300417116889
5.11 -18.8836400342394
5.12 -19.2396480625897
5.13 -19.597975566599
5.14 -19.9585376086762
5.15 -20.3212540578617
5.16 -20.6860491669184
5.17 -21.0528512021035
5.18 -21.4215921172813
5.19 -21.7922072655971
5.2 -22.1646351431558
5.21 -22.5388171601295
5.22 -22.9146974354972
5.23 -23.2922226122513
5.24 -23.6713416904135
5.25 -24.0520058756214
5.26 -24.4341684413868
5.27 -24.8177846034087
5.28 -25.2028114045564
5.29 -25.589207609335
5.3 -25.9769336068049
5.31 -26.3659513210653
5.32 -26.7562241285287
5.33 -27.1477167813081
5.34 -27.5403953361245
5.35 -27.9342270882142
5.36 -28.3291805097729
5.37 -28.7252251925309
5.38 -29.1223317940973
5.39 -29.5204719877485
5.4 -29.9196184153763
5.41 -30.3197446433339
5.42 -30.7208251209504
5.43 -31.1228351415052
5.44 -31.525750805472
5.45 -31.9295489858642
5.46 -32.3342072955266
5.47 -32.7397040562324
5.48 -33.1460182694597
5.49 -33.5531295887284
5.5 -33.9610182933939
5.51 -34.3696652637984
5.52 -34.7790519576913
5.53 -35.189160387836
5.54 -35.5999731007287
5.55 -36.0114731563587
5.56 -36.4236441089463
5.57 -36.8364699886
5.58 -37.2499352838362
5.59 -37.6640249249124
5.6 -38.078724267926
5.61 -38.4940190796346
5.62 -38.9098955229582
5.63 -39.3263401431239
5.64 -39.7433398544187
5.65 -40.1608819275177
5.66 -40.5789539773559
5.67 -40.9975439515159
5.68 -41.4166401191037
5.69 -41.8362310600877
5.7 -42.2563056550778
5.71 -42.6768530755217
5.72 -43.097862774297
5.73 -43.5193244766822
5.74 -43.9412281716839
5.75 -44.3635641037076
5.76 -44.7863227645516
5.77 -45.2094948857119
5.78 -45.633071430981
5.79 -46.0570435893292
5.8 -46.4814027680533
5.81 -46.9061405861822
5.82 -47.3312488681271
5.83 -47.7567196375651
5.84 -48.1825451115464
5.85 -48.6087176948152
5.86 -49.0352299743343
5.87 -49.4620747140063
5.88 -49.8892448495803
5.89 -50.3167334837391
5.9 -50.7445338813576
5.91 -51.1726394649253
5.92 -51.6010438101266
5.93 -52.0297406415727
5.94 -52.4587238286781
5.95 -52.8879873816758
5.96 -53.317525447767
5.97 -53.7473323073984
5.98 -54.1774023706623
5.99 -54.6077301738162
6 -55.0383103759144
};
\addlegendentry{$\log\lvert X_\mathrm{p}(t)\lvert=\log\lvert y(1,t)\lvert$}
\addplot [black, dashed]
table {%
0 -0.301029995663981
0.01 -0.301093782341582
0.02 -0.301152401600521
0.03 -0.301207340647329
0.04 -0.301259374631588
0.05 -0.301309002113482
0.06 -0.301356579283043
0.07 -0.301402377260102
0.08 -0.301446611243139
0.09 -0.301489457080852
0.1 -0.301531061462458
0.11 -0.301571548545407
0.12 -0.301611024456582
0.13 -0.301649580449614
0.14 -0.301687295174595
0.15 -0.3017242363391
0.16 -0.301760461937685
0.17 -0.301796021166588
0.18 -0.3018309551022
0.19 -0.301865297198116
0.2 -0.301899073639236
0.21 -0.301932459439331
0.22 -0.301966301349076
0.23 -0.302002873884391
0.24 -0.302046390267302
0.25 -0.302102881675934
0.26 -0.302179462821623
0.27 -0.302283328180317
0.28 -0.302420870268253
0.29 -0.302597150136859
0.3 -0.302815748949175
0.31 -0.303078899943447
0.32 -0.303387763435189
0.33 -0.303742732082715
0.34 -0.304143699335741
0.35 -0.304590264666661
0.36 -0.305081875210782
0.37 -0.305617915574749
0.38 -0.306197760452579
0.39 -0.306820802919881
0.4 -0.307486467989156
0.41 -0.308194217834742
0.42 -0.308943552640617
0.43 -0.309734009347819
0.44 -0.310565159525146
0.45 -0.311436606965321
0.46 -0.312347985265038
0.47 -0.313298955477856
0.48 -0.314289203866266
0.49 -0.315318439775318
0.5 -0.316386393667388
0.51 -0.317492815370077
0.52 -0.318637472584117
0.53 -0.319820149676028
0.54 -0.32104064674977
0.55 -0.322298778963309
0.56 -0.323594376037611
0.57 -0.324927281899674
0.58 -0.326297354406144
0.59 -0.327704465105732
0.6 -0.329148499012593
0.61 -0.330629354375706
0.62 -0.332146942439196
0.63 -0.333701187194899
0.64 -0.335292025131835
0.65 -0.336919404988227
0.66 -0.338583287511372
0.67 -0.340283645229606
0.68 -0.342020462239321
0.69 -0.34379373400891
0.7 -0.345603467200535
0.71 -0.347449679509974
0.72 -0.349332399524366
0.73 -0.351251666597383
0.74 -0.353207530741261
0.75 -0.355200052535027
0.76 -0.35722930304831
0.77 -0.35929536378011
0.78 -0.361398326612003
0.79 -0.363538293775261
0.8 -0.365715377831475
0.81 -0.36792970166629
0.82 -0.37018139849594
0.83 -0.372470611886313
0.84 -0.374797495784309
0.85 -0.377162214561339
0.86 -0.379564943068819
0.87 -0.382005866705561
0.88 -0.384485181497021
0.89 -0.387003094186377
0.9 -0.389559822337459
0.91 -0.392155594449587
0.92 -0.394790650084395
0.93 -0.397465240004778
0.94 -0.40017962632609
0.95 -0.402934082679808
0.96 -0.405728894389851
0.97 -0.408564358661835
0.98 -0.411440784785513
0.99 -0.414358494350759
1 -0.417317821477409
1.01 -0.420319113059376
1.02 -0.42336272902345
1.03 -0.426449042603243
1.04 -0.429578440628799
1.05 -0.432751323832375
1.06 -0.435968107171025
1.07 -0.43922922016657
1.08 -0.442535107263661
1.09 -0.445886228206651
1.1 -0.449283058436043
1.11 -0.452726089505359
1.12 -0.456215829519315
1.13 -0.459752803594246
1.14 -0.463337554341804
1.15 -0.466970642377015
1.16 -0.470652646851839
1.17 -0.474384166015492
1.18 -0.478165817802831
1.19 -0.481998240452221
1.2 -0.485882093154383
1.21 -0.489818056733842
1.22 -0.493806834364677
1.23 -0.497849152322417
1.24 -0.501945760774038
1.25 -0.506097434608166
1.26 -0.510304974307714
1.27 -0.514569206867359
1.28 -0.51889098675843
1.29 -0.523271196943947
1.3 -0.527710749946764
1.31 -0.532210588973969
1.32 -0.536771689100932
1.33 -0.541395058518641
1.34 -0.546081739848215
1.35 -0.550832811526808
1.36 -0.555649389269412
1.37 -0.560532627611402
1.38 -0.565483721537068
1.39 -0.570503908199764
1.4 -0.575594468739746
1.41 -0.580756730206264
1.42 -0.585992067590994
1.43 -0.591301905980459
1.44 -0.596687722835748
1.45 -0.602151050408484
1.46 -0.607693478302797
1.47 -0.613316656193829
1.48 -0.619022296714272
1.49 -0.624812178521369
1.5 -0.630688149557971
1.51 -0.636652130522404
1.52 -0.642706118563289
1.53 -0.648852191216894
1.54 -0.655092510606258
1.55 -0.661429327923157
1.56 -0.667864988215965
1.57 -0.674401935508743
1.58 -0.681042718279367
1.59 -0.687789995327308
1.6 -0.694646542064773
1.61 -0.701615257268437
1.62 -0.70869917033288
1.63 -0.715901449071262
1.64 -0.723225408113701
1.65 -0.730674517959426
1.66 -0.738252414745062
1.67 -0.745962910798595
1.68 -0.753810006056648
1.69 -0.761797900431965
1.7 -0.769931007228539
1.71 -0.778213967713876
1.72 -0.786651666971711
1.73 -0.79524925117439
1.74 -0.804012146432456
1.75 -0.812946079400181
1.76 -0.822057099840334
1.77 -0.831351605380149
1.78 -0.840836368723787
1.79 -0.850518567625853
1.8 -0.86040581797658
1.81 -0.870506210403862
1.82 -0.880828350862126
1.83 -0.891381405755491
1.84 -0.902175152235832
1.85 -0.913220034429264
1.86 -0.924527226482576
1.87 -0.936108703491739
1.88 -0.947977321588107
1.89 -0.960146908729837
1.9 -0.972632368099588
1.91 -0.985449796482189
1.92 -0.998616620651278
1.93 -1.01215175574668
1.94 -1.02607579109927
1.95 -1.04041121143515
1.96 -1.05518266599378
1.97 -1.07041730783628
1.98 -1.08614524993238
1.99 -1.10240026040204
2 -1.11922115216742
2.01 -1.13665742066913
2.02 -1.15482918072907
2.03 -1.17371895498386
2.04 -1.19337603218918
2.05 -1.21385985884674
2.06 -1.23523842643121
2.07 -1.25758902410021
2.08 -1.28099969456895
2.09 -1.30557125076691
2.1 -1.3314199474691
2.11 -1.35868101163667
2.12 -1.38751334259775
2.13 -1.41810583512555
2.14 -1.45068598684969
2.15 -1.48553177524715
2.16 -1.52298831664509
2.17 -1.56349171545368
2.18 -1.60760409794049
2.19 -1.65606675601543
2.2 -1.70988401127845
2.21 -1.77046205598088
2.22 -1.83985180598919
2.23 -1.92120647436778
2.24 -2.01974019467461
2.25 -2.14505188064259
2.26 -2.31811637017067
2.27 -2.60227711939148
2.28 -3.62127203284888
2.29 -2.70593739146186
2.3 -2.38506399644564
2.31 -2.20655298134804
2.32 -2.08351465048972
2.33 -1.99030634971876
2.34 -1.9157928928637
2.35 -1.8541185721245
2.36 -1.80181971920627
2.37 -1.75667609944301
2.38 -1.71717830445499
2.39 -1.68225248131615
2.4 -1.65110562710819
2.41 -1.62313136964571
2.42 -1.5978399909273
2.43 -1.57474259417231
2.44 -1.55374431169165
2.45 -1.53462785749992
2.46 -1.51719849031578
2.47 -1.50128937230124
2.48 -1.48675798009352
2.49 -1.47348208811383
2.5 -1.46135629813922
2.51 -1.45028919187782
2.52 -1.44020102272842
2.53 -1.43102184400715
2.54 -1.42268998297639
2.55 -1.41515078690026
2.56 -1.40835558250737
2.57 -1.40226080232101
2.58 -1.39682724031883
2.59 -1.39201940545197
2.6 -1.38780494436772
2.61 -1.3841541022376
2.62 -1.38103917174507
2.63 -1.37843376081847
2.64 -1.37631174137348
2.65 -1.37464611993805
2.66 -1.37340838400843
2.67 -1.37256872916935
2.68 -1.37209710545586
2.69 -1.37196464299621
2.7 -1.3721449480911
2.71 -1.37261494055136
2.72 -1.37335515130332
2.73 -1.37434958009779
2.74 -1.3755852859373
2.75 -1.37705187258931
2.76 -1.3787409825442
2.77 -1.38064585960283
2.78 -1.38276100004088
2.79 -1.3850818885551
2.8 -1.38760480470135
2.81 -1.39032668329785
2.82 -1.39324501409514
2.83 -1.39635776920719
2.84 -1.39966334995415
2.85 -1.40316054734169
2.86 -1.40684851230115
2.87 -1.41072673313505
2.88 -1.41479501848899
2.89 -1.41905348472714
2.9 -1.42350254692038
2.91 -1.42814291284487
2.92 -1.43297557949975
2.93 -1.43800183173587
2.94 -1.44322324267146
2.95 -1.44864167566543
2.96 -1.45425928771784
2.97 -1.46007853425839
2.98 -1.46610217535686
2.99 -1.4723332834394
3 -1.47877525262297
3.01 -1.48543180979194
3.02 -1.49230702754406
3.03 -1.49940533913287
3.04 -1.50673155553615
3.05 -1.51429088478621
3.06 -1.52208895370974
3.07 -1.53013183224119
3.08 -1.53842606049439
3.09 -1.54697867880073
3.1 -1.55579726094822
3.11 -1.56488995088511
3.12 -1.57426550318322
3.13 -1.58393332759306
3.14 -1.59390353806695
3.15 -1.60418700668323
3.16 -1.61479542298321
3.17 -1.62574135934702
3.18 -1.63703834320925
3.19 -1.64870093719137
3.2 -1.66074482867111
3.21 -1.6731869310357
3.22 -1.68604550007391
3.23 -1.69934027102008
3.24 -1.71309262539106
3.25 -1.72732580351329
3.26 -1.74206519229787
3.27 -1.75733874894172
3.28 -1.77317770501044
3.29 -1.78961798181875
3.3 -1.80670419759579
3.31 -1.8245200817471
3.32 -1.84329383839628
3.33 -1.86287707469174
3.34 -1.88327863357122
3.35 -1.90449762736317
3.36 -1.92651907825624
3.37 -1.94932242189579
3.38 -1.97290188799486
3.39 -1.99729197013151
3.4 -2.02258578957233
3.41 -2.04893849797866
3.42 -2.07655870200111
3.43 -2.10569749846753
3.44 -2.13664331992862
3.45 -2.16972628092946
3.46 -2.20533322633624
3.47 -2.24393633925388
3.48 -2.28614961257135
3.49 -2.33307063784837
3.5 -2.38547950349277
3.51 -2.44431293763751
3.52 -2.51127086885092
3.53 -2.58890493908935
3.54 -2.68122555870085
3.55 -2.79509230472065
3.56 -2.94392484942743
3.57 -3.16063714896614
3.58 -3.57759578818378
3.59 -3.85752524863636
3.6 -3.28633998387817
3.61 -3.06088499496771
3.62 -2.92305983782272
3.63 -2.82664886473377
3.64 -2.75467025097713
3.65 -2.69895648695709
3.66 -2.65494043590223
3.67 -2.6198056744835
3.68 -2.59169081245727
3.69 -2.56929853343317
3.7 -2.55168193262867
3.71 -2.53811748292773
3.72 -2.52802682950727
3.73 -2.52093051087517
3.74 -2.51642343287251
3.75 -2.51416364240589
3.76 -2.51386731818369
3.77 -2.51530511315747
3.78 -2.51829745733931
3.79 -2.52270833862434
3.8 -2.52843808026711
3.81 -2.53541586441779
3.82 -2.5435925806056
3.83 -2.55293444858765
3.84 -2.56341846500768
3.85 -2.5750346791491
3.86 -2.5878452788124
3.87 -2.6025001887321
3.88 -2.61789834011144
3.89 -2.63310957735326
3.9 -2.64789597783137
3.91 -2.66304134576148
3.92 -2.67946697431818
3.93 -2.69784161167422
3.94 -2.71886661310312
3.95 -2.74172145873857
3.96 -2.76640241372169
3.97 -2.79301067676591
3.98 -2.82169844156834
3.99 -2.85266230623803
4 -2.88614689091573
4.01 -2.92245377773847
4.02 -2.9619553076539
4.03 -3.00511434221538
4.04 -3.05251281980169
4.05 -3.10489696052596
4.06 -3.16326672061042
4.07 -3.22903667830228
4.08 -3.3039107172211
4.09 -3.39694390860096
4.1 -3.50590868172464
4.11 -3.63942239876484
4.12 -3.81228870278865
4.13 -4.06071673704673
4.14 -4.52672976264247
4.15 -4.7477736266829
4.16 -4.25273268721474
4.17 -4.07155175547437
4.18 -3.97663151202766
4.19 -3.92408814986472
4.2 -3.89740583872343
4.21 -3.88860948440349
4.22 -3.89337353387435
4.23 -3.90927718045898
4.24 -3.93505800074685
4.25 -3.97026189303161
4.26 -4.01508496300325
4.27 -4.07032640736599
4.28 -4.1373660984874
4.29 -4.21805241326579
4.3 -4.31446570655736
4.31 -4.42864265288609
4.32 -4.56235896813877
4.33 -4.71698767709203
4.34 -4.8933699658437
4.35 -5.09156253014397
4.36 -5.31019373659088
4.37 -5.54490496446256
4.38 -5.78513050041378
4.39 -6.00968351966754
4.4 -6.18788473880866
4.41 -6.29801021258705
4.42 -6.34833434876159
4.43 -6.36645927190908
4.44 -6.37568765235915
4.45 -6.38817821262013
4.46 -6.40864881153674
4.47 -6.43838393340659
4.48 -6.47740837880548
4.49 -6.52543440777398
4.5 -6.58222154937141
4.51 -6.64768512848483
4.52 -6.72190350817551
4.53 -6.8050874817646
4.54 -6.89753844999159
4.55 -6.9996069957579
4.56 -7.11165727871637
4.57 -7.23403978866572
4.58 -7.36707331151794
4.59 -7.51103583395762
4.6 -7.66616342252226
4.61 -7.83265584256236
4.62 -8.01068778866247
4.63 -8.20042501361097
4.64 -8.40204532293406
4.65 -8.61576537381696
4.66 -8.84187564036151
4.67 -9.08078817969434
4.68 -9.3331058387699
4.69 -9.59972927205401
4.7 -9.8820343931397
4.71 -10.1821903207989
4.72 -10.5037840771795
4.73 -10.8532032531113
4.74 -11.2432584038666
4.75 -11.7056209948024
4.76 -12.3656681192654
4.77 -12.9593108842659
4.78 -12.6261434025302
4.79 -12.6668728539344
4.8 -12.8026018789334
4.81 -12.9853889482003
4.82 -13.1977363149604
4.83 -13.4312209700356
4.84 -13.6811084841477
4.85 -13.9444563551787
4.86 -14.2193032207963
4.87 -14.5042745346478
4.88 -14.798372357504
4.89 -15.1008555660511
4.9 -15.4111683355618
4.91 -15.7288963353018
4.92 -16.0537400131068
4.93 -16.3854993276543
4.94 -16.7240670284462
4.95 -17.0694293062384
4.96 -17.4216739760686
4.97 -17.7810077459838
4.98 -18.1477860280494
4.99 -18.522561898387
5 -18.9061666848837
5.01 -19.2998466231259
5.02 -19.7055065842033
5.03 -20.126177054367
5.04 -20.5670020765422
5.05 -21.0376467556594
5.06 -21.5595947609645
5.07 -22.1990965574855
5.08 -23.7184238693936
5.09 -22.9780108854319
5.1 -23.0209425559165
5.11 -23.2025090521061
5.12 -23.4401788221118
5.13 -23.709073007122
5.14 -23.9982675423695
5.15 -24.301965153884
5.16 -24.6167046556382
5.17 -24.9402442836874
5.18 -25.2710420993687
5.19 -25.6079870778794
5.2 -25.9502484070068
5.21 -26.2971856264324
5.22 -26.6482923058597
5.23 -27.0031592657832
5.24 -27.3614497192085
5.25 -27.7228819707861
5.26 -28.0872170672686
5.27 -28.4542497854782
5.28 -28.8238019265681
5.29 -29.1957172393076
5.3 -29.5698575166353
5.31 -29.9460995521169
5.32 -30.3243327366112
5.33 -30.7044571384276
5.34 -31.086381953391
5.35 -31.4700242412825
5.36 -31.8553078864408
5.37 -32.2421627355835
5.38 -32.6305238770705
5.39 -33.0203310340086
5.4 -33.4115280497699
5.41 -33.8040624490452
5.42 -34.1978850611249
5.43 -34.5929496947534
5.44 -34.9892128560076
5.45 -35.3866335022847
5.46 -35.7851728267528
5.47 -36.1847940686525
5.48 -36.5854623456332
5.49 -36.9871445049689
5.5 -37.3898089910199
5.51 -37.7934257267388
5.52 -38.1979660073645
5.53 -38.6034024047312
5.54 -39.0097086808704
5.55 -39.4168597097639
5.56 -39.8248314062781
5.57 -40.2336006614442
5.58 -40.6431452833655
5.59 -41.0534439431282
5.6 -41.4644761251774
5.61 -41.8762220816796
5.62 -42.2886627904697
5.63 -42.7017799162149
5.64 -43.1155557744807
5.65 -43.529973298418
5.66 -43.9450160078286
5.67 -44.3606679803811
5.68 -44.7769138247937
5.69 -45.1937386558017
5.7 -45.6111280707586
5.71 -46.029068127734
5.72 -46.4475453249813
5.73 -46.8665465816633
5.74 -47.2860592197392
5.75 -47.7060709469181
5.76 -48.1265698405994
5.77 -48.5475443327257
5.78 -48.9689831954783
5.79 -49.390875527759
5.8 -49.8132107423982
5.81 -50.2359785540392
5.82 -50.6591689676542
5.83 -51.0827722676469
5.84 -51.5067790075038
5.85 -51.9311799999616
5.86 -52.3559663076491
5.87 -52.7811292341835
5.88 -53.2066603156874
5.89 -53.6325513127009
5.9 -54.0587942024661
5.91 -54.4853811715647
5.92 -54.9123046088804
5.93 -55.339557098877
5.94 -55.7671314151662
5.95 -56.1950205143546
5.96 -56.6232175301507
5.97 -57.0517157677217
5.98 -57.4805086982805
5.99 -57.9095899538991
6 -58.3389533225296
};
\addlegendentry{$\log\lvert y(0,t)\lvert$}
\end{axis}

\end{tikzpicture}

%% file: DRW_130121_Marc210215_depotjournal.bbl
\begin{thebibliography}{33}
\addcontentsline{toc}{section}{References}
\bibitem{ABCR} B. d'Andr\'ea-Novel, F. Boustany, F. Conrad, and B. P. Rao, {\em Feedback stabilization of a hybrid PDE-ODE system: application to an overhead crane}, MCSS Journal {\bf 7} (1994), 1--22. 
\bibitem{AC} B. d'Andr\'ea-Novel, J.-M. Coron, {\em Exponential stabilization of an overhead crane with flexible cable via a back-stepping approach}, Automatica {\bf 36} (2000), no. 4, 587--593.    
\bibitem{VLFC} B. d'Andr\'ea-Novel, J.-M. Coron, {\em Stabilization of an overhead crane with a variable length  flexible cable}, Computational and Applied Mathematics, Vol. 21 No. 1, (2002), 101--134.
\bibitem{dandrea2010acoustic} B. d'Andr\'ea-Novel, B. Fabre, J.-M. Coron, {\em An acoustic model for
automatic control of a slide flute}, Acta Acustica united with Acustica {\bf 96} (2010), no. 4, 713--721.   
\bibitem{DMR}  B. d'Andr\'ea-Novel, I. Moyano, L. Rosier, {\em Finite-time stabilization of an overhead crane with a flexible cable}, 
Math. Control Signals Systems {\bf 31} (2019). no. 2,  Art. 6, 19 pp.
\bibitem{APR} F. Alabau-Boussouira, V. Perrollaz, L. Rosier, {\em Finite-time stabilization of a network of strings}, Mathematical Control and Related Fields
{\bf 5} (2015), No. 4, 721--742. 
\bibitem{AA} H. Anfinsen, O. M. Aamo, {\em Adaptative Control of Hyperbolic PDEs}, Springer, 2019.  
\bibitem{BR} A. Bacciotti, L. Rosier, Liapunov Functions and Stability in Control Theory, 2nd Edition, Springer-Verlag Berlin Heidelberg 2005. 
\bibitem{BB} S. P. Bhat, D. S. Bernstein, {\em Finite-time stability of continuous autonomous systems}, SIAM Journal Control Optim. {\bf 38} (2000), no. 3, 751--766. 
\bibitem{coron} J.-M. Coron. Control and Nonlinearity. Mathematical Surveys and Monographs, Vol. 136, American Mathematical Society, 2007. 
\bibitem{BBeam} J.-M. Coron, B. d'Andr\'ea-Novel, {\em Stabilization of a rotating body-beam without damping}, IEEE Transactions on Automatic Control, 43(5), (1998), 608--618. \
\bibitem{CVKB} J.-M. Coron, R. Vazquez, M. Krstic, G. Bastin, {\em Local exponential $H^2$ stabilization of a $2\times 2$ quasilinear hyperbolic system using backstepping}, SIAM J. Control Optim. {\bf 51} (2013), no. 3, pp 2005--2035. 
\bibitem{daafouz2014nonlinear} J. Daafouz, M. Tucsnak, J. Valein, {\em Nonlinear control of a coupled PDE/ODE system modeling a switched power converter with a transmission line}, Systems $\&$ Control Letters {\bf 70} (2014), 92--99. 
\bibitem{GPR}  M. Gugat, V.  Perrollaz, L. Rosier, {\em Boundary stabilization of quasilinear hyperbolic 
systems of balance laws: exponential decay for small source terms}, J. Evol. Equ. 18  (2018), no. 3, 1471--1500.
\bibitem{haimo} V. T. Haimo, {\em Finite time controllers}, SIAM J. Control Optim., {\bf 24} (1986), No. 4, 760-770. 
Conference on Decision and Control, Phoenix, 1999, 1302--1307.  
\bibitem{LS} G. Leugering, E. J. P. G. Schmidt, {\em On the modelling and stabilization of flows in networks of open canals}, SIAM J. Control Optim. {\bf 41}  (2002), 164--180.
\bibitem{MP} A. Majda, R. Phillips, {\em Disappearing solutions for the dissipative wave equations}, Indiana University Mathematics Journal {\bf 24} 
(1975), no. 12, 1119--1133.  
\bibitem{MRR} P. Martin, L. Rosier, P. Rouchon, {\em Null controllability of one-dimensional parabolic equations by the flatness approach}, SIAM J. Control Optim. {\bf 54} (2016), no. 1, 198--220.
\bibitem{mattioni2020stabilisation} A. Mattioni, Y. Wu, Y. Le Gorrec, H. Zwart, {\em Stabilisation of a rotating beam clamped on a moving inertia with strong dissipation feedback}, Proceedings of the 59th IEEE Conference on Decision and Control, Jeju Island, 2020, 5056--5061.
\bibitem{mifdal} A. Mifdal, {\em Stabilisation uniforme d'un syst\`eme hybride}, C. R. Acad. Sci. Paris, t. 324, S\'erie I 
(1997), 37--42. 
\bibitem{more} J. J. Mor\'e, B. S. Garbow, K. E. Hillstrom, {\em User guide for MINPACK-1}, Argonne National Laboratory, 1980.
\bibitem{MRC} O. Morg\"ul, B. P.  Rao, F. Conrad, {\em On the stabilization of a cable with a tip mass}, 
IEEE Transactions on Automatic Control {\bf 39}  (1994), No. 10,  440--454.   
\bibitem{PR} V. Perrollaz, L. Rosier. {\em Finite-time stabilization of $2\times 2$ hyperbolic systems on tree-shaped networks}, SIAM J. Control Optim.  {\bf 52} (2014), No. 1, 143--163.
\bibitem{polyakov2019sliding} A. Polyakov. {\em Sliding mode control design using canonical homogeneous norm}, Int. J. Robust and Nonlinear Control, {\bf 29} (2019), No. 3, 682--701.
\bibitem{PEB} A. Polyakov, D. Efimov, B. Brogliato. {\em Consistent discretization of finite-time and fixed-time stable systems}, SIAM J. Control Optim. {\bf 57} (2019), No. 1, 78--103.
\bibitem{Sainsaulieu} L. Sainsaulieu. Calcul scientifique. Masson, 1996. 
\bibitem{Sepulchre} R. Sepulchre, M. Jankovi\'{c}, P. Kokotovi\'{c}. Constructive Nonlinear Control. Springer-Verlag, (1997).
\bibitem{SLX} Y. Shang, D. Liu, G. Xu, {\em Super-stability and the spectrum of one-dimensional wave equations on general feedback controlled 
networks}, IMA J. Math. Control and Inform. {\bf 31} (2014), 73--90.  
\bibitem{trinh2017design} N.-T. Trinh, V. Andrieu, C.-Z. Xu, {\em Design of integral controllers for nonlinear systems governed by scalar hyperbolic partial differential equations}, IEEE Transactions on Automatic Control {\bf 62} (2017), no. 9, 4527--4536.   
\bibitem{VK} R. Vazquez, M. Krstic, {\em Marcum $Q$-functions and explicit kernels for stabilization of $2\times 2$ linear hyperbolic systems with constant coefficients}, Systems $\&$ Control Letters {\bf 68}  (2014), 33--42.  
\end{thebibliography}
